\documentclass{article}

\usepackage[utf8]{inputenc} 
\usepackage[T1]{fontenc} 
\usepackage{mathpazo} 
\usepackage{comment}
\usepackage[dvipsnames]{xcolor}
\usepackage[final]{pdfpages}
\usepackage{amsmath}
\usepackage{amssymb}
\usepackage{mathtools}
\usepackage{amsthm}
\usepackage{thmtools}
\usepackage{thm-restate}
\usepackage{newtxmath}
\usepackage{stmaryrd} 
\usepackage[mathscr]{euscript}
\let\euscr\mathscr \let\mathscr\relax
\usepackage[scr]{rsfso}
\usepackage{comment}
\usepackage[colorlinks=true, allcolors=blue]{hyperref}
\usepackage{cleveref}
\usepackage[notintoc]{nomencl}
\usepackage{enumitem}
\usepackage{anyfontsize}
\counterwithout{footnote}{section}
\usepackage{float}
\usepackage{silence}
\usepackage{etoolbox}
\usepackage{tikz}
\usetikzlibrary{arrows.meta,positioning}
\usetikzlibrary{shapes.geometric}

\numberwithin{equation}{section}
\theoremstyle{plain}
\newtheorem{proposition}{Proposition}[section]
\newtheorem{lemma}[proposition]{Lemma}
\newtheorem{theorem}[proposition]{Theorem}
\newtheorem{corollary}[proposition]{Corollary}

\theoremstyle{definition}
\newtheorem{definition}[proposition]{Definition}

\theoremstyle{remark}
\newtheorem{remark}[proposition]{Remark}
\AtEndEnvironment{remark}{\hfill$\triangle$}

\newtheorem{example}[proposition]{Example}
\AtEndEnvironment{example}{\hfill$\triangledown$}
\newtheorem{conjecture}[proposition]{Conjecture}
\renewcommand\thmcontinues[1]{continued}
\newcommand{\olsi}[1]{\,\overline{\!{#1}}} 
\newcommand{\closedb}{\olsi{B}}
\newcommand{\openb}{\mathring{B}}
\newcommand{\diam}{\mathop{\mathrm{diam}}\nolimits}

\newcommand{\val}{\mathop{\mathrm{val}}\nolimits}
\newcommand{\len}{\mathop{\mathrm{len}}\nolimits}
\newcommand{\supp}{\mathop{\mathrm{supp}}\nolimits}
\newcommand{\pack}{\mathop{\mathrm{pack}}\nolimits}
\newcommand{\Dya}{\mathbb{Y}}
\newcommand{\powerset}{\euscr{P}}
\newcommand{\N}{\mathbb{N}}

\newcommand{\R}{\mathbb{R}}

\newcommand{\kbits}{{\{0,1\}}^k}
\newcommand{\canopy}{\left\llbracket T\right\rrbracket}
\newcommand{\Dyacanopy}{\left\llbracket \Dya\right\rrbracket}
\newcommand{\abs}[1]{\left\lvert #1 \right\rvert}
\newcommand{\cyl}[1]{\llbracket #1 \rrbracket} 
\newcommand{\seqAp}{\langle A_n,p_n\mid n\in \mathbb{N}\rangle}

\newcommand{\seqApn}[1]{\langle A_n,p_n\mid 0\le n\le {#1}\rangle}

\newcommand{\hist}[1]{\langle A_0,p_0,\dots,A_{#1}\rangle}
\newcommand{\Prb}{\mathbb{P}}
\newcommand{\Unif}{\mathrm{Unif}}

\makeatletter
\newcommand*{\defeq}{\mathrel{\rlap{%
                     \raisebox{0.3ex}{$\m@th\cdot$}}%
                     \raisebox{-0.3ex}{$\m@th\cdot$}}%
                     =}
\makeatother

\makeatletter
\DeclareRobustCommand\widecheck[1]{{\mathpalette\@widecheck{#1}}}
\def\@widecheck#1#2{%
    \setbox\z@\hbox{\m@th$#1#2$}%
    \setbox\tw@\hbox{\m@th$#1%
       \widehat{%
          \vrule\@width\z@\@height\ht\z@
          \vrule\@height\z@\@width\wd\z@}$}%
    \dp\tw@-\ht\z@
    \@tempdima\ht\z@ \advance\@tempdima2\ht\tw@ \divide\@tempdima\thr@@
    \setbox\tw@\hbox{%
       \raise\@tempdima\hbox{\scalebox{1}[-1]{\lower\@tempdima\box
\tw@}}}%
    {\ooalign{\box\tw@ \cr \box\z@}}}
\makeatother

\title{Determining the Winner in Alternating-Move Games\thanks{We thank Tom Meyerovitch for suggesting Theorem~\ref{FlipCoinTheorem}, Tushar Das and David Simmons for helpful discussions and for suggesting an approach for proving Theorem~\ref{MainSchmidt}, and Omri N. Solan for his useful insight that led to the proof of Theorem \ref{unrestricted_upper_bound}. We are grateful to Eilon Solan, Barak Weiss, and Ehud Lehrer for their guidance, advice, and helpful comments on earlier versions of this paper. This work was supported by the Israel Science Foundation (grants No. 591/21 and No. 211/22). Refine.ink was used to check the paper for consistency and clarity.}}

\author{Itamar Bella\"iche\thanks{Department of Economics, New York University, itamarbellaiche@nyu.edu.} \and Auriel Rosenzweig\thanks{School of Mathematical Sciences, Tel Aviv University, aurielr@tauex.tau.ac.il.}}

\begin{document}
\maketitle

\begin{abstract}
We provide a criterion for determining the winner in two-player win-lose alternating-move games on trees, in terms of the Hausdorff dimension of the target set. We focus our study on special cases, including the Gale-Stewart game on the complete binary tree and a family of Schmidt games, generalizing a result of Schmidt from Hilbert spaces to arbitrary complete metric spaces. Building on the Hausdorff dimension games originally introduced by Das, Fishman, Simmons, and Urba\'nski, which provide a game-theoretic approach for computing Hausdorff dimensions, we employ a generalized family of these games to obtain lower bounds on the Hausdorff dimensions of target sets whenever Player I can guarantee a win.
\end{abstract}

\begin{sloppypar}

\section{Introduction}
\label{SectionIntroduction}
\paragraph{Win-lose alternating-move games.}
Win-lose alternating-move games over a tree were first introduced by Gale and Stewart (1953, \cite{GaleStewart1953}).\footnote{These games are also known as \emph{Gale-Stewart games} or \emph{Borel games}.} Such a game is denoted $G\defeq\left(T,W\right)$, where $T$ is a tree, \emph{i.e.}, a set of finite sequences, closed under prefixes, with the property that every sequence $p\in T$ has an extension in $T$, and $W$ is a set of $\omega$-sequences whose every finite prefix is in $T$. The game is played in stages $n\in \mathbb{N}\defeq\left\{0,1,2,...\right\}$ as follows:
\begin{itemize}
\item At $n=0$, Player I chooses $a_0$ such that $\langle a_0\rangle \in T$.
\item If $n$ is odd, and $a_0,...,a_{n-1}$ have already been chosen such that $\langle a_0,...,a_{n-1}\rangle\in T $, then Player II chooses $a_n$ such that $\langle a_0,...,a_{n-1},a_n\rangle \in T$.
\item If $2\le n$ is even, and $a_0,...,a_{n-1}$ have already been chosen such that $\langle a_0,...,a_{n-1}\rangle\in T $, then Player I chooses $a_n$ such that $\langle a_0,...,a_{n-1},a_n\rangle \in T$.
\end{itemize}
Player I wins if $\langle a_0,a_1,...\rangle\in W $, Player II wins otherwise.\\
The set $W$ is called \emph{determined} if one of the players has a winning strategy in the game $\left(T,W\right)$; namely, a strategy that guarantees that the player wins, whatever strategy her opponent adopts.

Gale and Stewart (1953, \cite{GaleStewart1953}) proved that when endowing the canopy of the tree $T$, which is the space $\canopy$ of all $\omega$-sequences whose every finite prefix is in $T$, with the product topology, all open or closed sets are determined. They also proved that there exist non-determined subsets of the canopy of the \emph{complete binary tree}, \emph{i.e.}, the set of all finite binary sequences, which we denote throughout the paper by $\Dya\defeq\left\{0,1\right\}^{<\omega }$. Their result was extended to higher Borel hierarchies, until Martin (1975, \cite{Martin75Paper}) proved that every Borel subset of the canopy of every tree is determined.
\paragraph{Who wins the game?}
While determinacy proves that \emph{one} of the players has a winning strategy, it does not provide information about \emph{which} of the players has it. The goal of this paper is to shed light on the latter question. We will begin by studying games over $\Dya$. These games were called \emph{binary games} by Gale and Stewart (1953, \cite{GaleStewart1953}) and $2$-\emph{digit-games} by Schmidt (1966, \cite{Schmidt1966}).

\begin{example}[label=CantorExample]
Let $\mathcal{C}\subseteq \left[0,1\right]$ be the Cantor-like set of numbers which have a base 5 expansion containing only the digits 1 and 3, \emph{i.e.}, 
$$\mathcal{C}\defeq\left\{\sum_{n\in \mathbb{N}}\frac{a_n}{5^{n+1}}\in \left[0,1\right]\ \middle|\ a_n\in \left\{1,3\right\}\forall n\in \mathbb{N}\right\}.$$

Let
$$W_\mathcal{C}\defeq\left\{\langle a_0,a_1,...\rangle \in \Dyacanopy\ \middle|\ \sum_{n\in \mathbb{N}}\frac{a_n}{2^{n+1}}\in \mathcal{C}\right\}.$$
Thus, in the game $\left(\Dya, W_{\mathcal{C}}\right)$, the two players alternately choose bits $0$ or $1$, and Player I wins if the number $x$, whose binary expansion coincides with the sequence of bits chosen by the players, lies in $\mathcal{C}$. The set $W_\mathcal{C}$ is Borel, and thus by Martin (1975, \cite{Martin75Paper}), determined. The reader may try to decide which of the two players has a winning strategy in the game $\left(\Dya,W_\mathcal{C}\right)$.
\end{example}

Intuitively, as the target set $W$ gets smaller, Player II is more likely to have a winning strategy. For instance, one size criterion of a set $W$ is its cardinality. In fact, if $\left|W\right|<\mathfrak{c}$, where $\mathfrak{c}$ is the cardinality of the continuum, then Player II can guarantee a win in the game $\left(\Dya,W\right)$. Indeed, each strategy $s_{II}$ of Player II generates the set $\llbracket s_{II}\rrbracket$ of possible $\omega $-sequences that can be constructed in the game if Player II is playing according to $s_{II}$. Player II can guarantee a win in the game $\left(\Dya,W\right)$ if and only if there exists a strategy $s_{II}$ of Player II such that $\llbracket s_{II}\rrbracket \subseteq W^c$. Since there exists a collection $\left\{s_{II}^\alpha \ \middle|\ \alpha <\mathfrak{c} \right\}$ of strategies of the cardinality of the continuum satisfying $\llbracket s^\alpha _{II}\rrbracket\cap \llbracket s^\beta _{II}\rrbracket=\emptyset $ if $\alpha \ne \beta $, if $\left|W\right|<\mathfrak{c}$, then there exists $\alpha<\mathfrak{c}$ such that $ s^\alpha _{II}$ is a winning strategy of Player II in the game $\left(\Dya,W\right)$.

The question that arises is whether there exists a topological-metric sufficient criterion for Player II to have a winning strategy. One topological criterion for classifying a set ``small'' is Baire category. However, as the next example shows, this criterion turns out to be irrelevant in the sense of classifying a set to be a winning set of Player II.
\begin{example}
\label{Y^0} 
Define
$$\llbracket \Dya^0 \rrbracket\defeq\left\{\langle a_0,a_1,...\rangle\in \Dyacanopy \ \middle|\ a_{2n}=0,\forall n\in \mathbb{N}\right\}\subseteq \Dyacanopy.$$
Player I can guarantee a win in the game $\left(\Dya,\llbracket \Dya^0 \rrbracket\right)$ by choosing $0$ in every step of the game. The set $\llbracket \Dya^0\rrbracket$ is a closed nowhere dense set in $\Dyacanopy$, and, in particular, of Baire category 1.
\end{example}
\paragraph{Main result.}
In this paper we will measure the size of a set by its \emph{Hausdorff dimension}. The topological space $\Dyacanopy$ is a metric space with the metric
$$d_2\left(x,y\right)\defeq\begin{cases}
0, & x=y,\\
2^{-\min\left\{n\in \mathbb{N}\ \middle|\  a_n\ne b_n\right\}}, &x\ne y,
\end{cases}$$
where $x=\langle a_0,a_1,...\rangle $ and $y=\langle b_0,b_1,...\rangle $. Our first main result states that Player II has a winning strategy if the Hausdorff dimension of $W$ with respect to the metric $d_2$ is smaller than $\frac{1}{2}$.
\begin{theorem}
\label{MainTheoremSmall} 
Let $W\subseteq \Dyacanopy$ such that $\dim_{\mathcal{H}}\left(W\right)<\frac{1}{2}$. Then Player II can guarantee a win in the game $G=\left(\Dya,W\right)$, and, in particular, $W$ is determined.
\end{theorem}
Analogously to the proof of the claim that if $\left|W\right|<\mathfrak{c}$ then Player II can guarantee a win in the game $G=\left(\Dya,W\right)$, the main idea of the proof of Theorem~\ref{MainTheoremSmall} is that for $W$ not to be a winning set of Player II, it must contain a sufficiently large subset. We start by showing that for every strategy $s_I$ of player I, the Hausdorff dimension of the set $\llbracket s_I\rrbracket $ of possible $\omega $-sequences that can be constructed in the game if Player I is playing according to $s_I$, is $\frac{1}{2}$. Thus, Player I cannot have a winning strategy $s_I$ in the game $G$ if $\dim_{\mathcal{H}}\left(W\right)<\frac{1}{2}$, since this will imply that $\llbracket s_I\rrbracket \subseteq W$, and by monotonicity of the Hausdorff dimension
$$\frac{1}{2}=\dim_{\mathcal{H}}\left(\llbracket s_I\rrbracket \right)\le \dim_{\mathcal{H}}\left(W\right)<\frac{1}{2},$$
a contradiction. Using the Borel regularity of the Hausdorff dimension\footnote{The Borel regularity of the Hausdorff dimension states that for every set $A$, there exists a Borel set $B$ such that $A\subseteq B$, and $\dim_{\mathcal{H}}\left(A\right)=\dim_{\mathcal{H}}\left(B\right)$.} and combining it with Martin \cite{Martin75Paper}, we deduce that Player II has a winning strategy.
\begin{example}[continues=CantorExample]
Since $\dim_{\mathcal{H}}\left(W_\mathcal{C}\right)=\log_{5}\left(2\right)<\frac{1}{2}$, Theorem~\ref{MainTheoremSmall} implies that Player II has a winning strategy in the game $G=\left(\Dya,W_\mathcal{C}\right)$.
\end{example}
\begin{remark}
Theorem~\ref{MainTheoremSmall} is proved without any measurability assumptions regarding the set $W$, and thus this theorem not only gives us information about who can guarantee a win in the game, but also extends the class of known determined subsets of the canopy of the complete binary tree. Furthermore, in Subsection~\ref{MonotonicFunctions} we show how the method used in the proof enables us to employ alternative size criteria given by monotone, Borel regular functions.
\end{remark}
We will also present a method to construct examples of target sets of Hausdorff dimension $\delta $, for every $\delta \in \left[\frac{1}{2},1\right]$, which are winning sets for either of the players, or non-determined. Consequently, the Hausdorff dimension of a target set gives us information about the determinacy of the game, only when it is strictly less than $\frac{1}{2}$. This result is also analogous to the cardinality criterion, since a set of the cardinality of the continuum may be a winning set for either of the players, or it may be non-determined.
\paragraph{Schmidt games.}
The use of Hausdorff dimension as a criterion for determining whether a player has a winning strategy in a win-lose alternating-move game is not new. It previously appeared in the context of \emph{Schmidt games on complete metric spaces}. Schmidt games, in the broad sense, were first formulated by Schmidt (1966, \cite{Schmidt1966}). Schmidt introduced a general class of win-lose alternating-move games ``played on'' mathematical spaces, usually paired with an additional structure such as a topology, a metric, or a measure.

A Schmidt game on a space $X$ is a win-lose alternating-move game played with respect to a predetermined target set $S\subseteq X$. At each stage of the game, the player whose turn it is to move needs to choose one subset from a collection of available subsets of the previously chosen set. Player I wins if the intersection of the resulting nested sequence of chosen sets is contained in $S$. Player II wins otherwise. From this perspective, these games can be viewed as a special case of the games introduced by Gale and Stewart, where the tree of the game consists of finite sequences of nested subsets of the space $X$.
\paragraph{$\left(\alpha ,\beta \right)$-Schmidt games.}
A particular family of Schmidt games is the \emph{$\left(\alpha ,\beta \right)$-Schmidt games}, introduced in \cite{Schmidt1966}. Let $\left(X,d\right)$ be a complete metric space, let $S\subseteq X$, and let $\alpha,\beta \in \left(0,1\right)$. The $\left(\alpha ,\beta \right)$-Schmidt game on the space $X$ with target set $S$ is played as follows:
\begin{itemize}
\item At the beginning of the game, Player II chooses a closed ball $B_{-1}\subseteq X$ of radius $\rho_{-1}$.\footnote{We follow the convention from Schmidt's paper, where Player II (Black, in the original paper) is the one making the first choice. This is in contrast with the usual game-theoretic convention where Player I is the first to move. We can think of these games as ones where Player II chooses the game, and then the players play a standard win-lose alternating-move game.}
\item If $n$ is even, and $B_{n-1}$ is a closed ball chosen by Player II in the previous stage, at stage $n$ Player I chooses a closed ball $B_n\subseteq B_{n-1}$ of radius $\rho_n=\alpha \cdot \rho_{n-1}$.
\item If $n$ is odd, and $B_{n-1}$ is a closed ball chosen by Player I in the previous stage, at stage $n$ Player II chooses a closed ball $B_n\subseteq B_{n-1}$ of radius $\rho_n=\beta \cdot \rho_{n-1}$.
\end{itemize}
Player I wins if $\bigcap_{n\in \mathbb{N}}B_n\subseteq S$, Player II wins otherwise.

Schmidt introduced this family of games as a new geometric and dynamical tool for studying fundamental sets in Diophantine approximation, such as the set of badly approximable numbers. Although such sets are small in the measure-theoretic sense (having Lebesgue measure zero), they are invariant under many natural transformations, suggesting a different notion of largeness.

The $\left(\alpha ,\beta \right)$-Schmidt game gave Schmidt a framework to formalize this notion of largeness. Schmidt called a set $S\subseteq X$ an \emph{$\alpha$-winning set} if Player I has a winning strategy in the $\left(\alpha ,\beta \right)$-Schmidt game with a target set $S$, for every $\beta \in \left(0,1\right)$. His key insight was that $\alpha $-winning sets have strong stability properties, such as density and full Hausdorff dimension. Furthermore, the class of $\alpha$-winning sets is closed under countable intersections. Among these results, Schmidt proved the following theorem.

\begin{theorem}[\cite{Schmidt1966}, Corollary 1]
\label{SchmidtHausdorff} 
Let $X$ be a Hilbert space, let $\alpha,\beta \in \left(0,1\right)$, and let $N(\beta)$ be such that every ball of radius $\rho $ in $X$ contains $N(\beta)$ balls of radius $\beta \cdot \rho $ with pairwise disjoint interiors, for every $\rho >0$. Let $S\subseteq X$ such that 
$$\dim_{\mathcal{H}}\left(S\right)<\log_{\left(\alpha \beta\right)^{-1} }\left(N(\beta \right)).$$
Then Player I does not have a winning strategy in the $\left(\alpha ,\beta \right)$-Schmidt game with target set $S$. 
\end{theorem}
Combining Theorem~\ref{SchmidtHausdorff} with Martin's theorem \cite{Martin75Paper} and the Borel regularity of the Hausdorff dimension, we obtain that in the setting of Theorem~\ref{SchmidtHausdorff}, Player II has a winning strategy in the $\left(\alpha ,\beta \right)$-Schmidt game with target set $S$.

From a topological point of view, win-lose alternating-move games on trees can be interpreted as Schmidt games on the metric space corresponding to the tree's canopy, where at each stage a player chooses a closed ball contained in the previously chosen one. In this sense, Schmidt games on general complete metric spaces can be viewed as a generalization of win-lose alternating-move games on trees.

Indeed, for every $W\subseteq \Dyacanopy$, the game $\left(\Dya,W\right)$ can be interpreted as a subgame of the $\left(\frac{1}{2} ,\frac{1}{2} \right)$-Schmidt game with target set $W$, in which Player II is restricted to choose the ball of radius $1$ at the beginning of the game, and the threshold appearing in Theorem~\ref{MainTheoremSmall} aligns with the one in Theorem~\ref{SchmidtHausdorff}. Note that Theorem~\ref{MainTheoremSmall} does not follow from Theorem~\ref{SchmidtHausdorff}, since $\Dyacanopy$ is not a Hilbert space. In Section~\ref{SchmidtSection} we study \emph{Schmidt games on complete metric spaces}, and state a generalization of Theorem~\ref{SchmidtHausdorff} to the setting of an arbitrary complete metric space rather than a Hilbert space.
\paragraph{Hausdorff dimension games.}
The proof of Theorem~\ref{SchmidtHausdorff} in \cite{Schmidt1966} relies on geometric properties of Hilbert spaces, such as half-spaces defined using an inner product, to obtain a lower bound for the Hausdorff dimension of the target set, when Player I can guarantee a win in the game. In this paper, one of the techniques used to obtain the lower bound employs \emph{Hausdorff dimension games}. Dimension games were first introduced by Das, Fishman, Simmons, and Urba\'nski (2024, \cite{DimensionPaper}) as part of their development of a variational principle in the parametric geometry of numbers. Their framework provides a new method for computing the Hausdorff and packing dimension of Borel subsets of a complete doubling\footnote{A complete doubling metric space is a complete metric space with a constant $D\in\N\setminus\{0\}$ such that for every radius $r\in(0,\infty)$ and every point $x\in X$, the closed ball $\closedb(x,r)$ can be covered by at most $D$ closed balls of radius $\tfrac{r}{2}$.} metric space. In particular, they used these games to compute the Hausdorff and packing dimension of sets arising in Diophantine approximation, including the set of singular matrices.

Let $\left(X,d\right)$ be a complete doubling metric space, and let $\emptyset\ne S\subseteq X$. The Hausdorff dimension game with a target set $S$ is played as follows: In each step of the game, Player 1 presents Player 2 with a non-empty collection of closed balls of equal radius, all contained in the previously chosen ball, and Player 2 chooses one of them. The radii of the balls offered in each step decay exponentially, at a fixed rate chosen by Player 1 at the beginning of the game. At the end of the game, if the intersection of the resulting nested sequence of closed balls chosen by Player 2 is contained in $S$, Player 2 pays Player 1 the limit inferior of the long-run average of the logarithm of the number of balls Player 1 offers Player 2. Otherwise, Player 1 pays $1$ to Player 2.

In \cite{DimensionPaper}, it was proved that, under suitable restrictions on the collection of balls that Player 1 may offer, the Hausdorff dimension game with a target $S$ has a value whenever $S$ is Borel. Furthermore, the value of the game coincides with the Hausdorff dimension of $S$. In Section~\ref{HDG}, we formally define the Hausdorff dimension games in the setting of the complete binary tree's canopy, and, applying the result from \cite{DimensionPaper}, use these games to compute the Hausdorff dimension of various subsets. This forms the main technical part in the proof of Theorem~\ref{MainTheoremSmall}.

In Section~\ref{GeneralHDG} we introduce a more general family of Hausdorff dimension games, and state a generalization of the result from \cite{DimensionPaper}. This generalization enables us to apply the same techniques used in the proof of Theorem~\ref{MainTheoremSmall} to establish an analogous version of Theorem~\ref{SchmidtHausdorff} in the setting of complete doubling metric spaces which is applicable to a broad class of Schmidt games.
\paragraph{}
The paper is structured as follows. In Section~\ref{Preli}, we formally present all the definitions and main theorems regarding win-lose alternating-move games and the Hausdorff dimension. In Section~\ref{HDG}, we present the Hausdorff dimension game in the setting of the complete binary tree. In Section~\ref{Main}, we state and prove a general criterion for determining the winner in win-lose alternating-move games, using monotone Borel regular functions, and apply it to prove Theorem~\ref{MainTheoremSmall}. In Section~\ref{Ext} we present some additional results regarding the connections between the Hausdorff dimension and determining the winner in win-lose alternating-move games.

In Section~\ref{SchmidtSection} we introduce a general formulation of Schmidt games as a family of win-lose alternating-move games on trees. We formally present the $\left(\alpha ,\beta \right)$-Schmidt game in the setting of an arbitrary complete metric space,  and state generalizations of Theorems~\ref{SchmidtHausdorff} and~\ref{MainTheoremSmall}. The proofs of these versions follow the same main ideas as the proof of Theorem~\ref{MainTheoremSmall}, with additional steps to accommodate the more general settings, and are provided in Appendix~\ref{DimensionGamesProofsSection}. In Section~\ref{GeneralHDG} we define the notion of imposed subtrees, introduce a more general family of Hausdorff dimension games, and state a generalization of the result in \cite{DimensionPaper}.

\section{Preliminaries}
\label{Preli} 
In this section we formally present all the definitions and the results we will need regarding alternating-move games and Hausdorff dimension.
\subsection{Alternating-Move Games}
In this subsection we define alternating-move games and state the known results regarding their determinacy. For further background, we refer the reader to Solan (\cite{BorelGamesBook}, 2025).

Let $\mathcal{A}$ be a set, termed the \emph{alphabet}. The set $\mathcal{A}^{<\omega}$ is the set of all finite sequences of elements in $\mathcal{A}$. An element in $\mathcal{A}^{<\omega}$ is denoted $p\defeq\langle a_0,a_1,...,a_n\rangle $. The \emph{length} of $p$ is $\len\left(p\right)\defeq n+1$. \emph{The empty sequence} in $\mathcal{A}^{<\omega}$ is denoted by $\langle \text{ }\rangle $, and $\len\left(\langle \text{ }\rangle \right)\defeq 0$.\\
For $p=\langle a_0,...,a_n\rangle ,p^\prime=\langle a^\prime_0,...,a^\prime_m\rangle \in \mathcal{A}^{<\omega}$, we define their concatenation by
$$p\frown p^\prime\defeq \langle a_0,...a_n,a^\prime_0,...,a^\prime_m\rangle .$$
If $p$ is a prefix of $q$, we say that $q$ \emph{extends} $p$, and write $p\preceq q$.

\begin{definition}[tree and positions]
\label{def:tree}
A \emph{tree} over an alphabet $\mathcal{A}$ is a non-empty set $T\subseteq\mathcal{A}^{<\omega}$, satisfying
\begin{itemize}
\item If $p\in T$ and $q\preceq p$, then $q\in T$.
\item For every $p\in T$ there exists $q\in T$ such that $p\preceq q$, and $\len\left(q\right)=\len\left(p\right)+1$.
\end{itemize}
An element of $T$ is called a $\emph{position}$. For every $p\in T$, denote by $\mathcal{A}_p\left(T\right)$ the set of available alphabet elements in position $p$, \emph{i.e.},
$$\mathcal{A}_p\left(T\right)\defeq\left\{a\in \mathcal{A}\ \middle|\  p\frown \langle a\rangle \in T\right\}.$$
A set $T'\subseteq T$ which is also a tree is called a \emph{subtree} of $T$. For every $p\in T$, denote by
$$T_p\defeq\left\{q\in T\ \middle|\ \text{$q\preceq p$ or $p\preceq q$}\right\}$$
the subtree of $T$ that consists of prefixes and extensions of $p$.
\end{definition}

We denote by $\mathcal{A}^{\omega}$ the set of $\omega $-sequences of elements in $\mathcal{A}$. An element of $\mathcal{A}^{\omega}$ is denoted $x=\langle a_0,a_1,...\rangle $.

\begin{definition}[plays and canopy]
Let $T$ be a tree over an alphabet $\mathcal{A}$. A \emph{play} in $T$  is an element of $\langle a_0,a_1,...\rangle\in \mathcal{A}^{\omega}$ satisfying $\langle a_0,...,a_{n-1}\rangle\in T$, for every $n\in \mathbb{N}$. We denote by $\canopy$, the set of all plays in $T$. $\canopy$ is called the \emph{canopy} of the tree $T$.\footnote{The canopy is also known as the branch space, body, or boundary of the tree.} We endow $\canopy$ with the product topology,\footnote{The product topology is defined based on the discrete topology on $\mathcal{A}$.} the topology generated by the basis of \emph{prefix cylinders}
$$\left\{\left\llbracket T_p\right\rrbracket\ \middle|\ p\in T\right\}.$$
We denote by $\mathcal{B}$ the Borel tribe, or Borel $\sigma$-algebra, generated by this topology.
\end{definition}
Note that in the case of the canopy of the \emph{complete binary tree}, $\Dya\defeq\left\{0,1\right\}^{<\omega }$, the topology induced by the dyadic metric $d_2$ is identical to the product topology, and the set of clopen balls under $d_2$ coincides with the set of prefix cylinders.

Gale and Stewart (1953, \cite{GaleStewart1953}) introduced zero-sum alternating-move games over a tree $T$.
\begin{definition}[zero-sum alternating-move game]
Let $\mathcal{A}$ be an alphabet, and let $T$ be a tree over $\mathcal{A}$. Let $f:\canopy\to \mathbb{R}$ be a function. A \emph{(zero-sum alternating-move) game} on $T$ is a pair $G\defeq\left(T,f\right)$. The game is played in stages $n\in \mathbb{N}$ as follows:
\begin{itemize}
\item In stage $0$, Player I chooses $a_0\in \mathcal{A}_{\langle \text{ }\rangle }\left(T\right)$.
\item If $n\in \mathbb{N}_{\text{odd}}=\left\{1,3,5,...\right\}$, and $p=\langle a_0,...,a_{n-1}\rangle\in T $ has been chosen, Player II chooses $a_n\in \mathcal{A}_p\left(T\right)$ in stage $n$.
\item If $2\le n\in \mathbb{N}_{\text{even}}=\left\{0,2,4,...\right\}$, and $p=\langle a_0,...,a_{n-1}\rangle\in T $ has been chosen, Player I chooses $a_n\in \mathcal{A}_p\left(T\right)$ in stage $n$.
\end{itemize}
Player II pays Player I $f\left(\langle a_0,a_1,...\rangle \right)$. The function $f$ is called the \emph{payoff function}. If the payoff function $f$ is a characteristic function of a set $W\subseteq \canopy$, we call the game $G=\left(T,W\right)\defeq \left(T,\vmathbb{1}_W\right)$ a \emph{win-lose alternating-move game}. If $T^\prime\subseteq T$ is a subtree, the pair $\left(T',f|_{\llbracket T'\rrbracket}\right)$ is called a \emph{subgame} of the zero-sum alternating-move game $\left(T,f\right)$.
\end{definition}
A \emph{pure strategy} $s_I$ of Player I is a function that assigns to every position $p\in T$ of even length, an element $s_I\left(p\right)\in \mathcal{A}_p\left(T\right)$. We denote the set of pure strategies of Player I by $S_I\left(T\right)$. A \emph{pure strategy} $s_{II}$ of Player II is a function that assigns to every position $p\in T$ of odd length, an element $s_{II}\left(p\right)\in \mathcal{A}_p\left(T\right)$. We denote the set of pure strategies of Player II by $S_{II}\left(T\right)$.\\
Given a pure strategy $s_I\in S_I\left(T\right)$ of Player I and a pure strategy $s_{II}\in S_{II}\left(T\right)$ of Player II, we denote by $\langle s_I,s_{II}\rangle$ the play generated by these two strategies, \emph{i.e.}, the sequence of actions taken by the players when choosing in each stage a pure action according to these strategies. Furthermore, we denote by
$$\llbracket s_I\rrbracket\defeq\left\{\langle s_I,s_{II}\rangle\ \middle|\ s_{II}\in S_{II}\left(T\right)\right\},$$
the subcanopy generated by the strategy $s_I$, and similarly 
$$\llbracket s_{II}\rrbracket\defeq\left\{\langle s_I,s_{II}\rangle\ \middle|\ s_I\in S_I\left(T\right)\right\},$$
the subcanopy generated by the strategy $s_{II}$.\\
Given a game $G=\left(T,f\right)$, we say that a pure strategy $\widehat{ s}_I$ of Player I guarantees $z$ in $G$ for Player I if for every pure strategy $s_{II}$ of Player II,
$$f\left(\langle\widehat{ s}_I,s_{II}\rangle\right)\ge z.$$
A pure strategy $\widehat{ s}_{II}$ of Player II guarantees $z$ in $G$ for Player II if for every pure strategy $s_I$ of Player I,
$$f\left(\langle s_I,\widehat{s}_{II}\rangle\right)\le z.$$
The real number $v$ is the \emph{value} of the zero-sum alternating-move game $G$ if for every $\varepsilon>0$ there exist pure strategies $ s^\varepsilon_I$ and $ s^\varepsilon_{II}$ that guarantee $v-\varepsilon$ and $v+\varepsilon$ for Player I and for Player II, respectively. If such a real number exists, then it is unique, and denoted by $\val\left(G\right)$. In such a case the game is said to \emph{have a value}, or is called \emph{determined}.\\
If $G=\left(T,W\right)$ is a win-lose alternating-move game, and $\val\left(G\right)=1$, we say that Player I can \emph{guarantee a win in the game $G$}. Note that this happens if and only if there exists $\widehat{ s}_I\in S_I\left(T\right)$ such that $\llbracket \widehat{s}_I\rrbracket\subseteq W $. Similarly, if $\val\left(G\right)=0$, we say that Player II can \emph{guarantee a win in the game $G$}. This happens if and only if there exists $\widehat{s}_{II}\in S_{II}\left(T\right)$ such that $\llbracket \widehat{ s}_{II}\rrbracket\subseteq W^c $. In either of these cases, we say that $W$ is \emph{determined}.

\begin{theorem}[Gale and Stewart (1953, \cite{GaleStewart1953})]
\label{TheoremGaleStewart} 
$\text{ }$
\begin{enumerate}
\item Let $T$ be a tree over an alphabet $\mathcal{A}$, and let $W\subseteq \canopy$ be either an open or a closed subset. Then the game $\left(T,W\right)$ is determined.
\item There exists a set $B\subseteq \Dyacanopy$ such that the game $\left(\Dya,B\right)$ is non-determined.
\end{enumerate}
\end{theorem}
Martin (1975, \cite{Martin75Paper}) extended the first statement in Theorem~\ref{TheoremGaleStewart}.
\begin{theorem}[Martin (1975, \cite{Martin75Paper})]
\label{Martin75} 
Let $T$ be a tree over an alphabet $\mathcal{A}$, and let $f:\canopy\to \mathbb{R}$ be a bounded Borel measurable function. Then the game $\left(T,f\right)$ is determined.
\end{theorem}

Let $T$ be a tree, and suppose that $\left|\mathcal{A}_p\left(T\right)\right|\le \aleph_0$, for every $p\in T$. A \emph{behavior strategy} $\sigma_I$ of Player I is a function that assigns to every position $p\in T$ of even length, a distribution $\sigma_I\left(p\right)\in \Delta\left(\mathcal{A}_p\left(T\right)\right)$. We denote the set of behavior strategies of Player I by $\Sigma_I\left(T\right)$. A \emph{behavior strategy} $\sigma_{II}$ of Player II is a function that assigns to every position $p\in T$ of odd length, a distribution $\sigma_{II}\left(p\right)\in \Delta\left(\mathcal{A}_p\left(T\right)\right)$. We denote the set of behavior strategies of Player II by $\Sigma_{II}\left(T\right)$.\\
Given a behavior strategy $\sigma_I\in \Sigma_I\left(T\right)$ of Player I and a behavior strategy $\sigma_{II}\in \Sigma_{II}\left(T\right)$ of Player II, we denote by $\Prb_{\sigma_I,\sigma_{II}}$ the probability measure on $\canopy$ induced by the pair $\left(\sigma_I,\sigma_{II}\right)$.

\subsection{Hausdorff Dimension}
In this subsection we define the Hausdorff dimension and list some of its properties.
\begin{definition}[$\delta $-dimensional Hausdorff outer measure]
Let $\left(X,d\right)$ be a metric space, let $A\subseteq X$, let $\delta \in \left[0,\infty\right)$, and let $\rho >0$. Define
$$\mathcal{H}_\rho^\delta\left(A\right)\defeq\inf\left\{\sum_{i\in \mathbb{N}}\left(\diam\left(A_i\right)\right)^\delta \ \middle|\ A\subseteq \bigcup_{i\in \mathbb{N}}A_i, \diam\left(A_i\right)<\rho,\forall i\in \mathbb{N}\right\}.$$
The \emph{$\delta $-dimensional Hausdorff outer measure corresponding to $d$} is
$$\mathcal{H}^\delta\left(A\right) \defeq\sup_{\rho>0}\mathcal{H}_\rho^\delta\left(A\right).$$
\end{definition}

\begin{definition}[Hausdorff dimension]
\label{HausdorffDimension} 
Let $\left(X,d\right)$ be a metric space, and let $A\subseteq X$. The \emph{Hausdorff dimension} of $A$ is
$$\dim_{\mathcal{H}}\left(A\right)\defeq \inf \left\{\delta \in \left[0,\infty \right)\ \middle|\ \mathcal{H}^\delta \left(A\right)=0\right\},$$
where, by convention, $\inf\emptyset =\infty $.
\end{definition}

\begin{remark}
If we endow $\Dyacanopy$ with the metric $d_2$, we have $\dim_{\mathcal{H}}\left(\Dyacanopy\right)=1$.
\end{remark}
A key property of the Hausdorff dimension is its monotonicity: if $A\subseteq B\subseteq X$, then $\dim_{\mathcal{H}}\left(A\right)\le \dim_{\mathcal{H}}\left(B\right)$.

By the Borel regularity of the $\delta $-dimensional Hausdorff measure $\mathcal{H}^\delta$, (\cite[Theorem 2.1]{HausdorffDimensionBook}, whose proof carries over from $\mathbb{R}^n$ to a general metric space $\left(X,d\right)$ without change), for every $A\subseteq X$ there exists a Borel set $A\subseteq Z\subseteq X$, such that $\mathcal{H}^\delta\left(Z\right)=\mathcal{H}^\delta\left(A\right)$. Thus, we conclude the following corollary.
\begin{corollary}
\label{BorelSameDimension} 
Let $\left(X,d\right)$ be a metric space, and let $A\subseteq X$. Then, there exists a Borel set $A\subseteq Z\subseteq X$, such that
$$\dim_{\mathcal{H}}\left(Z\right)=\dim_{\mathcal{H}}\left(A\right).$$
\end{corollary}
The proof of Corollary~\ref{BorelSameDimension} follows from the Borel regularity of the $\delta $-dimensional Hausdorff measure $\mathcal{H}^\delta$, by taking $Z\defeq X$ if $\dim_{\mathcal{H}}\left(A\right)=\infty$, and, otherwise, letting $Z_n\subseteq X$ be such that $A\subseteq Z_n$ and
$$\mathcal{H}^{\dim_{\mathcal{H}}\left(A\right)+\frac{1}{n}}\left(Z_n\right)=\mathcal{H}^{\dim_{\mathcal{H}}\left(A\right)+\frac{1}{n}}\left(A\right)=0,$$
and defining $Z\defeq\bigcap_{1\le n\in \mathbb{N}}Z_n$.

\section{Hausdorff Dimension Games}
\label{HDG} 
Hausdorff (and packing) dimension games were presented for the first time by Das, Fishman, Simmons, and Urba\'nski (2024, \cite{DimensionPaper}), and were since studied by Solan (2021, \cite{OmriSolan}), Crone, Fishman, and Jackson (2022, \cite{UnflodingDimensionGame}) and Bella\"iche (2023, \cite{BellaicheThesis}). In this section, we focus on a simple special case: the \emph{dyadic Hausdorff dimension game} on the complete binary tree. In Section \ref{GeneralHDG}, we study some more general variants.

The Hausdorff dimension quantifies the density of a set at small scales. That is, as we zoom in to a finer and finer resolution, how much freedom remains for constructing points in the set. The following game formalizes this intuition in the complete binary tree by tying together the Hausdorff dimension of a set and the limit average of the freedom available in the infinite process of refining a cover of the set by dyadic prefix cylinder sets.

\subsection{The Dyadic Hausdorff Dimension Game}
In this subsection we present the dyadic Hausdorff dimension game, and explain the intuition for the connection between this game and the Hausdorff dimension of subsets of $\Dyacanopy$.
Given a non-empty target set $\emptyset \ne S\subseteq\Dyacanopy$, the \emph{dyadic Hausdorff dimension game} is a zero-sum alternating-move game denoted by $G_S$ or $(G,p_{_S})$ and is played as follows:
\begin{itemize}
	\item In step 0, Player 1 chooses some $k\in\N\setminus\{0\}$ and a non-empty subset $A_0\subseteq\kbits$ (stage 0 of the game). Then, Player 2 chooses one element $b_0\in A_0$ (stage 1 of the game).
	\item In every step $n\in\N\setminus\{0\}$, Player 1 chooses a non-empty subset $A_n\subseteq\kbits$ (stage $2n$ of the game) and then Player 2 chooses one element $b_n\in A_n$ (stage $2n+1$ of the game).
\end{itemize}
According to this terminology, a step in the game consists of two stages: the first played by Player 1 and the second by Player 2.
Denote the projection
\begin{equation*}
	Y_\infty\colon\bigcup_{k\in\N\setminus\{0\}}\biggl(\{k\}\times\Bigl(\bigl(\powerset\left(\kbits\right)\setminus\{\emptyset\}\bigr)\times\kbits\Bigr)^\N\biggr)\longrightarrow\Dyacanopy,
\end{equation*}
where $\powerset\left(\kbits\right)$ is the power set of $\kbits$, by
\begin{equation}
    \label{projection} 
	\langle a_0,a_1,...\rangle\defeq  Y_\infty\bigl(\langle k,\langle A_n,b_n\mid n\in\N \rangle\rangle\bigr)  =b_0\frown b_1\frown ...
\end{equation}
The payoff function of the game, denoted
\begin{equation*}
	p_{_S}\colon\bigcup_{k\in\N\setminus\{0\}}\biggl(\{k\}\times\Bigl(\bigl(\powerset\left(\kbits\right)\setminus\{\emptyset\}\bigr)\times\kbits\Bigr)^\N\biggr)\longrightarrow\mathbb{R},
\end{equation*}
is defined by
$$p_{_S}\bigl(\langle k,\langle A_n,b_n\mid n\in\N \rangle\rangle\bigr)\defeq \begin{cases}   \liminf_{N\rightarrow\infty}\tfrac{1}{N}\sum_{n=0}^{N-1}\log_{2^k}(|A_n|), & \langle a_0,a_1,...\rangle\in S,\\
			-1, & \text{otherwise}. \end{cases}$$

In words, in each step of the game, Player 1 offers Player 2 a non-empty collection of finite strings of the same length, and Player 2 chooses one of them. Each chosen string is concatenated to the previous chosen string, thereby constructing an infinite sequence. Player 1 has two objectives. The first goal is to ensure the resulting sequence is inside the target set $S$. Since $S\ne \emptyset $, she can always obtain this goal by offering a single string in each step. Her second goal is to maximize the long-run average number of options she offers in each step. To do so, she offers as many options as possible, while ensuring that no matter what option Player 2 chooses, Player 2 cannot force Player 1 to lower the long-run number of options Player 1 offers, in order to secure the first goal.

Fixing a strategy of Player 1, the payoff function, given by the limit average of the logarithm of the number of options Player 1 offers, records the exponential rate of the number of options Player 1 offers throughout the game. This construction mirrors the geometric intuition behind the Hausdorff dimension precisely because the progression of the game corresponds to the refinements of the cover of a subset of the target set $S$ by dyadic prefix cylinders, so that the average number of options presented by Player 1 to Player 2 corresponds to the freedom available at each finer scale.

\subsection{Value of the Dimension Game}
The nomenclature of ``dimension game'' was chosen in \cite{DimensionPaper} precisely because of the connection between the value of such games and the dimension of their target sets. Actually, the original reason behind the construction of the dimension game was providing a method of calculating the dimensions of some sets. The following theorem, proved originally in a more general setting in \cite{DimensionPaper}, states that the value of the dyadic Hausdorff dimension game is equal to the Hausdorff dimension of its target set, for the case of Borel target sets. For an explanation of how the dyadic setting we present here is a special case of \cite[Theorem 29.2.]{DimensionPaper}, the reader is referred to Bella\"iche (2023, \cite[Section 4.7]{BellaicheThesis}).
\begin{theorem}[Das, Fishman, Simmons, and Urba\'nski (2024, \cite{DimensionPaper})]
\label{ItamarTheorem} 
Let $\emptyset\neq S\in\mathcal{B}$ be a non-empty Borel target set. Then the dyadic Hausdorff dimension game with target set $S$ has a value, which is equal to the Hausdorff dimension of the set $S$, \emph{i.e.},
\begin{equation*}
    \val \left(G_S\right)=\dim_{\mathcal{H}}\left(S\right).
\end{equation*}
\end{theorem}
The value of the game captures the rate at which the freedom Player 1 has to offer options persists indefinitely, aligning with the geometric intuition of the Hausdorff dimension.

\subsection{Hausdorff Dimension as the Density of Unconstrained Coordinates}
Theorem~\ref{ItamarTheorem} suggests the intuition that the Hausdorff dimension of a subset $S\subseteq \Dyacanopy$ corresponds to the limiting density of coordinates at which sequences in $S$ may take either value $0$ or $1$. This intuition is formalized by the following family of closed subsets of $\Dyacanopy$.
\begin{definition}[the set $F_M$]
\label{F_M} 
Let $M\subseteq \mathbb{N}$. Define
$$F_M\defeq \left\{\langle a_0,a_1,...\rangle \in \Dyacanopy\ \middle|\ \text{$a_n = 0$, if $n\notin M$}\right\}.$$
\end{definition}
\begin{proposition}
\label{CalculateDimension} 
Let $M\subseteq \mathbb{N}$. Suppose that
$$\frac{\left|M\cap \left[0,n\right]\right|}{n+1}\xrightarrow[n\rightarrow\infty]{}\delta ,$$
where $\delta \in \left[0,1\right]$. Then $\dim_{\mathcal{H}}\left(F_M\right)=\delta $.
\end{proposition}
While Proposition~\ref{CalculateDimension} could be proved directly from the definition of the Hausdorff dimension, our proof uses the dimension game. This approach avoids some technical complications and highlights the intuition that the Hausdorff dimension of subsets of $\Dyacanopy$ corresponds to the density of unconstrained coordinates.
\begin{proof}
Consider the dyadic Hausdorff dimension game with a target set $F_M$, \emph{i.e.}, $G_{F_M}$. The set $F_M$ is closed. Indeed, for every $n\in \mathbb{N}$, the set $F_{\mathbb{N}\setminus \left\{n\right\}}$ is closed, and hence $F_M=\bigcap_{n\in M^c}F_{\mathbb{N}\setminus \left\{n\right\}}$ is closed. By Theorem~\ref{ItamarTheorem}, we have $\val\left(G_{F_M}\right)=\dim_{\mathcal{H}}\left(F_M\right)$. Thus, it is sufficient to show that $\val\left(G_{F_M}\right)=\delta $.\\
We start by proving that $\val\left(G_{F_M}\right)\ge \delta $. Let $\Gamma $ be the tree of the game $G_{F_M}$. Define a strategy $\widetilde{\sigma}_1\in S_1\left(\Gamma \right)$, as follows:
\begin{itemize}
\item Define
$$\widetilde{\sigma}_1\left(\langle \text{ }\rangle\right)\defeq \begin{cases}
    \langle 1,\left\{0,1\right\}\rangle , & 0\in M,\\
    \langle 1,\left\{0\right\}\rangle , & 0\notin M.
\end{cases}$$
\item Let $p=\langle \langle k,A_0\rangle ,b_0,A_1,b_1,...,A_{n-1},b_{n-1}\rangle \in \Gamma $. If $k=1$, define
$$\widetilde{\sigma}_1\left(p\right)\defeq \begin{cases}
    \left\{0,1\right\}, & n\in M,\\
    \left\{0\right\}, & n\notin M.
\end{cases}$$
Otherwise, define $\widetilde{\sigma}_1\left(p\right)$ arbitrarily.
\end{itemize}
Let $\sigma_{2}\in S_{2}\left(\Gamma \right)$ be any strategy of Player 2. Denote
$$\langle \langle 1,A_0\rangle ,b_0,A_1,b_1,...\rangle\defeq \langle \widetilde{\sigma}_1,\sigma_{2}\rangle  ,$$
and let $y\defeq \langle a_0,a_1,...\rangle \in \Dyacanopy$ be defined as in Equation~(\ref{projection}). Note that $y\in F_M$ by the definition of $\widetilde{\sigma}_1$. Therefore,
\begin{align*}
p_{F_M}\left(\langle\widetilde{\sigma}_1,\sigma_{2}\rangle \right)&=\liminf_{N\rightarrow\infty}\frac{1}{N}\cdot \sum_{n=0}^{N-1}\log_{2}\left(\left|A_n\right|\right)\\
&= \liminf_{N\rightarrow\infty}\frac{1}{N}\cdot \sum_{n=0}^{N-1}\vmathbb{1}_{\left\{n\in M\right\}}\\
&= \liminf_{N\rightarrow\infty}\frac{\left|M\cap \left[0,N-1\right]\right|}{N}\\
&= \lim_{N\rightarrow\infty}\frac{\left|M\cap \left[0,N-1\right]\right|}{N}=\delta .
\end{align*}
Thus, $\widetilde{\sigma}_1$ guarantees $\delta $.\\
Now, we show that $\val\left(G_{F_M}\right)\le \delta $. Define a strategy $\widetilde{\sigma}_{2}\in S_{2}\left(\Gamma \right)$ as follows:
\begin{itemize}
\item Let $p=\langle \langle k,A_0\rangle \rangle $. If $\widetilde{A}_0\ne \emptyset $, where
$$\widetilde{A}_0\defeq \left\{b\in A_0\ \middle|\ b\frown \langle 0,0,...\rangle \notin F_M\right\},$$
choose arbitrarily $b_0\in \widetilde{A}_0$. Otherwise, choose arbitrarily $b_0\in A_0$. Define $\widetilde{\sigma}_{2}\left(p\right)\defeq b_0$.
\item Let $p=\langle \langle k,A_0\rangle ,b_0,A_1,b_1,...,A_n\rangle \in \Gamma $. If $\widetilde{A}_n\ne \emptyset $, where
$$\widetilde{A}_n\defeq \left\{b\in A_n\ \middle|\ \langle\underset{\text{$k\cdot n$ times}}{\underbrace{0,0,...,0}}\rangle \frown b\frown \langle 0,0,...\rangle \notin F_M\right\},$$
choose arbitrarily $b_n\in \widetilde{A}_n$. Otherwise, choose arbitrarily $b_n\in A_n$. Define $\widetilde{\sigma}_{2}\left(p\right)\defeq b_n$.
\end{itemize}
In words, in every step, if Player 2 is offered finite sequences such that choosing any of them would lead to the constructed infinite sequence being in the complement of $F_M$, she chooses one of those sequences. Otherwise, she chooses an arbitrary finite sequence from the options provided. Let $\sigma_1\in S_1\left(\Gamma \right)$ be any strategy of Player I. Denote
$$\langle \langle k,A_0\rangle ,b_0,A_1,b_1,...\rangle\defeq \langle \sigma_1,\widetilde{\sigma}_{2}\rangle ,$$
and let $y\defeq \langle a_0,a_1,...\rangle \in \Dyacanopy$ be defined as in Equation~(\ref{projection}). If there exists $n\in \mathbb{N}$, such that $\widetilde{A}_n\ne \emptyset$, then $y\notin F_M$, and thus
$$p_{F_M}\left(\langle \sigma_1,\widetilde{\sigma}_{2}\rangle \right)=-1<\delta .$$
Otherwise, $\widetilde{A}_n=\emptyset$, for every $n\in \mathbb{N}$. Define $M_n\defeq M\cap \left[nk,nk+k-1\right]$, for every $n\in \mathbb{N}$. It follows that,
$$\left|A_n\right|\le 2^{\left|M_n\right|}.$$
Thus,
\begin{align*}
p_{F_M}\left(\langle \sigma_1,\widetilde{\sigma}_{2}\rangle \right)&=\liminf_{N\rightarrow\infty}\frac{1}{N}\cdot \sum_{n=0}^{N-1}\log_{2^k}\left(\left|A_n\right|\right)\\
&\le \liminf_{N\rightarrow\infty}\frac{1}{N}\cdot \sum_{n=0}^{N-1}\log_{2^k}\left(2^{\left|M_n\right|}\right)\\
&= \liminf_{N\rightarrow\infty}\frac{1}{N}\cdot \sum_{n=0}^{N-1}\frac{\left|M_n\right|}{k}\\
&= \liminf_{N\rightarrow\infty}\frac{\left|M\cap \left[0,kN-1\right]\right|}{kN}\\
&= \lim_{N\rightarrow\infty}\frac{\left|M\cap \left[0,N-1\right]\right|}{N}=\delta ,
\end{align*}
and thus $\widetilde{\sigma}_{2}$ guarantees $\delta $. Hence, $\val\left(G_{F_M}\right)=\delta $ and by Theorem~\ref{ItamarTheorem}, $\dim_{\mathcal{H}}\left(F_M\right)=\delta $.
\end{proof}

\section{Determining the Winner}
\label{Main} 
In this section we prove a general criterion for determining the winner in a win-lose alternating-move game, and apply it to prove Theorem~\ref{MainTheoremSmall}.
\subsection{Main Principle}
\label{MonotonicFunctions} 
In this subsection we present the main principle used in this paper for determining whether it is possible for Player I to guarantee a win in a win-lose alternating-move game $G=\left(T,W\right)$, and whether a winning strategy of Player II must exist. Our method consists of two steps.
\begin{enumerate}
\item Checking whether a strategy $s_I\in S_I\left(T\right)$ such that $\llbracket s_I\rrbracket\subseteq W$ may exist, \emph{i.e.}, whether a winning strategy of Player I may exist.
\item If not, and $W$ is a Borel set, or contained in a Borel set $Z\subseteq \canopy$ such that $\llbracket s_I\rrbracket\nsubseteq Z$, for every $s_I\in S_I\left(T\right)$, we may apply Theorem~\ref{Martin75} to deduce that Player II has a winning strategy.
\end{enumerate}
Refining these two steps yields the following general theorem.
\begin{theorem}
\label{MainMonotone} 
Let $T$ be a tree over an alphabet $\mathcal{A}$, and let $\xi:\powerset\left(\canopy\right)\to \left[0,\infty \right]$ be a monotone function, \emph{i.e.}, $\xi\left(C\right)\le \xi\left(Z\right)$ whenever $C\subseteq Z$. Let
$$t\defeq \inf_{s_I\in S_I\left(T\right)}\xi\left(\llbracket s_I\rrbracket \right).$$
If $W\subseteq \canopy$ such that $\xi\left(W\right)< t$, then Player I does not have a winning strategy in the game $\left(T,W\right)$.\\
Moreover, assume $W$ is a Borel set, or $\xi$ is Borel regular, \emph{i.e.}, for every $C\subseteq \canopy$, there exists a Borel set $Z\subseteq \canopy$ such that $C\subseteq Z$ and $\xi\left(C\right)=\xi\left(Z\right)$. Then Player II has a winning strategy in the game $\left(T,W\right)$.
\end{theorem}
Although the Hausdorff dimension is the focus of this paper, any monotone Borel regular function on the power set of the canopy of plays, which includes all Borel regular outer measures or dimensions, can be used as a criterion for determining the winner. For example, in Subsection~\ref{CorollariesMain} the $\tfrac{1}{2}$-dimensional Hausdorff outer measure is used as a criterion for determining the winner.
\begin{proof}[Proof of Theorem~\ref{MainMonotone}]
Let $W\subseteq \canopy$ such that $\xi\left(W\right)< t$. Suppose in contradiction that $\widetilde{s}_I\in S_I\left(T\right)$ is a winning strategy of Player I in the game $\left(T,W\right)$. Then $\llbracket \widetilde{s}_I\rrbracket\subseteq W$, and thus by the monotonicity of $\xi$, we get
$$\xi\left(\llbracket \widetilde{s}_I\rrbracket\right)\le \xi\left(W\right)< t,$$
a contradiction.\\
If $W$ is a Borel set, then by Theorem~\ref{Martin75}, the game $\left(T,W\right)$ is determined, and since Player I does not have a winning strategy in the game, Player II can guarantee a win in the game.\\
If $\xi$ is Borel regular, then there exists a Borel set $Z\subseteq \canopy$ such that $W\subseteq Z$, and 
$$\xi\left(Z\right)=\xi\left(W\right)<t.$$
Since $Z$ is a Borel set, the game $\left(T,Z\right)$ is determined. Since $\xi\left(Z\right)<t$, Player I does not have a winning strategy in the game, and thus there exists $\widetilde{s}_{II}\in S_{II}\left(T\right)$, a winning strategy of Player II in $\left(T,Z\right)$. Since $W\subseteq Z$, we have $\llbracket\widetilde{s}_{II}\rrbracket \subseteq Z^c\subseteq W^c$, and hence $\widetilde{s}_{II}$ is a winning strategy of Player II in $\left(T,W\right)$.
\end{proof}
\subsection{Proof of Theorem~\ref{MainTheoremSmall} Using Theorem~\ref{MainMonotone}}
\label{FirstProofMain} 
In this subsection we prove Theorem~\ref{MainTheoremSmall} by applying Theorem~\ref{MainMonotone}. We start by proving that $\dim_{\mathcal{H}}\left(\llbracket s_I\rrbracket \right)=\frac{1}{2}$, for every $s_I\in S_I\left(\Dya\right)$.

\begin{proposition}
\label{Strategy1/2} 
Let $s_I\in S_I\left(\Dya\right)$ be a pure strategy of Player I. Then $\dim_{\mathcal{H}}\left(\llbracket s_I\rrbracket \right)=\frac{1}{2}$. Similarly, let $s_{II}\in S_{II}\left(\Dya\right)$ be a pure strategy of Player II. Then $\dim_{\mathcal{H}}\left(\llbracket s_{II}\rrbracket \right)=\frac{1}{2}$.
\end{proposition}
\begin{proof}
We give the proof for strategies of Player I. The proof for strategies of Player II is analogous. Let
$$\llbracket \Dya^0 \rrbracket\defeq \left\{\langle a_0,a_1,...\rangle\in \Dyacanopy \ \middle|\ a_{2n}=0,\forall n\in \mathbb{N}\right\}\subseteq \Dyacanopy.$$
By Proposition~\ref{CalculateDimension}, $\dim_{\mathcal{H}}\left(\llbracket \Dya^0 \rrbracket\right)=\frac{1}{2}$. Proposition~\ref{Strategy1/2} follows since $\llbracket s_I\rrbracket$ is isometric to $\llbracket \Dya^0 \rrbracket$, for every $s_I\in S_I\left(\Dya\right)$. Indeed, define $\varphi: \llbracket s_I\rrbracket \to \llbracket \Dya^0 \rrbracket$ by $\varphi\left(\langle a_0,a_1,...\rangle \right)\defeq \langle b_0,b_1,...\rangle $, where
$$b_n\defeq \begin{cases}
    0, & n\in \mathbb{N}_{\text{even}},\\
    a_n, & n\in \mathbb{N}_{\text{odd}}.
\end{cases}$$
We prove that $\varphi$ is an isometry. First, we show that indeed, $\varphi$ is onto. Let $y=\langle b_0,b_1,...\rangle \in \llbracket \Dya^0 \rrbracket$. Define $s_{II}\in S_{II}\left(\Dya\right)$, a strategy of Player II by
$$s_{II}\left(p\right)=b_n,$$
for every $n\in \mathbb{N}_{\text{odd}}$, and $p\in \Dya$ such that $\len\left(p\right)=n$. Then $\varphi\left(\langle s_I,s_{II}\rangle \right)=y$.\\
We now show that $\varphi$ is distance-preserving, and, in particular, an injection. Let $x,x^\prime\in \llbracket s_I\rrbracket$ such that $x\ne x^\prime$. Denote $x=\langle a_0,a_1,...\rangle $, $x^\prime=\langle a^\prime_0,a^\prime_1,...\rangle $, $\varphi\left(x\right)=\langle b_0,b_1,...\rangle $, and $\varphi\left(x^\prime\right)=\langle b^\prime_0,b^\prime_1,...\rangle $. Let $m^*\in \mathbb{N}$ be the minimal $m\in \mathbb{N}$ such that $a_m\ne a^\prime_m$. Note that since $x,x^\prime\in \llbracket s_I\rrbracket$, it follows that $a_m=s_I\left(x_m\right)$, and $a^\prime_m=s_I\left(x^\prime_m\right)$ for every $m\in \mathbb{N}_{\text{even}}$, and thus $m^*\in \mathbb{N}_{\text{odd}}$. Therefore, by the definition of $\varphi$, it follows that $b_{m^*}\ne b^\prime_{m^*}$, and $b_m=b^\prime_m$, for every $0\le m < m^*$. This implies that $d\left(x,x^\prime\right)=d\left(\varphi\left(x\right),\varphi\left(x^\prime\right)\right)$, and thus $\varphi$ is indeed an isometry.
\end{proof}

We are now ready to prove Theorem~\ref{MainTheoremSmall}.
\begin{proof}[Proof of Theorem~\ref{MainTheoremSmall}]
By Proposition~\ref{Strategy1/2}, we have
$$\dim_{\mathcal{H}}\left(W\right)<\frac{1}{2}=\inf_{s_I\in S_I\left(\Dya\right)}\dim_{\mathcal{H}}\left(\llbracket s_I\rrbracket\right).$$
Since the Hausdorff dimension is a monotone Borel regular function, by Theorem~\ref{MainMonotone}, Player II can guarantee a win in the game $G=\left(\Dya,W\right)$.
\end{proof}

\subsection{Corollaries of Theorem~\ref{MainTheoremSmall}}
\label{CorollariesMain} 
In this subsection we present two corollaries of Theorem~\ref{MainTheoremSmall}.

The threshold given in Theorem~\ref{MainTheoremSmall} can be relaxed to a weaker form in terms of the $\frac{1}{2}$-dimensional Hausdorff measure. Direct computation shows that $\mathcal{H}^{\frac{1}{2}}\left(\llbracket \Dya^0 \rrbracket\right)=\sqrt{\frac{1}{2}}$.
Consequently, as in the proof of Proposition~\ref{Strategy1/2}, we deduce that $\mathcal{H}^{\frac{1}{2}}\left(\llbracket s_I \rrbracket\right)=\sqrt{\frac{1}{2}}$, for every $s_I\in S_I\left(\Dya\right)$. We therefore obtain the following corollary of Theorem~\ref{MainMonotone}.
\begin{corollary}
\label{HausdorffMeasureThreshold} 
Let $W\subseteq \Dyacanopy$ such that $\mathcal{H}^{\frac{1}{2}}\left(W\right)<\sqrt{\frac{1}{2}}$. Then Player II can guarantee a win in the game $G=\left(\Dya,W\right)$, and, in particular, $W$ is determined.
\end{corollary}

Exchanging the role of the two players we obtain the following corollary.
\begin{corollary}
Let $W\subseteq \Dyacanopy$ such that $\dim_{\mathcal{H}}\left(W^c\right)<\frac{1}{2}$. Then Player I can guarantee a win in the game $G=\left(\Dya,W\right)$, and, in particular, $W$ is determined.
\end{corollary}

\subsection{An Alternative Proof for Theorem~\ref{MainTheoremSmall}}
In this subsection we provide an alternative proof of Theorem~\ref{MainTheoremSmall}, using a probabilistic method.

Let $W\subseteq \Dyacanopy$ such that $\dim_{\mathcal{H}}\left(W\right)<\frac{1}{2}$. Theorem~\ref{MainTheoremSmall} guarantees the existence of a winning pure strategy of Player II in the game $G=\left(\Dya,W\right)$. However, the proof is non-constructive and does not offer an explicit winning strategy. The following theorem states that by flipping a fair coin in every step of the game and choosing her action accordingly, Player II wins with probability 1.\footnote{We thank Tom Meyerovitch for suggesting Theorem~\ref{FlipCoinTheorem}.} We will then show that this yields an alternative proof for Theorem~\ref{MainTheoremSmall}.

\begin{theorem}
\label{FlipCoinTheorem} 
Let $W\subseteq \Dyacanopy$ such that $\dim_{\mathcal{H}}\left(W\right)<\frac{1}{2}$. Let $s_I\in S_I\left(\Dya\right)$ be any pure strategy of Player I. Let $\sigma^U_{II}$ be the strategy of Player II defined by $\sigma^U_{II}\left(p\right)\sim \mathcal{U}\left(\left\{0,1\right\}\right)$, for every $p\in \Dya$ such that $\len\left(p\right)$ is odd. Let $\Prb \defeq \Prb_{s_I,\sigma^U_{II}}$ be the probability measure on $\Dyacanopy$ corresponding to the pair of strategies $\left(s_I,\sigma^U_{II}\right)$. Then $\Prb \left(W\right)=0$.
\end{theorem}
\begin{remark}
    By Theorem~\ref{FlipCoinTheorem}, the condition $\dim_{\mathcal{H}}\left(W\right)<\frac{1}{2}$ for $W\subseteq \Dyacanopy$ is sufficient for the uniform strategy of Player II to guarantee a win with probability $1$. However, this condition is not necessary. Indeed, let $M\subseteq\mathbb{N}_{\text{odd}}$ be an infinite set of density $0$, \emph{i.e.},
    $$\frac{\left|M\cap \left[0,n\right]\right|}{n+1}\xrightarrow[n\rightarrow\infty]{}0 .$$
    Let $F_{M^c}\subseteq \Dyacanopy$ be defined as in Definition~\ref{F_M}. By Proposition~\ref{CalculateDimension}, we have $\dim_{\mathcal{H}}\left(F_{M^c}\right)=1$. Nevertheless, since $M$ is infinite, if Player II plays uniformly at random at each stage, then with probability $1$ she chooses the digit $1$ at least once during stages from $M$, and therefore wins the game.
\end{remark}
\begin{lemma}
\label{Lemma1} 
Let $\emptyset \ne A\subseteq \Dyacanopy$ be a set such that $0<\diam\left(A\right)\le \frac{1}{2}$. Then there exists a position $p\in \Dya$ such that $\len\left(p\right)$ is even, $A\subseteq \left\llbracket \Dya_p\right\rrbracket$, and 
$$\emph{diam}\left(\left\llbracket \Dya_p\right\rrbracket\right)\le 2\cdot \emph{diam}\left(A\right).$$
\end{lemma}
\begin{proof}
Denote $n=-\log_2 \left(\diam\left(A\right)\right)$, let $x=\langle a_0,a_1,...\rangle \in A$, and let $q=\langle a_0,...,a_{n-1}\rangle $. Since $\diam\left(A\right)=\frac{1}{2^n}$, we have $A\subseteq \left\llbracket \Dya_q\right\rrbracket$, and
$$\diam\left(\left\llbracket \Dya_q\right\rrbracket\right)= \diam\left(A\right).$$
If $q$ is even, let $p=q$. If $q$ is odd, let $p=\langle a_0,...,a_{n-2}\rangle $. Then $A\subseteq \left\llbracket \Dya_q\right\rrbracket\subseteq \left\llbracket \Dya_p\right\rrbracket$, and
$$\diam\left(\left\llbracket \Dya_p\right\rrbracket\right)\le 2\cdot \diam\left(\left\llbracket \Dya_q\right\rrbracket\right) = 2\cdot \diam\left(A\right).$$
\end{proof}

\begin{lemma}
\label{Lemma2} 
Let $p\in \Dya$ such that $\len\left(p\right)$ is even. Then $\Prb \left(\left\llbracket \Dya_p\right\rrbracket\right)\le 2^{-\len\left(p\right)/2}$.
\end{lemma}
\begin{proof}
Denote $p=\langle a_0,...,a_{n-1}\rangle $. If there exists $k\in \mathbb{N}_{\text{even}}$ such that $0\le k\le n-1$ and $s_I\left(\langle a_0,...,a_{k-1}\rangle\right)\ne a_k$, then $\Prb \left(\left\llbracket \Dya_p\right\rrbracket\right) = 0$. Otherwise, Player II has chosen $\len\left(p\right)/2$ of the indices in $p$ uniformly and independently, and thus $\Prb \left(\left\llbracket \Dya_p\right\rrbracket\right)= 2^{-\len\left(p\right)/2}$.
\end{proof}

\begin{proof}[Proof of Theorem~\ref{FlipCoinTheorem}]
Let $\varepsilon\in \left(0,1\right)$ and let $\delta \in \left(\dim_{\mathcal{H}}\left(W\right),\frac{1}{2}\right)$. By the definition of the Hausdorff dimension, there exists a cover $\left(A_i\right)_{i\in \mathbb{N}}$ of $W$, such that 
$$\sum_{i\in \mathbb{N}}\left(\diam\left(A_i\right)\right)^\delta <\varepsilon.$$
By Lemma~\ref{Lemma1}, for every $i\in\mathbb{N}$ such that $0<\diam\left(A_i\right)\le \frac{1}{2}$ there exists $p_i\in\Dya$ such that $\len\left(p_i\right)$ is even, $A_i\subseteq \llbracket \Dya_{p_i}\rrbracket$, and 
$$\diam\left(\llbracket \Dya_{p_i}\rrbracket\right)\le 2\cdot \diam\left(A_i\right).$$
Let $I\subseteq\mathbb{N}$ be the set of $i\in \mathbb{N}$ such that $\left|A_i\right|=1$. Thus,
\begin{align*}
\Prb \left(W\right)&\le \Prb \left(\bigcup_{i\in \mathbb{N}}A_i\right)\\
&\le \sum_{i\in \mathbb{N}}\Prb \left(A_i\right)\\
&= \sum_{i\in \mathbb{N}\setminus I}\Prb \left(A_i\right)\\
&\le \sum_{i\in \mathbb{N}\setminus I}\Prb \left(\left\llbracket \Dya_{p_i}\right\rrbracket\right)\\
&\le \sum_{i\in \mathbb{N}\setminus I} 2^{-\len\left(p_i\right)/2}&&\text{Lemma~\ref{Lemma2}}\\
&=\sum_{i\in \mathbb{N}\setminus I} \left(2^{-\len\left(p_i\right)}\right)^{\frac{1}{2}}\\
&\le \sum_{i\in \mathbb{N}\setminus I} \left(2^{-\len\left(p_i\right)}\right)^\delta \\
&= \sum_{i\in \mathbb{N}\setminus I} \left(\diam\left(\left\llbracket \Dya_{p_i}\right\rrbracket\right)\right)^\delta\\
&\le 2^\delta\cdot \sum_{i\in \mathbb{N}\setminus I} \left(\diam\left(A_i\right)\right)^\delta\\
&=2^\delta\cdot \sum_{i\in \mathbb{N}} \left(\diam\left(A_i\right)\right)^\delta<2^\delta \cdot \varepsilon.
\end{align*}
Since $\varepsilon$ is arbitrary, $\Prb \left(W\right)=0$.
\end{proof}

\begin{proof}[Alternative proof for Theorem~\ref{MainTheoremSmall}]
Let $Z\subseteq \Dyacanopy$ be a Borel set given by Corollary~\ref{BorelSameDimension}, \emph{i.e.}, $W\subseteq Z$, and $\dim_{\mathcal{H}}\left(Z\right)=\dim_{\mathcal{H}}\left(W\right)<\frac{1}{2}$. Since $Z$ is a Borel set, by Theorem~\ref{Martin75}, the game $\left(\Dya,Z\right)$ is determined. We argue that Player I cannot guarantee a win in the game. Indeed, suppose by contradiction that $\widetilde{s}_I\in S_I\left(\Dya\right)$ is a winning strategy. Let $\sigma^U_{II}$ be as in Theorem~\ref{FlipCoinTheorem}, and let $\Prb =\Prb_{\widetilde{s}_I,\sigma^U_{II}}$ be the probability measure corresponding to the pair of strategies $\left(\widetilde{s}_I,\sigma^U_{II}\right)$. By Theorem~\ref{FlipCoinTheorem}, $\Prb \left(Z\right)=0$. However, since $\widetilde{s}_I$ is a winning strategy of Player I, $\llbracket \widetilde{s}_I\rrbracket \subseteq Z$, and thus
$$1=\Prb \left(\llbracket \widetilde{s}_I\rrbracket\right)\le \Prb \left(Z\right)=0,$$
a contradiction.\\
Therefore, Player II can guarantee a win in the game $\left(\Dya,Z\right)$. Let $\widetilde{s}_{II}\in S_{II}\left(\Dya\right)$ be a winning strategy in $\left(\Dya,Z\right)$. Since $W\subseteq Z$, we have $\llbracket\widetilde{s}_{II}\rrbracket \subseteq Z^c\subseteq W^c$, and hence $\widetilde{s}_{II}$ is a winning strategy in $G=\left(\Dya,W\right)$.
\end{proof}

\section{Tightness of Theorem~\ref{MainTheoremSmall} and Some Extensions}
\label{Ext} 
In this section we study the tightness of the conditions of Theorem~\ref{MainTheoremSmall}. As reflected in both proofs of Theorem~\ref{MainTheoremSmall}, the number $\frac{1}{2}$ serves as an important threshold for classifying a set as small enough to guarantee the existence of a winning pure strategy of Player II. In Subsection~\ref{LargeDimesion}, we show that once the Hausdorff dimension of a set is at least $\frac{1}{2}$, whether the game is determined, and, if it is, the identity of the winner, is not determined by the Hausdorff dimension. In Subsection~\ref{IncompleteSubsection} we show that the results in this paper may not hold if we replace the complete binary tree with a proper subtree. In Subsection~\ref{mAdicSubsection} we present a straightforward generalization for the results in the paper, when replacing $\Dya$ with the complete $m$-adic tree $\left\{0,1,...,m-1\right\}^{<\omega }$.

\subsection{Large Dimension Sets}
\label{LargeDimesion} 
In this subsection we show that once the Hausdorff dimension of a set is at least $\frac{1}{2}$, whether the game is determined, and, if it is, the identity of the winner, is not determined by the Hausdorff dimension.
\begin{theorem}
\label{LargeSets} 
Let $\delta \in \left[\frac{1}{2},1\right]$. There exist $\widetilde{W}^{II}\subset \widetilde{W}\subset \widetilde{W}^I\subseteq\Dyacanopy$, such that
\begin{itemize}
\item $\dim_{\mathcal{H}}\left(\widetilde W\right)=\dim_{\mathcal{H}}\left(\widetilde W^I\right)=\dim_{\mathcal{H}}\left(\widetilde W^{II}\right)=\delta $.
\item The set $\widetilde W$ is non-determined, Player I can guarantee a win in the win-lose alternating-move game $G_I=\left(\Dya,\widetilde W^I \right)$, and Player II can guarantee a win in the win-lose alternating-move game $G_{II}=\left(\Dya,\widetilde W^{II} \right)$.
\end{itemize}
\end{theorem}

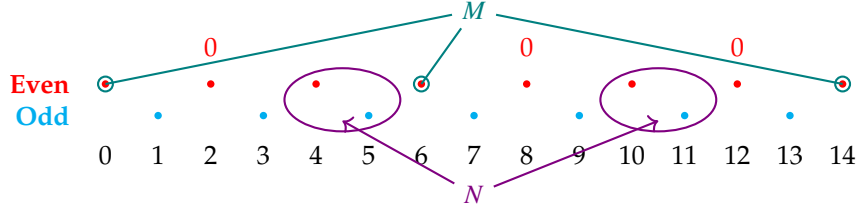
\begin{figure}
\centering
\begin{tikzpicture}[scale=0.7]
\def\yeven{0.6}
\def\yodd{0}
\foreach \x in {0,...,14} {
  \node[below] at (\x,-0.4) {\x};
}
\foreach \x in {1,3,5,7,9,11,13} {
  \fill[cyan] (\x,\yodd) circle (2pt);
}
\foreach \x in {0,2,4,6,8,10,12,14} {
  \fill[red] (\x,\yeven) circle (2pt);
}
\foreach \x in {0,6,14} {
  \draw[teal, thick] (\x,\yeven) circle (4pt);
}
\foreach \x in {2,8,12} {
\node[above=7pt, red] at (\x,\yeven) {$0$};
}
\draw[violet, thick] (4.5,0.3) ellipse (1.1cm and 0.6cm);
\draw[violet, thick] (10.5,0.3) ellipse (1.1cm and 0.6cm);
\node[teal] (M) at (7,2) {$M$};
\draw[teal, thick] (M) -- (0,\yeven);
\draw[teal, thick] (M) -- (6,\yeven);
\draw[teal, thick] (M) -- (14,\yeven);
\node[violet] (N) at (7,-1.5) {$N$};
\draw[violet, thick, ->] (N) -- (4.5,-0.1);
\draw[violet, thick, ->] (N) -- (10.5,-0.1);
\node[left,red] at (-0.5,\yeven) {\textbf{Even}};
\node[left,cyan] at (-0.5,\yodd) {\textbf{Odd}};
\end{tikzpicture}
\caption{A visualization of the set $W_\delta$: A play $\langle a_0,a_1,...\rangle \in W_\delta $ if all actions chosen by Player I at stages in $\mathbb{N}_{\text{even}}\setminus \left(M\cup N\right)$ are $0$, and all actions in stages in $N$ form a play in $W$.}
\label{WdFig}
\end{figure}

\begin{definition}[the set $W_\delta $]
\label{Wd} 
Let $W\subseteq \Dyacanopy$ be a set, such that $\langle 0,0,...\rangle \in W$, and let $\delta \in \left[\frac{1}{2},1\right]$. Let $N=\left\{n_0,n_1,...\right\}\subseteq \mathbb{N}$ be an infinite set of natural numbers that satisfy the following:
\begin{enumerate}
\item $n_{2\ell }\in \mathbb{N}_{\text{even}}$, and $n_{2\ell +1}=n_{2\ell }+1$, for every $\ell \in \mathbb{N}$.
\item The set $N$ has density $0$, \emph{i.e.},
$$\frac{\left|N\cap \left[0,k\right]\right|}{k+1}\xrightarrow[k\rightarrow\infty]{}0.$$
\end{enumerate}
Let $M\subseteq \mathbb{N}_{\text{even}}\setminus N$ be an infinite set of even natural numbers such that
$$\frac{\left|M\cap \left[0,k\right]\right|}{k+1}\xrightarrow[k\rightarrow\infty]{}\delta -\frac{1}{2}.$$
Define
\begin{align*}
    W_\delta&\defeq 
    \begin{aligned}[t]
    \Big\{
    \langle a_0,a_1,...\rangle \in \Dyacanopy\mid \text{$a_k=0$, if }&\text{$k\in \mathbb{N}_{\text{even}}\setminus \left(M\cup N\right)$,} \\
    & \text{and $\langle a_k\mid k\in N\rangle \in W$}
    \Big\}.
    \end{aligned}
\end{align*}
\end{definition}
The set $N$ has density $0$, and it contains pairs of consecutive numbers. At stages which are in $N$, the players can be thought of as playing the game $\left(\Dya,W\right)$. The set $M$ is a set of even non-negative integers of density $\delta -\frac{1}{2}$. It determines which coordinates outside $N$ are allowed to be $1$. The set $W_\delta $ is the set of binary sequences that ``look like'' a sequence from $W$ when we zoom in on the coordinates in $N$, and outside the set $N$, each coordinate in $\mathbb{N}_{\text{even}}\setminus M$ must have the value $0$. Since the possible values of coordinates in $M\cup \left(\mathbb{N}_{\text{odd}}\setminus N\right)$ of a play in $W_\delta $ are not limited, and since $M\cup \left(\mathbb{N}_{\text{odd}}\setminus N\right)$ is a set of density $\delta $, we will show, using Proposition~\ref{CalculateDimension}, that $\dim_{\mathcal{H}}\left(W_\delta \right)=\delta $. The set $W_\delta$ is displayed graphically in Figure~\ref{WdFig}.

\begin{remark}
\label{RemarkWd} 
$\text{ }$
\begin{enumerate}
\item The condition $\langle 0,0,...\rangle \in W$ can be changed to $W\ne \emptyset $. The play $\langle 0,0,...\rangle$ was chosen to match the definition of the set $F_M$.
\item Denote $\mathbb{N}\setminus N=\left\{\ell_0,\ell_1,...\right\}$. The game $G=\left(\Dya,W_\delta \right)$ can be described as follows. Player I and Player II are playing two win-lose alternating-move games, intermittently. On stages in $N$, the players perform moves in the game $G^{\prime} =\left(\Dya,W\right)$. On stages in $\mathbb{N}\setminus N$, the players perform moves in the game $G^{\prime\prime} =\left(\Dya,D\right)$, where
$$D\defeq \left\{\langle a_0,a_1,...\rangle \in \Dyacanopy\ \middle|\ \text{$a_n=0$, if $\ell_n \in \mathbb{N}_{\text{even}}\setminus M$}\right\}.$$
To win the game $G$, Player I (respectively, Player II) must win both games $G^{\prime}$ and $G^{\prime\prime}$ (respectively, at least one of those games). Since the games are independent, and since Player I can guarantee a win in $G^{\prime\prime}$ (by playing $a_n=0$ when $\ell_n \in \mathbb{N}_{\text{even}}\setminus M$), each player can guarantee a win in the game $G$ if and only if she can guarantee a win in the game $G^{\prime}$.
\item The game is essentially played in the stages that are in $N$, a set of zero density. The definition of the target set in stages in $\mathbb{N}\setminus N$ controls only the Hausdorff dimension of the set $W_\delta $. To construct an infinite sequence in $\mathbb{N}\setminus N$, there is freedom to choose each element, outside a set of density $1-\delta $. Therefore, $\dim_{\mathcal{H}}\left(W_\delta \right)=\delta $.
\end{enumerate}
\end{remark}
We now state these claims formally.
\begin{theorem}
\label{Wdequivalent} 
Let $\delta \in \left[\frac{1}{2},1\right]$, let $N,M\subseteq \mathbb{N}$, and let $W_\delta \subseteq \Dyacanopy$ as in Definition~\ref{Wd}. Then 
\begin{enumerate}
\item Player I has a winning strategy in the game $G=\left(\Dya,W_\delta \right)$ if and only if she has a winning strategy in the game $G^\prime=\left(\Dya,W\right)$.
\item Player II has a winning strategy in the game $G=\left(\Dya,W_\delta \right)$ if and only if she has a winning strategy in the game $G^\prime=\left(\Dya,W\right)$.
\end{enumerate}
\end{theorem}

\begin{theorem}
\label{WdDimd} 
Let $\delta \in \left[\frac{1}{2},1\right]$, let $N,M\subseteq \mathbb{N}$, and let $W_\delta \subseteq \Dyacanopy$ as in Definition~\ref{Wd}. Then $\dim_{\mathcal{H}}\left(W_\delta \right)=\delta $.
\end{theorem}

Items 2 and 3 of Remark~\ref{RemarkWd} give us the intuition for the proof of Theorems~\ref{Wdequivalent} and~\ref{WdDimd}, respectively. Formal proofs are given in Appendix~\ref{WdProofs}.

\begin{proof}[Proof of Theorem~\ref{LargeSets}]
Let $W$ be a non-determined set such that $\langle 0,0,...\rangle\in W $. Denote $\widetilde W=W_\delta $, $\widetilde W^I=\Dyacanopy_\delta $, and $\widetilde W^{II}=\left\{\langle 0,0,...\rangle\right\}_\delta$. By Theorems~\ref{Wdequivalent} and~\ref{WdDimd}, these three sets satisfy the desired requirements.
\end{proof}

\subsection{Incomplete Binary Tree}
\label{IncompleteSubsection} 
In this subsection we show that the results of Section~\ref{Main} do not necessarily hold if we replace the complete binary tree $\Dya$ with a proper subtree $T\subset \Dya$.

The following example shows that Player I can guarantee a win in a game with target set of Hausdorff dimension $0$, even when the Hausdorff dimension of the canopy is $1$. 
\begin{example}
\label{LargeTreeExample} 
Define
$$T\defeq \left\{\langle a_0,...,a_{n-1}\rangle\in \Dya \ \middle|\ \text{$a_0=1$, or $a_k=0$, for every $0\le k\le n-1$} \right\},$$
\emph{i.e.}, the tree with a canopy that consists exactly of the sequence $\langle 0,0,...\rangle $, and all the sequences that start with $1$. Note that $\dim_{\mathcal{H}}\left(\canopy\right)=1$.\\
Let $W=\left\{\langle 0,0,...\rangle \right\}\subseteq \canopy$. Even though $\dim_{\mathcal{H}}\left(W\right)=0$, Player I can guarantee a win in the game $G=\left(T,W\right)$ by choosing $0$ in every stage.
\end{example}
Applying Theorem~\ref{MainMonotone}, we obtain the following version of Theorem~\ref{MainTheoremSmall}.
\begin{theorem}
\label{IncompleteDimensionSmall} 
Let $T\subseteq \Dya$ be a subtree, and let $W\subseteq \canopy$ such that
$$\dim_{\mathcal{H}}\left(W\right)<\inf_{s_I\in S_I\left(T\right)}\dim_{\mathcal{H}}\left(\llbracket s_I\rrbracket \right).$$
Then Player II can guarantee a win in the win-lose alternating-move game $G=\left(T,W\right)$, and, in particular, $W$ is determined.
\end{theorem}
The proof is analogous to the proof of Theorem~\ref{MainTheoremSmall}. The following example shows that if $T\subset \Dya$, then the conclusion of Theorem~\ref{LargeSets} may fail. In particular, a non-determined set of dimension $\frac{1}{2}$ may not exist.
\begin{example}
\label{SmallTreeExample} 
Let $T$ be the subtree of $\Dya$ that contains all positions in which Player II always plays $0$. Formally,
$$T\defeq \left\{\langle a_0,...,a_{n-1}\rangle\in \Dya \ \middle|\ \text{$a_k=0$, for every $k\in \left[0,n-1\right]\cap \mathbb{N}_{\text{odd}}$} \right\}.$$
By Proposition~\ref{CalculateDimension}, $\dim_{\mathcal{H}}\left(\canopy\right)=\frac{1}{2}$. Let $W\subseteq \canopy$. Player I can guarantee a win in the game $G=\left(T,W\right)$ if and only if $W\ne \emptyset $. This is because Player II has a single strategy in the game. Thus, if $W\ne \emptyset $, the strategy of Player I that follows $x\in W$ is a winning strategy. In particular, if $\dim_{\mathcal{H}}\left(W\right)=\frac{1}{2}$, Player I can guarantee a win in the game $G$, and thus $W$ cannot be non-determined.
\end{example}
Example~\ref{SmallTreeExample} shows that if the dimension of the canopy of $T$ is at most $\frac{1}{2}$, one of the players might become a dummy player.

The following conjecture is an analogous version of Theorem~\ref{LargeSets} in the case where $\canopy$ is a strict subset of $\Dyacanopy$.
\begin{conjecture}
Let $T\subseteq \Dya$ be a subtree and suppose that $d=\dim_{\mathcal{H}}\left(\canopy\right)>\frac{1}{2}$. Let $\delta \in \left[\frac{1}{2},d\right]$. Then there exist $W,W^I,W^{II}\subseteq \canopy$, such that
\begin{itemize}
\item $\dim_{\mathcal{H}}\left(W\right)=\dim_{\mathcal{H}}\left(W^I\right)=\dim_{\mathcal{H}}\left(W^{II}\right)=\delta $.
\item The set $W$ is non-determined, Player I can guarantee a win in the win-lose alternating-move game $G_I=\left(T,W^I \right)$, and Player II can guarantee a win in the win-lose alternating-move game $G_{II}=\left(T,W^{II} \right)$.
\end{itemize}
\end{conjecture}

\subsection{\texorpdfstring{$m$}{Lg}-adic Trees}
\label{mAdicSubsection} 
In this subsection we extend the results of Sections~\ref{HDG} and~\ref{Main} to the situation where the underlying tree is not a complete binary tree but rather a \emph{complete $m$-adic tree}; namely, each player has $m$ actions in each position she controls.

Let $2\le m\in \mathbb{N}$ and let $T=\left\{0,1,...,m-1\right\}^{<\omega}$ be the complete $m$-adic tree. Define a metric on $\canopy$ as follows:
$$d_m\left(x,y\right)\defeq \begin{cases}
0, & x=y,\\
m^{-\min\left\{n\in \mathbb{N}\ \middle|\  a_n\ne b_n\right\}}, &x\ne y,
\end{cases}$$
where $x=\langle a_0,a_1,...\rangle $, and $y=\langle b_0,b_1,...\rangle $.\\
Denote by $\dim_{\mathcal{H}_m}$ the Hausdorff dimension corresponding to the metric $d_m$. The analogs of Theorems~\ref{MainTheoremSmall} and~\ref{LargeSets} are the following.
\begin{theorem}
\label{mDimensionSmall} 
Let $T=\left\{0,1,...,m-1\right\}^{<\omega}$. Let $A\subseteq \canopy$, such that $\dim_{\mathcal{H}_m}\left(A\right)<\frac{1}{2}$. Then Player II can guarantee a win in the win-lose alternating-move game $G=\left(T,A\right)$, and, in particular, the game is determined.
\end{theorem}
\begin{theorem}
\label{mDimensionLarge} 
Let $T=\left\{0,1,...,m-1\right\}^{<\omega}$, and let $\delta \in \left[\frac{1}{2},1\right]$. Then there exist $W,W^I,W^{II}\subseteq \canopy$, such that
\begin{itemize}
\item $\dim_{\mathcal{H}_m}\left(W\right)=\dim_{\mathcal{H}_m}\left(W^I\right)=\dim_{\mathcal{H}_m}\left(W^{II}\right)=\delta $.
\item The set $W$ is non-determined, Player I can guarantee a win in the win-lose alternating-move game $G_I=\left(T,W^I \right)$, and Player II can guarantee a win in the win-lose alternating-move game $G_{II}=\left(T,W^{II} \right)$.
\end{itemize}
\end{theorem}

The proofs of Theorems~\ref{mDimensionSmall} and~\ref{mDimensionLarge} are analogous to the proofs of Theorems~\ref{MainTheoremSmall} and~\ref{LargeSets}, respectively, and hence omitted. One needs to replace every occurrence of the number $2$ by the number $m$ and of the set $\left\{0,1\right\}$ by the set $\left\{0,1,...,m-1\right\}$ (and adapt the proofs accordingly).

\section{Alternating-Move Games on Complete Metric Spaces}
\label{SchmidtSection} 
In the previous sections, we considered win-lose alternating-move games played over trees, following the framework of Gale and Stewart \cite{GaleStewart1953} and Martin \cite{Martin75Paper}. In that setting, the sufficient criterion for classifying a set as small in Theorem~\ref{MainTheoremSmall} corresponds to the metric and topology induced by the tree. However, in some situations the game is already associated with a different topology and metric, and it is natural to ask whether there is an analogous criterion, for determining whether a set is small, corresponding to the given metric and topology.
\begin{example}[continues=CantorExample]
Recall $\mathcal{C}\subseteq \left[0,1\right]$, the Cantor-like set of numbers which have a base 5 expansion containing only the digits 1 and 3. We endow the interval $\left[0,1\right]$ with the standard metric $d$ inherited from $\mathbb{R}$. Consider the following win-lose alternating-move game, played in stages $n\in \mathbb{N}$:
\begin{itemize}
\item At $n=0$, Player I chooses an interval $B_0\in \left\{\left[0,\frac{1}{2}\right],\left[\frac{1}{2},1\right]\right\}$.
\item If $n$ is odd, and an interval $B_{n-1}=\left[a,b\right]$ has already been chosen, Player II chooses an interval $B_n\in \left\{\left[a,\frac{a+b}{2}\right],\left[\frac{a+b}{2},b\right]\right\}$.
\item If $n$ is even, and an interval $B_{n-1}=\left[a,b\right]$ has already been chosen, Player I chooses an interval $B_n\in \left\{\left[a,\frac{a+b}{2}\right],\left[\frac{a+b}{2},b\right]\right\}$.
\end{itemize}
In words, the players alternately bisect the current interval and choose one of the halves. Player I wins if $\bigcap_{n\in \mathbb{N}}B_n\subseteq \mathcal{C}$. We call this game the \emph{dyadic intervals intersection game with target set} $\mathcal{C}$. The Hausdorff dimension of $\mathcal{C}$ with respect to $d$ is $\dim_{\mathcal{H}_d}\left(\mathcal{C}\right)=\log_{5}\left(2\right)<\frac{1}{2}$. This raises a natural question: does this fact imply that Player II has a winning strategy in the game? The answer is positive, as the following proposition proves.
\end{example}
\begin{proposition}
\label{DyadicIntervalsProposition}
    Let $S\subseteq [0,1]$ such that $\dim_{\mathcal{H}}\left(S\right)<\frac{1}{2}$, where $\dim_{\mathcal{H}}$ denotes the Hausdorff dimension corresponding to the usual Euclidean metric on $\mathbb R$. Then, Player II can guarantee a win in the dyadic intervals intersection game with target set $S$.
\end{proposition}
We do not prove this proposition, as it is a special case of Corollary~\ref{madicTheorem}, which will be stated below. However, note that Proposition~\ref{DyadicIntervalsProposition} is almost identical to Theorem~\ref{MainTheoremSmall}, with two main differences: The tree of the game is not the complete binary tree, but the tree of nested dyadic intervals (in which each position of length $n$ represents a dyadic interval of length $2^{-n}$ inside the unit interval); and the monotone Borel regular function defined on the power set of the canopy of plays is not the Hausdorff dimension corresponding to the metric of the tree, but the Hausdorff dimension (corresponding to the Euclidean metric) of the projection of the tree's canopy on the unit interval.
\begin{remark}
\label{CantorExampleRemark}
The two sets defined in Example~\ref{CantorExample}, $W_\mathcal{C}\subseteq \Dyacanopy$ and $\mathcal{C}\subseteq \left[0,1\right]$, have the same Hausdorff dimension, corresponding to the metrics $d_2$ on $\Dyacanopy$ and the usual Euclidean metric $d$ on $\mathbb R$, respectively. This phenomenon is not unique to this example. Indeed, let $A\subseteq \left[0,1\right]$ be a subset of the unit interval, and let
$$W_A\defeq \left\{\langle a_0,a_1,...\rangle \in \Dyacanopy\ \middle|\ \sum_{n\in \mathbb{N}}\frac{a_n}{2^{n+1}}\in A\right\}.$$
Then $A$ and $W_A$ have the same Hausdorff dimension, corresponding to the metric $d_2$ on $\Dyacanopy$ and $d$ on $\mathbb R$, respectively (see, for example, Theorem 3.4 in \cite{Baek2009}). While the games $\left(\Dya,W_A\right)$ and the dyadic intervals intersection game with target set $A$ can be naturally identified (see Subsection \ref{SchmidtGamesSubsection}), the fact that they are played on nonisometric metric spaces makes the equality of the Hausdorff dimension of the two target sets nontrivial.
\end{remark}

In this section we present the formal definition of a Schmidt game and a generalization of Theorem~\ref{SchmidtHausdorff} to the setting of an arbitrary complete metric space. We also prove an analogous version of the theorem in the case of subgames of the $(\alpha,\beta)$-Schmidt games in a complete doubling metric space, and apply it to prove an analogous version of $m$-adic intersection games in $\mathbb{R}^d$, which Proposition \ref{DyadicIntervalsProposition} is a special case of.

\subsection{Schmidt Games as Win-Lose Alternating-Move Games}
\label{SchmidtGamesSubsection}
Given a mathematical space $X$, a Schmidt game on $X$ in the broad sense of \cite[Section 2]{Schmidt1966},\footnote{In Schmidt's original terminology, these are called $(\mathfrak{F},\mathfrak{G})$-games.} is a win-lose alternating-move game over a tree $T$, in which each position $p\in T$ represents a subset of $X$, and each extension $p\preceq q$ represents a subset of the one represented by $p$. Given a set $S\subseteq X$, the target set of the Schmidt game is the set of all plays which encode nested sequences of sets whose intersection is contained in $S$. In other words, the game goes as follows: the current position represents a subset of $X$, and the set of admissible actions is some collection of subsets of the current position. In each step, the player whose turn it is to play chooses one of the admissible sets. Player I wins if the intersection of all the chosen sets is contained in $S$.

Many topological games, including Banach-Mazur and Choquet games (see for example \cite{KechrisBook} or \cite{telgarsky}), are Schmidt games. These games transform abstract spatial properties into tangible strategic procedures, revealing both intuitions behind a property and methods for verifying it. In games like those of Banach-Mazur and Choquet, when one player has a winning strategy, this directly corresponds to a classical topological property, such as being a Baire space, Choquet space, or strongly Choquet space.

Although Schmidt alluded to the structure of win-lose alternating-move games as games over trees (see \cite[Section 13]{Schmidt1966}), he did not use such a formalism in his definition of the class of games nowadays known as Schmidt games. Since we find the connection between Schmidt games and games over trees to be useful, we follow the formal definition of Schmidt games written in the language of trees, which is due to Fishman, Ly, and Simmons (2014, \cite{FLS2014}).
\begin{definition}[Schmidt game]
    Given a space $X$, a target set $S\subseteq X$, and an alphabet $\mathcal{A}$, a \emph{Schmidt game on $X$ with target set $S$} is a pair $\left(\left\{(T^i,\varphi_i)\right\}_{i\in \mathcal{I}},S\right)$, where $\left\{T^i\right\}_{i\in \mathcal{I}}$ is a non-empty collection of trees over $\mathcal{A}$ and for every $i\in \mathcal{I}$ the function $\varphi_i\colon T^i\to\powerset(X)$ is a monotonic decreasing coding function satisfying $\varphi_i(p)\supseteq\varphi_i(q)$, for every $p,q\in T_i$ such that $p\preceq q$. The game is played as follows:
    \begin{itemize}
    \item Player II chooses a tree $T^i$.
    \item The players are playing the win-lose alternating-move game $\left(T^i,\Phi_i^{-1}\left[\powerset(S)\right]\right)$, where $\Phi_i\colon\llbracket T^i\rrbracket\to\powerset(X)$ is defined by
    \begin{align*}
        \Phi_i\left(\langle a_0,a_1,...\rangle\right)\defeq\bigcap_{n\in \mathbb{N}}\varphi_i\left(\langle a_0,...,a_n\rangle\right).
    \end{align*}
    \end{itemize}
\end{definition}
\begin{remark}
From a game-theoretic point of view, it may appear somewhat unusual that Player II makes the first move in Schmidt games. The move, however, could be interpreted not as a strategic move within the game itself, but rather as a choice of the game's setting, from a prescribed family of admissible settings. Equivalently, the game can be formulated as a standard win-lose alternating-move game, in which in step $0$ Player I has a single available move, in step $1$ Player II chooses the tree $T^i$, and then the Players proceed by playing the game $\left(T^i,\Phi_i^{-1}\left[\powerset(S)\right]\right)$. An analogous situation arises in tennis, where the player who does not serve first is the one who chooses which side of the court each player occupies in the first game of the match.
\end{remark}

On the one hand, one can always take $\mathcal{A}\subseteq\powerset(X)$ to be the set of subsets that the players can choose (at some point in the game), set $\varphi_i\left(\langle\text{ } \rangle \right)$ to be the first subset of $X$ chosen by Player II at the beginning of the game, and $\varphi_i\left(p \right)=a_n$, for every $p=\langle a_0,...,a_n\rangle\in T^i$, as is done in most of the literature. This is how, for example, the Banach-Mazur game is usually described, and how we treat the $(\alpha,\beta)$-Schmidt game in this paper. On the other hand, every win-lose alternating-move game $(T,W)$ can be seen as the Schmidt game $\left(\left\{\left(T,\varphi_{\llbracket T\rrbracket}\right)\right\},W\right)$ on the metric space $\left(\llbracket T\rrbracket,d_2\right)$ where the encoding function $\varphi_{\llbracket T\rrbracket}\colon T\to\powerset\left(\llbracket T\rrbracket\right)$ defined by
\begin{align*}
    \varphi_{\llbracket T\rrbracket}(p)\defeq \llbracket T_p\rrbracket.
\end{align*}
The advantage of formalizing Schmidt games using trees is that sometimes different Schmidt games, played on different spaces, can be described using the same tree (and therefore the same set of strategies), but with different coding functions.
\begin{example}[continues=CantorExample]
    The dyadic interval intersection game with target set $\mathcal{C}$ is a Schmidt game on the unit interval that can be formalized as $\left(\left\{\left(\Dya,\varphi_{[0,1]}\right)\right\},\mathcal C\right)$ with
    \begin{align*}
        \varphi_{[0,1]}\left(\langle a_0,...,a_{n-1}\rangle\right)\defeq
        \left[\frac{1}{2^{n}}\cdot \left(\sum_{i=0}^{n-1}2^{n-1-i}\cdot a_i\right),\frac{1}{2^n}\cdot \left(1+\sum_{i=0}^{n-1}2^{n-1-i}\cdot a_i\right)\right].
    \end{align*}
    The win-lose alternating-move game $(\Dya,W_\mathcal{C})$ can be seen as a Schmidt game $\left(\left\{\left(\Dya,\varphi_\Dya\right)\right\},W_\mathcal{C}\right)$, with
    \begin{align*}
        \varphi_\Dya(p)\defeq \llbracket T_p\rrbracket.
    \end{align*}
    In both cases, we have $W_\mathcal{C}=\Phi_{[0,1]}^{-1}\left[\mathcal \powerset\left(C\right)\right]=\Phi_\Dya^{-1}\left[\powerset\left(W_\mathcal{C}\right)\right]$, and thus these two Schmidt games are equivalent as win-lose alternating-move games over the tree $\Dya$, although they are Schmidt games on nonisometric metric spaces (see Remark \ref{CantorExampleRemark}).
\end{example}
The following example presents another case of games played on nonisometric metric spaces, which are equivalent as win-lose alternating-move games over the same tree.
\begin{example}[label=QuaternaryExample]
    Consider the alphabet $\mathcal{A}=\{0,1,2,3\}$ and the complete quaternary tree $T=\{0,1,2,3\}^{<\omega }$. We could define three different Schmidt games, each on a different metric space, played over the tree $T$.
    \begin{enumerate}
        \item A Schmidt game on the canopy of plays $(\llbracket T\rrbracket,d_2)$. In this case, we define as usual
        \begin{align*}
            \varphi_{\llbracket T\rrbracket}(p)\defeq\llbracket T_p\rrbracket.
        \end{align*}
        Alternatively, we could treat it as a Schmidt game on the canopy of the tree endowed with the $4$-adic metric, \emph{i.e.}, $(\llbracket T\rrbracket,d_4)$.
        \item A Schmidt game on the unit interval with the usual metric, \emph{i.e.}, $([0,1],d)$. In this case
        \begin{align*}
            \varphi_{[0,1]}\left(\langle a_0,...,a_{n-1}\rangle\right)\defeq
            \left[\frac{1}{4^n}\cdot \left(\sum_{i=0}^{n-1}4^{n-1-i}\cdot a_i\right),\frac{1}{4^n}\cdot \left(1+\sum_{i=0}^{n-1}4^{n-1-i}\cdot a_i\right)\right].
        \end{align*}
        That is, in this game, the players are choosing nested tetradic closed intervals inside the unit interval.
        \item A Schmidt game on the unit square with the $\ell_\infty$ norm, \emph{i.e.}, $([0,1]^2, \Vert\cdot\Vert_\infty)$. In this case, each chosen square is divided into $4$ two-dimensional dyadic closed squares, and the players are alternately choosing the next square, constructing a nested sequence of dyadic closed squares of shrinking size. Formally,
        \begin{align*}
        & \varphi_{[0,1]^2}\left(\langle a_0,...,a_{n-1}\rangle\right)\defeq\\
        & \left[\frac{1}{2^n}\cdot \left(\sum_{i=0}^{n-1}2^{n-1-i}\cdot(a_i\pmod 2)\right),\frac{1}{2^n}\cdot \left(\sum_{i=0}^{n-1}2^{n-1-i}\cdot(a_i\pmod 2)+1\right)\right]\times\\
        & \left[\frac{1}{2^n}\cdot \left(\sum_{i=0}^{n-1}2^{n-1-i}\cdot\lfloor\tfrac{a_i}{2}\rfloor\right),\frac{1}{2^n}\cdot \left(1+\sum_{i=0}^{n-1}2^{n-1-i}\cdot\lfloor\tfrac{a_i}{2}\rfloor\right)\right].
        \end{align*}
    \end{enumerate}
\end{example}
In \cite[Section 3]{Schmidt1966}, Schmidt describes some games on metric spaces that are special cases of the broad class we have formalized, in which the diameters of the sets chosen by the players decrease at a controlled rate. The $(\alpha,\beta)$-game, introduced there, is played by choosing nested closed balls whose radii are shrinking in each turn by $\alpha\in(0,1)$ or by $\beta\in(0,1)$, depending on the player. Target sets which Player I can always reach (called ``winning sets'' in this literature) possess a remarkable combination of properties: despite being small from measure-theoretic perspectives, they are stable under countable intersections, achieve full Hausdorff dimension, and since Player II chooses the first subset in the game, are also dense. This framework provides a powerful set of tools for Diophantine approximation (for instance, badly approximable vectors form a ``winning set''), enabling proofs that certain sets have full dimension while remaining closed under countable intersections and various geometric transformations.

The game's applications extend well beyond classical approximation into dynamical systems, through the Dani correspondence (see, for example, \cite{DimensionPaper}) and related techniques, as Schmidt-style arguments demonstrate that certain exceptional orbit collections (such as points whose orbits under toral endomorphisms avoid specific targets, or sets associated with bounded orbits in homogeneous dynamics) are winning sets with full dimension, and their intersections remain winning even within fractal structures. Subsequent refinements have further enhanced this framework's power and scope. For a review on different Schmidt games played on metric spaces the reader is referred to Badziahin, Harrap, Nesharim, and Simmons (2024, \cite{SchmidtGamesCantorWinning}).

In the following subsections we study two Schmidt games with a constant shrinking rate on complete metric spaces: the first is the $(\alpha,\beta)$-Schmidt game, and the second is a class of subgames of the $(\alpha,\beta)$-Schmidt game on Euclidean spaces, which generalizes the dyadic interval intersection game.
\subsection{\texorpdfstring{$(\alpha,\beta)$}{(alpha,beta)}-Schmidt Games}
We now turn to formally present the $(\alpha,\beta)$-Schmidt game from \cite[Section 3]{Schmidt1966}.

Given a complete metric space $(X,d)$, shrinking rates $\alpha,\beta\in(0,1)$, and a target set $S\subseteq X$, the \emph{$(\alpha,\beta)$-Schmidt game with target set $S$} is played as follows:
\begin{itemize}
\item At the beginning of the game, Player II chooses a closed ball $B_{-1}\subset X$ of radius $\rho_{-1}\in(0,\infty)$.
\item If $n$ is even, and a closed ball $B_{n-1}$ of radius $\rho_{n-1}$ has already been chosen, Player I chooses a closed ball $B_n\subset B_{n-1}$ of radius $\rho_n\defeq\rho_{n-1}\cdot \alpha$.
\item If $n$ is odd, and a closed ball $B_{n-1}$ of radius $\rho_{n-1}$ has already been chosen, Player II chooses a closed ball $B_n\subset B_{n-1}$ of radius $\rho_n\defeq\rho_{n-1}\cdot \beta$.
\end{itemize}
In words, the players alternately choose nested closed balls of shrinking radii, where the shrinking rates are $\alpha$ for Player I and $\beta$ for Player II. Player I wins if $\bigcap_{n\in \mathbb{N}}B_n\subseteq S$.
\begin{remark}
\label{AlphaBetaSchmidt} 
Formally, an $\left(\alpha ,\beta \right)$-Schmidt game on $\left(X,d\right)$ with target set $S$ can be viewed as the Schmidt game $\left(\left\{\left(T^i,\varphi_i\right)\right\}_{i\in \mathcal{I}},S\right)$, where $\mathcal{I}$ represents all closed balls with a positive radius in $X$, the empty sequence is coded $\varphi_i(\langle\rangle)=i$, a closed ball of radius $\rho_{-1}$, and the tree $T^i$ consists of sequences $\langle a_0,...,a_n\rangle $ of nested closed balls such that $a_n$ is of radius $\alpha^{\lceil\frac{n+1}{2}\rceil}\cdot \beta^{\lfloor\frac{n+1}{2}\rfloor}\cdot \rho_{-1}$ and $\varphi_i(\langle a_0,...,a_n\rangle)=a_n$.
\end{remark}
Note that in the $(\alpha,\beta)$-Schmidt game, the set of possible admissible actions of each player in each move is potentially infinite. This potentially makes the tree of the $(\alpha,\beta)$-Schmidt game much more complex than the complete binary tree, and it is not obvious, therefore, that the following version of Theorem~\ref{MainTheoremSmall} holds for the case of $(\alpha,\beta)$-Schmidt games as well.

\begin{theorem}
\label{MainSchmidt}
Let $(X,d)$ be a complete metric space, let $\dim_{\mathcal{H}}$ be the Hausdorff dimension with respect to the metric $d$, and let $\alpha,\beta \in \left(0,1\right)$. Suppose that there exist a closed ball $\widetilde{B}$ of radius $\widetilde{\rho}>0$ and $m\in\N\setminus\{0\}$ such that for every $\rho\in\left(0,\widetilde{\rho}\right)$, every closed ball $B\subseteq \widetilde{B}$ of radius $\rho$, contains at least $m$ closed balls with pairwise disjoint interiors of radius $\beta \cdot \rho$. Let $S\subseteq X$ such that 
$$\dim_{\mathcal{H}}\left(S\right)<\log_{\left(\alpha \beta\right)^{-1} }\left(m\right).$$
Then Player II has a winning strategy in the $\left(\alpha ,\beta \right)$-Schmidt game with target set $S$.
\end{theorem}
The proof of Theorem~\ref{MainSchmidt} appears in Appendix \ref{MainSchmidtProofSubsection}.\footnote{Note that we do not need the restriction $S\subseteq\widetilde{B}$ because Player II chooses the first ball.}

\subsection{Subgames of \texorpdfstring{$(\alpha,\beta)$}{(alpha,beta)}-Schmidt Games}
The game $G=\left(\Dya,W\right)$, where $W\subseteq \Dyacanopy$, studied in Section~\ref{Main}, can be viewed as a variant of the $\left(\frac{1}{2} ,\frac{1}{2} \right)$-Schmidt game with target set $W$, in which Player II is restricted to choose a closed ball of radius $1$ in stage $-1$ of the game. Theorem~\ref{MainTheoremSmall} implies that despite this restriction, Player II can still guarantee a win in the game whenever $\dim_{\mathcal{H}}\left(W\right)<\frac{1}{2}$.

In this subsection we study subgames of the $\left(\alpha,\beta\right)$-Schmidt games, in which at each stage, the players cannot necessarily choose every closed ball of a given radius, contained in the previously chosen closed ball, but can only choose a ball from a collection of admissible options. Even when restricting the players to a subgame, the sufficient condition for Player II to have a winning strategy in the game coincides with that in Theorem~\ref{MainSchmidt}, as long as the collection of admissible options is large enough.
\begin{example}[continues=QuaternaryExample]
The game $(\{0,1,2,3\}^{<\omega },\varphi_{[0,1]})$ is a $\left(\frac{1}{4},\frac{1}{4}\right)$-Schmidt game, in which Player II has a single choice in stage $-1$ of the game, the closed ball of radius $\frac{1}{2}$ around $\frac{1}{2}$. Then Player I can choose a closed ball of radius $\frac{1}{8}$ out of
$$\left\{\left[0,\frac{1}{4}\right],\left[\frac{1}{4},\frac{1}{2}\right],\left[\frac{1}{2},\frac{3}{4}\right],\left[\frac{3}{4},1\right]\right\},$$
and so on. Similarly, the game $(\{0,1,2,3\}^{<\omega },\varphi_{[0,1]^2})$, is a $\left(\frac{1}{2},\frac{1}{2}\right)$-Schmidt game, in which Player II has a single choice in stage $-1$ of the game, the closed ball of radius $\frac{1}{2}$ around the point $\left(\frac{1}{2},\frac{1}{2}\right)$. Then Player I can choose a closed ball of radius $\frac{1}{4}$ out of
$$\left\{\left[0,\frac{1}{2}\right]^2,\left[0,\frac{1}{2}\right]\times\left[\frac{1}{2},1\right],\left[\frac{1}{2},1\right]\times\left[0,\frac{1}{2}\right],\left[\frac{1}{2},1\right]^2\right\}$$
and so on.
\end{example}
We now formally define a subgame of an $(\alpha,\beta)$-Schmidt game, using our formalism of Schmidt games.
\begin{definition}[$(\alpha,\beta)$-Schmidt subgame]
Let $(X,d)$ be a complete metric space, $S\subseteq X$ a target set, $\alpha, \beta\in(0,1)$ shrinking rates, and $\left(\left\{\left(T^i,\varphi_i\right)\right\}_{i\in \mathcal{I}},S\right)$ an $(\alpha,\beta)$-Schmidt game on $(X,d)$. We say that $\left(\left\{(\widetilde T^i,\widetilde\varphi_i)\right\}_{i\in \widetilde{\mathcal{I}}},S\right)$ is a subgame of $\left(\left\{\left(T^i,\varphi_i\right)\right\}_{i\in \mathcal{I}},S\right)$ if $\emptyset\ne \widetilde{\mathcal I}\subseteq\mathcal I$, and for every $i\in\widetilde{\mathcal I}$ it holds that $\widetilde T^i\subseteq T^i$ is a subtree and $\widetilde\varphi_i=\varphi_i|_{\widetilde T^i}$.
\end{definition}
Less formally, given a complete metric space $(X,d)$, shrinking rates $\alpha,\beta\in(0,1)$, and a target set $S\subseteq X$, we can describe a subgame of an $(\alpha,\beta)$-Schmidt game with target set $S$ as follows:
\begin{itemize}
\item At the beginning of the game, Player II chooses a closed ball $B_{-1}\subset X$ of radius $\rho_{-1}\in(0,\infty)$ out of a non-empty collection $\mathcal{C}_{-1}$.
\item If $n$ is even, and a closed ball $B_{n-1}$ of radius $\rho_{n-1}$ has already been chosen, Player I chooses a closed ball $B_n\subset B_{n-1}$ of radius $\rho_n\defeq\alpha\cdot\rho_{n-1}$ out of a non-empty collection $\mathcal{C}\left(B_{-1},B_0,...,B_{n-1}\right)$.
\item If $n$ is odd, and a closed ball $B_{n-1}$ of radius $\rho_{n-1}$ has already been chosen, Player II chooses a closed ball $B_n\subset B_{n-1}$ of radius $\rho_n\defeq\beta\cdot\rho_{n-1}$ out of a non-empty collection $\mathcal{C}\left(B_{-1},B_0,...,B_{n-1}\right)$.
\end{itemize}
Player I wins if $\bigcap_{n\in \mathbb{N}}B_n\subseteq S$.

Thus, the games $(\{0,1,2,3\}^{<\omega },\varphi_{[0,1]})$ and $(\{0,1,2,3\}^{<\omega },\varphi_{[0,1]^2})$ from Example \ref{QuaternaryExample} are a $(\tfrac{1}{4},\tfrac{1}{4})$-Schmidt subgame and a $(\tfrac{1}{2},\tfrac{1}{2})$-Schmidt subgame, respectively.

The analogous result in the case of subgames of $\left(\alpha,\beta\right)$-Schmidt games is proved in the setting of a doubling metric space.
\begin{definition}[doubling metric space]
	\label{DoublingDefinition}
	A metric space $(X,d)$ is a \emph{doubling} metric space if there exists a \emph{doubling constant} $D\in\N\setminus\{0\}$ such that for every radius $r\in(0,\infty)$ and for every point $x\in X$, the closed ball $\closedb(x,r)$ can be covered by at most $D$ closed balls of radius $\tfrac{r}{2}$.
\end{definition}

Examples of doubling metric spaces include $\mathbb{R}^d$ with any norm, and an $m$-adic tree with the $d_2$ metric. Laakso spaces \cite{LaaksoSpace}, a class of fractal metric spaces, provide nontrivial examples of doubling metric spaces.

We are now ready to state the analogous version of Theorem~\ref{MainSchmidt} for $\left(\alpha,\beta\right)$-Schmidt subgames.
\begin{theorem}
\label{GeneralAlphaBeta} 
Let $(X,d)$ be a complete doubling metric space, let $\dim_{\mathcal{H}}$ be the Hausdorff dimension with respect to the metric $d$, and let $\alpha,\beta \in \left(0,1\right)$. Let $S\subseteq X$ and let $\left(\left\{(T^i,\varphi_i)\right\}_{i\in \mathcal{I}},S\right)$ be a subgame of the $\left(\alpha ,\beta \right)$-Schmidt game with target set $S$. Suppose that there exist $\widetilde{i}\in \mathcal{I}$ and integers $\widetilde{n}\in \mathbb{N}$ and $m\in\N\setminus\{0\}$ such that for every $n\ge \widetilde{n}$ and for every $p\in T^{\widetilde{i}}$ of length $2n+1$, there exist $a_0,...,a_{m-1}\in \mathcal{A}_p(T^{\widetilde i})$ such that 
$$\left\{\varphi_{\widetilde{i}}\left(p\frown \langle a_0\rangle\right),..., \varphi_{\widetilde{i}}\left(p\frown \langle a_{m-1}\rangle\right)\right\},$$
is a collection of $m$ closed balls with pairwise disjoint interiors. Suppose that
$$\dim_{\mathcal{H}}\left(S\right)<\log_{\left(\alpha \beta\right)^{-1} }\left(m\right).$$
Then Player II has a winning strategy in the $\left(\alpha ,\beta \right)$-Schmidt subgame $\left(\left\{(T^i,\varphi_i)\right\}_{i\in \mathcal{I}},S\right)$.
\end{theorem}
The proof of Theorem~\ref{GeneralAlphaBeta} appears in Appendix \ref{GeneralAlphaBetaProofSubsection}.

\begin{example}[continues=QuaternaryExample]
Applying Theorem~\ref{GeneralAlphaBeta} to the game $(\{0,1,2,3\}^{<\omega },\varphi_{[0,1]^2})$ with $\alpha =\beta =\frac{1}{2}$, we obtain that a sufficient condition for Player II to have a winning strategy in the game is that the Hausdorff dimension of the target set is strictly less than $\log_{4}\left(4\right)=1$. In contrast, applying the theorem to the game $(\{0,1,2,3\}^{<\omega },\varphi_{[0,1]})$ with $\alpha =\beta =\frac{1}{4}$, implies that the Hausdorff dimension of the target set must be strictly less than $\log_{16}\left(4\right)=\frac{1}{2}$ to guarantee the existence of a winning strategy for Player II. This is consistent with the fact that $\left[0,1\right]^2$ has Hausdorff dimension $2$, twice the Hausdorff dimension of $\left[0,1\right]$.
\end{example}

We conjecture that Theorem~\ref{GeneralAlphaBeta} remains valid for arbitrary complete metric spaces, without assuming the doubling property.

\subsection{\texorpdfstring{$m$}{m}-adic Cubes Intersection Games}
In this subsection we focus on a special case of $\left(\alpha ,\beta \right)$-Schmidt subgames played on the space $\mathbb{R}^d$, named \emph{$m$-adic cubes intersection games}, which generalize the dyadic intervals intersection game presented earlier in this section. In addition, we state the analogous version of Proposition~\ref{DyadicIntervalsProposition} that follows from Theorem~\ref{GeneralAlphaBeta}.

Dyadic intervals on the real line, and more generally, dyadic cubes in $\mathbb R^d$, are used in harmonic analysis, geometric measure theory, and more (see for example \cite{HarmonicAnalysisBook} or \cite{GMTBook}). For every level $\ell\in\mathbb{Z}$, the set of all half-open dyadic cubes of side-length $2^{-\ell}$, denoted
\begin{align*}
    \mathcal{D}^2_\ell\defeq \left\{2^{-\ell}\cdot\left([0,1)^d+b\right)\,|\, b\in\mathbb{Z}^d\right\},
\end{align*}
partitions the space, and the set of all half-open dyadic cubes of all levels, denoted by
\begin{align*}
    \mathcal{D}^2\defeq\bigcup_{\ell\in\mathbb Z}\mathcal{D}^2_\ell,
\end{align*}
forms a semiring of sets.\footnote{A \emph{semiring} of sets is a family of sets $\mathcal{S}$ such that: (1) $\emptyset\in\mathcal S$, (2) $A,B\in\mathcal S$ implies $A\cap B\in\mathcal S$, and (3) $A,B\in \mathcal{S}$ implies $A\setminus B=\bigcup_{i=1}^nC_i$ for some $n\in\N$ and disjoint $C_1,..., C_n\in\mathcal S$.} Taking the closures or the interiors of these dyadic cubes gives rise to alternative systems of dyadic cubes, such that each level forms a covering or a packing, respectively. Whichever system one chooses, the interiors of dyadic cubes are either disjoint or one is contained in the other, and the boundaries of dyadic cubes are of Hausdorff dimension $d-1$, and hence are negligible with respect to Lebesgue's measure. Furthermore, under the $\left\lVert\cdot \right\rVert_\infty$ norm, the closure and interior of a dyadic cube of diameter $2^{-\ell}$ correspond to a closed ball and an open ball of radius $2^{-(\ell+1)}$, respectively. In this case, when considering these cubes as balls, the dyadic closed cubes system corresponds to a \emph{splitting structure}, as defined by Badziahin and Harrap (2017, \cite{SplittingStructure}).\footnote{Note that in this case, the number of cubes partitioning their ``parent'' cube (the cube which is one level bigger and contains them) does not depend on the level nor on the specific cubes (it is uniform).} There is nothing unique about the choice of base 2, as one could consider $m$-adic cubes instead, for any $2\le m\in\N$, which give rise to different partitions of $\mathbb{R}^d$ with
\begin{align*}
    \mathcal{D}^m_\ell\defeq \left\{m^{-\ell}\cdot\left([0,1)^d+b\right)\,|\, b\in\mathbb{Z}^d\right\},
\end{align*}
and
\begin{align*}
    \mathcal{D}^m\defeq\bigcup_{\ell\in\mathbb Z}\mathcal{D}^m_\ell.
\end{align*}

Given $m\in \mathbb{N}\setminus \left\{0,1\right\}$ and a subset $S\subseteq \mathbb{R}^d$, the \emph{$m$-adic cubes intersection game with target set $S$} is played as follows:
\begin{itemize}
\item At the beginning of the game, Player II chooses a level $\ell_{-1}\in\mathbb Z$, and a closed cube $\overline{Q_{-1}}$ such that $Q_{-1}\in\mathcal{D}^m_{\ell_{-1}}$.
\item If $n$ is even, and a cube $Q_{n-1}\in\mathcal{D}^m_{\ell_{-1}+n}$ has already been chosen, Player I chooses a closed cube $\overline{Q_n}$ such that $Q_n\in\mathcal{D}^m_{\ell_{-1}+n+1}$ and $Q_n\subseteq Q_{n-1}$.
\item If $n$ is odd, and a cube $Q_{n-1}\in\mathcal{D}^m_{\ell_{-1}+n}$ has already been chosen, Player II chooses a closed cube $\overline{Q_n}$ such that $Q_n\in\mathcal{D}^m_{\ell_{-1}+n+1}$ and $Q_n\subseteq Q_{n-1}$.
\end{itemize}
In words, the players alternately choose a closed cube inside the last chosen cube, from the following level. Player I wins if $\bigcap_{n\in \mathbb{N}}\overline{Q_n}\subseteq S$.

The games from items 2 and 3 in Example \ref{QuaternaryExample} are ``subgames'' (since Player II is restricted to a single choice of a tree at the beginning of the game) of $m$-adic cubes intersection games: $\left(\{0,1,2,3\}^{<\omega },\varphi_{[0,1]^2}\right)$ would correspond to dyadic intersection games, while $\left(\{0,1,2,3\}^{<\omega },\varphi_{[0,1]}\right)$ would correspond to $4$-adic intersection games.

Since under the $\left\lVert\cdot \right\rVert_\infty$ norm, the closure of an $m$-adic cube from level $\ell \in \mathbb{Z}$ corresponds to a closed ball of radius $\frac{1}{2}\cdot m^{-\ell}$, the $m$-adic cubes intersection game is a $\left(\tfrac{1}{m},\tfrac{1}{m}\right)$-Schmidt subgame. Thus, applying Theorem~\ref{GeneralAlphaBeta} to the $m$-adic cubes intersection game yields the following corollary.
\begin{corollary}
\label{madicTheorem} 
Let $m\in \mathbb{N}\setminus \left\{0,1\right\}$ and let $S\subseteq \mathbb{R}^d$. Suppose that $\dim_{\mathcal{H}}\left(S\right)<\frac{d}{2}$. Then Player II has a winning strategy in the $m$-adic cubes intersection game with target set $S$.
\end{corollary}

\section{Hausdorff Dimension Games in Doubling Metric Spaces}
\label{GeneralHDG} 
As in the proof of Theorem~\ref{MainTheoremSmall}, the proof of Theorem \ref{GeneralAlphaBeta} uses a Hausdorff dimension game as a tool to obtain a lower bound for the Hausdorff dimension of the target set, when Player I can guarantee a win in the game. In this section, we introduce a new family of Hausdorff dimension games, and state a theorem which generalizes the result from \cite{DimensionPaper}.

\subsection{Imposed Subgames}
The dimension games we shall need to prove Theorem \ref{GeneralAlphaBeta} do not exactly generalize the dimension game invented and presented in \cite{DimensionPaper}. To make the relation between these games clearer, we present the definition of an \emph{imposed subgame}, a subgame in which only one of the players is restricted to a smaller set of possible moves.

\begin{definition}[imposed subtree and subgame]
    Let $T$ be a tree over the alphabet $\mathcal A$. We say that $T^{\prime}\subseteq T$ is \emph{imposing restrictions only on Player I}\footnote{Alternatively, some say that $T^{\prime}$ is an \emph{I-imposed} subtree of $T$.} if for every position $p\in T^{\prime}$ of odd length we have
    \begin{align*}
        \left\{a\in \mathcal{A}\ \middle|\  p\frown \langle a\rangle \in T^{\prime}\right\}=\left\{a\in \mathcal{A}\ \middle|\  p\frown \langle a\rangle \in T\right\}.
    \end{align*}
    In such a case, given a zero-sum payoff function $f\colon\canopy\to\R$, we say that $(T^{\prime},f|_{\llbracket T^{\prime}\rrbracket})$ is a subgame of $(T,f)$ \emph{imposing restrictions only on Player I}.
\end{definition}
If $T^{\prime}$ is an I-imposed subtree of $T$, then then passing from $(T,f)$ to $(T^{\prime},f|_{\llbracket T^{\prime}\rrbracket})$, restricts only the admissible moves of Player I, while Player II retains all admissible moves at every node reachable in $T^{\prime}$. The following holds:
\begin{proposition}
\label{prop:subgame_value}
    Let $(T,f)$ be a zero-sum alternating-move game, and let $T^{\prime}\subseteq T$ be a subtree with restrictions imposed only on Player I. Then,
    \begin{itemize}
        \item if $z$ is guaranteed by Player I in the subgame $(T^{\prime},f|_{\llbracket T^{\prime}\rrbracket})$, then it is also guaranteed by Player I in the game $(T,f)$;
        \item if $z$ is guaranteed by Player II in the game $(T,f)$, then it is also guaranteed by Player II in the subgame $(T^{\prime},f|_{\llbracket T^{\prime}\rrbracket})$;
        \item if both $(T,f)$ and $(T^{\prime},f|_{\llbracket T^{\prime}\rrbracket})$ are determined, then
        \begin{align*}
            \val(T^{\prime},f|_{\llbracket T^{\prime}\rrbracket})\le\val(T,f).
        \end{align*}
    \end{itemize}
\end{proposition}
Finally, note that every pure strategy $s_I\in S_I(T)$ of Player I corresponds to a unique subtree $T_{s_I}\subset T$ that imposes restrictions only on Player I and leaves her no freedom of choice, \emph{i.e.}, a subtree imposing restrictions only on Player I where every position of even length $n$ has only one extension of length $n+1$. In other words, forcing Player I to play according to a pure strategy $s_I\in S_I(T)$ is equivalent to letting her play on the subtree $T_{s_I}$, and indeed, $\llbracket T_{s_I}\rrbracket=\llbracket s_I\rrbracket$.

\subsection{The Non-Overlapping Dimension Game}
In this subsection we introduce the \emph{associated non-overlapping dimension game} for a given  $\beta$-shrinking Schmidt game, for $\beta \in \left(0,1\right)$, and state a generalization of Theorem 29.2 from \cite{DimensionPaper}. We begin with the definition of the game.

Let $(X,d)$ be a complete metric space, $\emptyset\ne S\subseteq X$ a target set, $\beta\in(0,1)$ a shrinking rate, and $\left(\left\{(T^i,\varphi_i)\right\}_{i\in \mathcal{I}},S\right)$ a $(\beta,\beta)$-Schmidt subgame on $(X,d)$, with a function $r:\mathcal I\to\left(0,\infty \right)$ defined such that $\varphi_i\left(\langle \text{ }\rangle \right)$ is a closed ball of radius $r\left(i\right)$. The \emph{non-overlapping Hausdorff dimension game associated to the $(\beta,\beta)$-Schmidt subgame} $\left(\left\{(T^i,\varphi_i)\right\}_{i\in \mathcal{I}},S\right)$ \emph{adjusted to} $r$ is played as follows:
\begin{itemize}
	\item On step 0, Player 1 chooses a tree $T^i$, a shrinking exponent $k\in\N\setminus\{0\}$, and a non-empty subset $\emptyset\neq A_0\subseteq\{p\in T^i\,|\,\len(p)=k\}$, such that the balls in $\{\varphi_i(p)\,|\,p\in A_0\}$ have pairwise disjoint interiors (move 0 of the game). Then, Player 2 chooses one of the positions $p_0\in A_0$ (move 1 of the game).
	\item At every step $n\in\N\setminus\{0\}$, Player 1 chooses a non-empty subset of positions $\emptyset\neq A_n\subseteq\{p\in T^i_{p_{n-1}}\,|\,\len(p)=\len(p_{n-1})+k\}$ such that the balls in $\{\varphi_i(p)\,|\,p\in A_n\}$ have pairwise disjoint interiors (move $2n$ of the game). Then, Player 2 chooses one of the positions $p_n\in A_n$ (move $2n+1$ of the game).
\end{itemize}
The payoff function is defined as: 
\begin{align*}
	\begin{split}
		p^{\mathcal{H}}_{_S}&\bigl(\langle i, k,\langle A_n,p_n\mid n\in\N\rangle\rangle  \bigr)\defeq\\
        &\begin{cases}   \liminf_{N\rightarrow\infty}\tfrac{1}{N}\sum_{n=0}^{N-1}\log_{\beta^{-k}}\left(\abs{A_n}\right), & \text{if }\bigcap_{n\in \mathbb{N}}\varphi_i\left(p_n\right)\subseteq S,\\
			-1, & \text{otherwise}.
		\end{cases}
	\end{split}
\end{align*}
Denote the tree of this dimension game by $\Gamma $, and by $\Gamma_r\subseteq\Gamma $ the subtree with the following additional restrictions only on Player 1: in each step $n\in\N$, Player 1 can choose only collections $A_n$ such that $\{\varphi_i(p)\,|\,p\in A_n\}$ is an $r_n$-separated\footnote{A set $A$ of subsets is \emph{$\rho$-separated} if $d(x,y)>\rho$ for every $X,Y\in A$ such that $X\ne Y$ and every $x\in X$ and $y\in Y$.} collection of closed balls of radius $r_n$, where $r_n\defeq r\left(i\right)\cdot \beta^{k\left(n+1\right)}$. We call $\left(\Gamma_r,p^\mathcal{H}_S|_{\llbracket\Gamma_r\rrbracket}\right)$ the \emph{separated Hausdorff dimension subgame associated with the $(\beta,\beta)$-Schmidt subgame $\left(\left\{(T^i,\varphi_i)\right\}_{i\in \mathcal{I}},S\right)$}.
\begin{remark}
The dimension game in \cite{DimensionPaper} can be seen as a game in which Player 1 chooses $\beta $, and then the players play a separated Hausdorff dimension subgame associated with the $(\beta,\beta)$-Schmidt game in which Player 1 is restricted to choose $k=1$.
\end{remark}

\begin{theorem}[upper bound on the value - non-overlapping dimension game]
\label{unrestricted_upper_bound}
Let $(X,d)$ be a complete doubling metric space with doubling constant $D$, and let $\emptyset\ne S\subseteq X$ be a target set. Let $\beta\in(0,1)$ and let $\left(\left\{(T^i,\varphi_i)\right\}_{i\in \mathcal{I}},S\right)$ be a $\left(\beta,\beta \right)$-Schmidt subgame on $(X,d)$. If Player 1 has a pure strategy that guarantees $\delta$ in the non-overlapping Hausdorff dimension game associated with $\left(\left\{(T^i,\varphi_i)\right\}_{i\in \mathcal{I}},S\right)$, then $\dim_{\mathcal{H}}\left(S\right)\ge \delta $.
\end{theorem}
The proof of Theorem~\ref{GeneralAlphaBeta} requires only Theorem~\ref{unrestricted_upper_bound} to obtain a lower bound on the Hausdorff dimension of the target set, when Player I can win in the Schmidt game. To show that the Hausdorff dimension of the target set is also a lower bound on the value of the separated dimension subgame, we require two additional regularity conditions on the game.

\begin{definition}[countable-covering,$L$-colored Schmidt game]
Let $(X,d)$ be a complete doubling metric space with doubling constant $D$, and let $\emptyset\ne S\subseteq X$ be a target set. Let $\beta\in(0,1)$ and let $\left(\left\{(T^i,\varphi_i)\right\}_{i\in \mathcal{I}},S\right)$ be a $(\beta,\beta)$-Schmidt subgame. 
\begin{enumerate}
    \item We say that the game is \emph{countable-covering} if there exists $\mathcal{I}^\prime\subseteq\mathcal{I}$ with $\abs{\mathcal{I}^\prime}\le \aleph_0$, such that
$$\bigcup_{i\in\mathcal{I}^\prime}\bigcup_{y\in\llbracket T^i\rrbracket}\Phi_i\left(y\right)=X.$$
\item We say that the game is \emph{$L$-colored} if there exists a uniform integer $L\in\mathbb{N}\setminus\left\{0\right\}$ such that for every $i\in\mathcal{I}$, every $k\in\mathbb{N}\setminus\{0\}$, every $n\in \mathbb{N}$, and every $p\in  T^i$ such that $\len\left(p\right)=kn$:
\begin{itemize}
\item There exists a finite collection
$$E^p \subseteq \left\{q\in T^i_p \ \middle|\ \len\left(q\right)=nk+k\right\}$$
such that
$$\bigcup_{q\in E^p} \bigcup_{z\in\llbracket T^i_q\rrbracket}\Phi_i\left(z\right) = \bigcup_{y\in\llbracket T^i_p\rrbracket}\Phi_i\left(y\right),$$
\emph{i.e.}, extending $p$ to one of the sequences in $E^p$ does not restrict the set of possible game outcomes.
\item There are sub-collections $E^p_1,...,E_L^p$ (not necessarily non-empty) such that $\biguplus_{\ell=1}^L E^p_\ell=E^p$, and $\left\{\varphi_i\left(q\right) \ \middle|\ q\in E^p_\ell\right\}$ is $r_n$-separated, for every $1\le \ell \le L$, where $r_n= \beta^{k(n+1)}\cdot r\left(i\right)$ is defined as in the definition of the non-overlapping Hausdorff dimension game associated to $\left(\left\{(T^i,\varphi_i)\right\}_{i\in \mathcal{I}},S\right)$ and adjusted to $r$.
\end{itemize}
\end{enumerate}
\end{definition}
\begin{remark}
Let $S\subseteq X$ be a target set and let $\delta <\dim_{\mathcal{H}}\left(S\right)$. The countable-covering condition implies that
$$S=\bigcup_{i\in\mathcal{I}^\prime}\left(\left(\bigcup_{y\in\llbracket T^i\rrbracket}\Phi_i\left(y\right)\right)\cap S\right).$$
By a property of the Hausdorff dimension, this yields that
$$\dim_{\mathcal{H}}\left(S\right)=\sup_{i\in\mathcal{I}^\prime}\dim_{\mathcal{H}}\left(\left(\bigcup_{y\in\llbracket T^i\rrbracket}\Phi_i\left(y\right)\right)\cap S\right),$$
and in particular there exists $i\in \mathcal{I}^{\prime}$ such that
$$\dim_{\mathcal{H}}\left(\left(\bigcup_{y\in\llbracket T^i\rrbracket}\Phi_i\left(y\right)\right)\cap S\right)>\delta.$$
The $L$-colored condition is necessary, as the collection of admissible balls at each stage may be highly irregular and cannot, in general, be controlled only by the geometric properties of the space $X$. Both $\left(\beta,\beta\right)$-Schmidt games and the $m$-adic cubes intersection games satisfy these two conditions.
\end{remark}
\begin{theorem}[lower bound on the value - separated dimension subgame]
\label{restricted_associated_lower} 
Let $(X,d)$ be a complete doubling metric space with doubling constant $D$, and let $\emptyset\ne S\subseteq X$ be a target set. Let $\beta\in(0,1)$ and let $\left(\left\{(T^i,\varphi_i)\right\}_{i\in \mathcal{I}},S\right)$ be a $\left(\beta,\beta \right)$-Schmidt subgame on $(X,d)$ that is countable-covering and $L$-colored. If $S$ is Borel, then for every $\delta\in[0,\dim_\mathcal{H}(S))$ Player 1 has a strategy that guarantees $\delta$ in $\left(\Gamma_r,p^\mathcal{H}_S|_{\llbracket\Gamma_r\rrbracket}\right)$, the separated Hausdorff dimension subgame associated with $\left(\left\{(T^i,\varphi_i)\right\}_{i\in \mathcal{I}},S\right)$.
\end{theorem}
Both the proof of Theorem~\ref{unrestricted_upper_bound} and the proof of Theorem~\ref{restricted_associated_lower} follow the constructions in the proof of \cite[Theorem 29.1]{DimensionPaper}, adjusted to regular $\left(\beta,\beta \right)$-Schmidt subgames. Once the countable-covering and the $L$-colored conditions are included, the proof of Theorem~\ref{restricted_associated_lower} is a straightforward adaptation of the proof of the upper bound on the Hausdorff dimension in Theorem 29.1 of \cite{DimensionPaper}, and is therefore omitted.

However, while in \cite{DimensionPaper} Player 1 is allowed to offer only $r_n$-separated sets of balls, in Theorem~\ref{unrestricted_upper_bound} Player 1 is allowed to offer a collection of admissible balls with pairwise disjoint interiors, according to the corresponding $\left(\beta,\beta \right)$-Schmidt subgame. Nevertheless, Theorem~\ref{unrestricted_upper_bound} implies that even with this additional freedom, Player 1 does not have a strategy that guarantees a payoff exceeding the Hausdorff dimension of the target set. The proof of Theorem~\ref{unrestricted_upper_bound} appears in Appendix~\ref{unrestricted_upper_bound_proof_subsection}.

Combining Theorems~\ref{unrestricted_upper_bound} and~\ref{restricted_associated_lower}, and using Proposition \ref{prop:subgame_value}, we obtain the following corollary.
\begin{corollary}[value of the Hausdorff dimension subgames]
\label{unrestricted_value}
Let $(X,d)$ be a complete doubling metric space with doubling constant $D$, and let $\emptyset\ne S\subseteq X$ be a Borel target set. Let $\beta\in(0,1)$ and $\left(\left\{(T^i,\varphi_i)\right\}_{i\in \mathcal{I}},S\right)$ be a $\left(\beta,\beta \right)$-Schmidt subgame on $(X,d)$ that is countable-covering and $L$-colored. Denote by $\Gamma $ the tree of the non-overlapping Hausdorff dimension game associated with $\left(\left\{(T^i,\varphi_i)\right\}_{i\in \mathcal{I}},S\right)$. Let $\Gamma_r\subseteq\Gamma^{\prime}\subseteq\Gamma $ be a subtree of $\Gamma $ with restrictions imposed only on Player 1. Then, 
\begin{align*}
    \val\left(\Gamma^{\prime},p_S^\mathcal{H}|_{\llbracket\Gamma^{\prime}\rrbracket}\right)=\dim_\mathcal{H}(S).
\end{align*}
\end{corollary}

\appendix
\section{Proofs from Subsection~\ref{LargeDimesion}}
\label{WdProofs} 
In this appendix we provide formal proofs of Theorems~\ref{Wdequivalent} and~\ref{WdDimd}, stated in Subsection~\ref{LargeDimesion}.
\subsection{Proof of Theorem~\ref{Wdequivalent}}
\begin{proof}[Proof of Theorem~\ref{Wdequivalent}]
We prove Item 1. The proof of Item 2 is analogous. Let $s_I\in S_I\left(\Dya\right)$ be a winning strategy of Player I in the game $G$. Note that $s_I$ is a winning strategy in $G$ if and only if $\llbracket s_I\rrbracket \subseteq W_\delta $. For every $x\defeq \langle a_0,a_1,...\rangle \in \llbracket s_I\rrbracket$, define $x^\prime=\langle a^\prime_0,a^\prime_1,...\rangle $ as follows. For every $n\in \mathbb{N}$,
$$a^\prime_n \defeq \begin{cases}
    0, & \text{if $n\in \mathbb{N}_{\text{even}}\setminus N$},\\
    a_n, & \text{otherwise}.
\end{cases}$$
Note that $\langle a_{n_k}\,|\,k\in\N\rangle=\langle a^\prime_{n_k}\,|\,k\in\N\rangle \in W$, since $x\in \llbracket s_I\rrbracket\subseteq W_\delta $. Therefore, we may assume without loss of generality that $s_I$ is a winning strategy that satisfies
$$s_I\left(\langle a_0,a_1,...,a_{n-1}\rangle\right)=0,$$
for every $n\in \mathbb{N}_{\text{even}}\setminus N$, and
$$s_I\left(\langle \text{ }\rangle \right)=0,$$
if $0\notin N$. Define $\widehat{s}_I\in S_I\left(\Dya\right)$, a strategy in the game $G^{\prime}$, as follows:
\begin{itemize}
\item $\widehat{s}_I\left(\langle \text{ }\rangle \right)\defeq s_I\left(\langle\underset{\text{$n_0$ times}}{\underbrace{0,0,...,0}}\rangle\right)$.
\item For every $q=\langle a_0,a_1,...,a_{m -1}\rangle \in \Dya$, such that $m \in \mathbb{N}_{\text{even}}$, let $p\defeq \langle b_0,b_1,...,b_{n_{m }-1}\rangle \in \Dya$ defined by
$$b_j\defeq \begin{cases}
    a_{k}, & \text{if $j=n_k$ for $0\le k\le m-1$},\\
    0, & \text{otherwise}.
\end{cases}$$
Define $\widehat{s}_I\left(q\right)\defeq s_I\left(p\right)$.
\end{itemize}
Suppose in contradiction that there exists $\widehat{s}_{II}\in S_{II}\left(\Dya\right)$, such that $\langle \widehat{s}_I,\widehat{s}_{II}\rangle\notin W$. Define $s_{II}\in S_{II}\left(\Dya\right)$, a strategy of Player II in the game $G$, as follows. For every $p=\langle a_0,...,a_{m-1}\rangle \in \Dya$, such that $m\in \mathbb{N}_{\text{odd}}$,
$$s_{II}\left(p\right)\defeq \begin{cases}
    \widehat{s}_{II}\left(\langle a_{n_0},...,a_{n_{k-1}}\rangle \right), & \text{if there exists $k\in \mathbb{N}$, such that $m=n_k$},\\
    0, &\text{otherwise}.
\end{cases}$$
Denote $\langle a_0,a_1,...\rangle\defeq \langle s_I,s_{II}\rangle$. Note that $\langle a_{n_k}\mid k\in \mathbb{N}\rangle= \langle \widehat{s}_I,\widehat{s}_{II}\rangle\notin W$. Therefore, $\langle s_I,s_{II}\rangle\notin W_\delta $, which implies that $s_I$ is not a winning strategy of Player I in $G$, a contradiction. Thus, Player I has a winning strategy in the game $G^\prime$.\\
In the other direction, let $s^{\prime}_I\in S_I\left(\Dya\right)$ be a winning strategy of Player I in the game $G^\prime$. Define a strategy $\widecheck{s}_I\in S_I\left(\Dya\right)$ of Player I in the game $G$, as follows. For every $p=\langle a_0,...,a_{m-1}\rangle \in \Dya$, such that $m\in \mathbb{N}_{\text{even}}$,
$$\widecheck{s}_I\left(p\right)\defeq \begin{cases}
    s^{\prime}_I\left(\langle a_{n_0},...,a_{n_{k-1}}\rangle \right), & \text{if there exists $k\in \mathbb{N}$, such that $m=n_k$},\\
    0, &\text{otherwise}.
\end{cases}$$
Suppose in contradiction that there exists $\widecheck{s}_{II}\in S_{II}\left(\Dya\right)$, such that $\langle \widecheck{s}_I,\widecheck{s}_{II}\rangle\notin W_\delta$. As before, we may assume without loss of generality that
$$\widecheck{s}_{II}\left(\langle a_0,a_1,...,a_{n-1}\rangle\right)=0,$$
for every $n\in \mathbb{N}_{\text{odd}}\setminus N$. Define $s^{\prime}_{II}\in S_{II}\left(\Dya\right)$, a strategy of Player II in the game $G^\prime$, as follows: for every $q=\langle a_0,a_1,...,a_{m -1}\rangle \in \Dya$, such that $m \in \mathbb{N}_{\text{odd}}$, let $p\defeq \langle b_0,b_1,...,b_{n_{m }-1}\rangle \in \Dya$ defined by
$$b_j\defeq \begin{cases}
    a_{k}, & \text{if there exists $0\le k\le m-1$ such that $j=n_k$},\\
    0, & \text{otherwise}.
\end{cases}$$
Define $s^{\prime}_{II}\left(q\right)\defeq \widecheck{s}_{II}\left(p\right)$. Denote $\langle a_0,a_1,...\rangle\defeq \langle \widecheck{s}_I,\widecheck{s}_{II}\rangle \notin W_\delta $. By the definition of $\widecheck{s}_I$ and $\widecheck{s}_{II}$, we know that $a_n=0$, for every $n\in \mathbb{N}\setminus N$. Thus, 
$$\langle s^{\prime}_I,s^{\prime}_{II} \rangle =\langle a_{n_k}\mid k\in \mathbb{N}\rangle \notin W,$$
which implies that $s^{\prime}_I$ is not a winning strategy of Player I in $G^\prime$, a contradiction. Thus, Player I has a winning strategy in the game $G$.
\end{proof}

\subsection{Proof of Theorem~\ref{WdDimd}}
\begin{proof}[Proof of Theorem~\ref{WdDimd}]
Define
$$\widehat{W}_\delta \defeq \left\{\langle a_0,a_1,...\rangle \in \Dyacanopy\ \middle|\ \text{$a_n=0$, if $n\in \mathbb{N}_{\text{even}}\setminus \left(M\cup N\right)$} \right\},$$
and
$$\widecheck{W}_\delta \defeq \left\{\langle a_0,a_1,...\rangle \in \Dyacanopy\ \middle|\ \text{$a_n=0$, if $n\in \left(\mathbb{N}_{\text{even}}\setminus M\right)\cup N$} \right\}.$$
Note that $\widehat{W}_\delta =F_K$, for $K=\mathbb{N}_{\text{odd}}\cup M\cup N$. Moreover,
\begin{align*}
    \frac{\left|K\cap \left[0,n\right]\right|}{n+1}&=\frac{\left|\mathbb{N}_{\text{odd}}\cap \left[0,n\right]\right|}{n+1}+\frac{\left|M\cap \left[0,n\right]\right|}{n+1}+\frac{\left|\left(\mathbb{N}_{\text{even}}\cap N\right)\cap \left[0,n\right]\right|}{n+1}\\
    & \xrightarrow[n\rightarrow\infty]{}\frac{1}{2}+\left(\delta -\frac{1}{2}\right)+0=\delta .
\end{align*}
Thus, by Proposition~\ref{CalculateDimension}, $\dim_{\mathcal{H}}\left(\widehat{W}_\delta \right)=\delta $. Similarly, $\widecheck{W}_{\delta}=F_L$, for $L=\left(\mathbb{N}_{\text{odd}}\cup M\right)\cap N^c$, and
\begin{align*}
\frac{\left|L\cap \left[0,n\right]\right|}{n+1} &=\frac{\left|\left(\mathbb{N}_{\text{odd}}\cup M\right)\cap \left[0,n\right]\right|}{n+1}-\frac{\left|\left(\mathbb{N}_{\text{odd}}\cup M\right)\cap N\cap \left[0,n\right]\right|}{n+1}\\
&=\frac{\left|\mathbb{N}_{\text{odd}}\cap \left[0,n\right]\right|}{n+1}+\frac{\left|M\cap \left[0,n\right]\right|}{n+1}-\frac{\left|\left(\mathbb{N}_{\text{odd}}\cup M\right)\cap N\cap \left[0,n\right]\right|}{n+1}\\
&\xrightarrow[n\rightarrow\infty]{}\frac{1}{2}+\left(\delta -\frac{1}{2}\right)-0=\delta .
\end{align*}
Thus, by Proposition~\ref{CalculateDimension}, $\dim_{\mathcal{H}}\left(\widecheck{W}_\delta \right)=\delta $.\\
Finally, by the definitions of $W_\delta $, $\widehat{W}_\delta $, and $\widecheck{W}_\delta $, and since $\langle 0,0,...\rangle \in W$,
$$\widecheck{W}_\delta \subseteq W_\delta  \subseteq \widehat{W}_\delta .$$
Thus, by monotonicity of the Hausdorff dimension,
$$\dim_{\mathcal{H}}\left(W_\delta \right)=\delta.$$
\end{proof}

\section{Proofs from Sections~\ref{SchmidtSection} and~\ref{GeneralHDG}}
\label{DimensionGamesProofsSection}
In this appendix we begin by proving Theorem~\ref{GeneralAlphaBeta} using Theorem~\ref{unrestricted_upper_bound}, and then prove Theorem~\ref{unrestricted_upper_bound} itself. Then, we prove Theorem~\ref{MainSchmidt}, whose proof consists of ideas from the proofs of Theorems~\ref{GeneralAlphaBeta} and~\ref{unrestricted_upper_bound}.

\subsection{Proof of Theorem~\ref{GeneralAlphaBeta}}
\label{GeneralAlphaBetaProofSubsection}
\begin{proof}[Proof of Theorem~\ref{GeneralAlphaBeta}]
Let $T$ be a tree defined as follows:
\begin{itemize}
\item In stage $0$, Player I has a single choice, $X$. This is a dummy stage defined so that Player I is formally the first to make a play in the game.
\item In stage $1$, Player II chooses a closed ball $B_{-1}$ of radius $\rho_{-1}\in \left(0,\infty \right)$ from the collection $\mathcal{C}_{-1}$.
\item If $2\le n$ is even, and a closed ball $B_{n-3}$ of radius $\rho_{n-3}$ has previously been chosen, Player I chooses in stage $n$ a closed ball $B_{n-2}\subset B_{n-3}$ of radius $\rho_{n-2}\defeq\rho_{n-3}\cdot \alpha$ from the collection $\mathcal{C}\left(B_{-1},B_0,...,B_{n-3}\right)$.
\item If $3\le n$ is odd, and a closed ball $B_{n-3}$ of radius $\rho_{n-3}$ has previously been chosen, Player II chooses in stage $n$ a closed ball $B_{n-2}\subset B_{n-3}$ of radius $\rho_{n-2}\defeq\rho_{n-3}\cdot \beta$ from the collection $\mathcal{C}\left(B_{-1},B_0,...,B_{n-3}\right)$.
\end{itemize}
Let $\psi\colon \llbracket T\rrbracket\to X$ be the projection of a play $\langle X,\langle B_n\mid n\in \mathbb{N}\cup \left\{-1\right\}\rangle\rangle $ on the element $x\in X$ such that $\{x\}=\bigcap_{n\in\N\cup \left\{-1\right\}}B_n$. Note that each player has a winning strategy in the general $\left(\alpha ,\beta \right)$-Schmidt subgame with target set $S$ if and only if she has a winning strategy in the game $\left(T,\psi^{-1}\left[S\right]\right)$, where $\psi^{-1}\left[F\right]\defeq \left\{y\in \canopy\ \middle|\ \psi \left(y\right)\in F\right\}$, for every $F\subseteq X$. Let $\xi:\powerset\left(\canopy\right)\to \left[0,\infty \right]$ defined by
$$\xi\left(W\right)\defeq \dim_{\mathcal{H}}\left(\psi\left[W\right]\right)= \dim_{\mathcal{H}}\left(\left\{\psi\left(y\right)\ \middle|\ y\in W\right\}\right),$$
for every $W\in \powerset\left(\canopy\right)$. The function $\xi$ is monotone by its definition. Moreover, $\xi$ is Borel regular. Indeed, let $W\in \powerset\left(\canopy\right)$, and let $F=\psi\left[W\right]$. By the Borel regularity of $\dim_{\mathcal{H}}$, there exists $B\subseteq X$ such that $F\subseteq B$, and $\dim_{\mathcal{H}}\left(B\right)=\dim_{\mathcal{H}}\left(F\right)$. Let $Z=\psi^{-1}\left[B\right]$ be a Borel set since $\psi$ is continuous. Then $W\subseteq Z$, and thus
$$\xi\left(W\right)\le \xi\left(Z\right)=\dim_{\mathcal{H}}\left(\psi\left[\psi^{-1}\left[B\right]\right]\right)\le \dim_{\mathcal{H}}\left(B\right)=\dim_{\mathcal{H}}\left(F\right)=\xi\left(W\right).$$
Since
$$\xi\left(\psi^{-1}\left(S\right)\right)\le \dim_{\mathcal{H}}\left(S\right),$$ 
and
$$\xi\left(\llbracket s_I\rrbracket\right)=\dim_{\mathcal{H}}\left(\psi\left[\llbracket s_I\rrbracket\right]\right),$$
for every $s_I\in S_I(T)$, by Theorem~\ref{MainMonotone}, if 
$$\dim_{\mathcal{H}}\left(S\right)<\inf_{s_I\in S_I(T)}\dim_{\mathcal{H}}\left(\psi\left[\llbracket s_I\rrbracket\right]\right),$$
then Player II has a winning strategy in the game. Thus, it is sufficient to prove that 
$$\dim_{\mathcal{H}}\left(\psi\left[\llbracket s_I\rrbracket\right]\right)\ge \log_{\left(\alpha \beta \right)^{-1}} \left(m\right)>\dim_{\mathcal{H}}\left(S\right),$$
for every $s_I\in S_I(T)$.\\
Let $s_I\in S_I(T)$ be a strategy of Player I, let $F\defeq \psi \left[\llbracket s_I\rrbracket\right]$. Let $\left(\left\{\left(\widehat{T}^i,\widehat{\varphi}_i\right)\right\}_{i\in \mathcal{I}},F\right)$ be a $\left(\alpha \cdot \beta ,\alpha \cdot \beta\right)$-Schmidt subgame with the target $F$ played as follows.
\begin{itemize}
\item Player $\widehat{\text{II}}$ chooses a closed ball $B_{-1}\in \mathcal{C}_{-1}$ of radius $r_{-1}\in \left(0,\infty \right)$.\footnote{In order to distinguish between the $\left(\alpha ,\beta \right)$-Schmidt subgame and the $\left(\alpha \cdot \beta ,\alpha \cdot \beta\right)$-Schmidt subgame, we call the players $\widehat{\text{I}}$ and $\widehat{\text{II}}$ in the description of this game.}
\item In stage $0$, Player $\widehat{\text{I}}$ temporarily chooses a closed ball $\widehat{B}_0\in \mathcal{C}\left(B_{-1}\right)$, and then chooses a closed ball $B_0$ such that $B_0\in \mathcal{C}\left(B_{-1},\widehat{B}_0\right)$ of radius  $r_0\defeq \alpha \cdot \beta \cdot r_{-1}$.
\item If $n$ is odd, and a closed ball $B_{n-1}$ of radius $r_{n-1}$ has previously been chosen, Player $\widehat{\text{II}}$ temporarily chooses in stage $n$ a closed ball $\widehat{B}_n\in \mathcal{C}\left(B_{-1},\widehat{B}_0,B_0,...,B_{n-1}\right)$, which corresponds to a move of Player I in the original $\left(\alpha,\beta\right)$-Schmidt subgame, and then chooses a closed ball $B_n\in \mathcal{C}\left(B_{-1},\widehat{B}_0,B_0,...,B_{n-1},\widehat{B}_n\right)$ of radius $r_n\defeq \alpha \cdot \beta \cdot r_{n-1}$, which corresponds to a move of Player II in the original $\left(\alpha,\beta\right)$-Schmidt subgame.
\item If $2\le n$ is even, and a closed ball $B_{n-1}$ of radius $r_{n-1}$ has previously been chosen, Player $\widehat{\text{I}}$ temporarily chooses in stage $n$ a closed ball $\widehat{B}_n\in \mathcal{C}\left(B_{-1},\widehat{B}_0,B_0,...,B_{n-1}\right)$, which corresponds to a move of Player I in the original $\left(\alpha,\beta\right)$-Schmidt subgame, and then chooses a closed ball $B_n\in \mathcal{C}\left(B_{-1},\widehat{B}_0,B_0,...,B_{n-1},\widehat{B}_n\right)$ of radius $r_n\defeq \alpha \cdot \beta \cdot r_{n-1}$, which corresponds to a move of Player II in the original $\left(\alpha,\beta\right)$-Schmidt subgame.
\item The functions $\left\{\widehat{\varphi}\right\}_{i\in \mathcal{I}}$ are defined as in Remark~\ref{AlphaBetaSchmidt}.
\end{itemize}
Consider the game $\left(\Gamma,p^\mathcal{H}_F\right)$, the non-overlapping Hausdorff dimension game associated to the $\left(\alpha \cdot \beta ,\alpha \cdot \beta\right)$-Schmidt subgame with the target $F$, adjusted to $r$, where $r\left(i\right)=r_{-1} $, if the choice of Player II at the beginning of the $\left(\alpha \cdot \beta ,\alpha \cdot \beta\right)$-Schmidt game is a closed ball of radius $r_{-1}$. Define $\widetilde{\sigma}_1\in S_1\left(\Gamma \right)$, a pure strategy of Player 1 in this dimension game, as follows:
\begin{itemize}
\item Let $\widetilde{i}\in \mathcal{I}$ and $\widetilde{n}\in \mathbb{N}$ such that for every $n\ge \widetilde{n}$ and for every $p\in T^{\widetilde{i}}$ of length $2n+1$, the collection
$$\left\{\varphi_{\widetilde{i}}\left(p\frown \langle a\rangle\right)\ \middle|\ a\in\mathcal{A}_p(T^{\widetilde i}) \right\},$$
contains a subcollection of $m$ closed balls with pairwise disjoint interiors. Player 1 chooses the tree $\widehat{T}_{\widetilde{i}}$ and $k=1$. Denote by $\closedb\left(x_{-1}, r_{-1}\right)\defeq \widehat{\varphi}_{\widetilde{i}}\left(\langle \text{ }\rangle \right)$.
\item Denote by $\closedb\left(x_0,\alpha\cdot r_{-1}\right)\defeq s_I\left(\langle X,\closedb\left(x_{-1}, r_{-1}\right) \rangle\right)$. Let
$$A_0\subseteq \mathcal{C}\left(\closedb\left(x_{-1}, r_{-1}\right),\closedb\left(x_0,\alpha\cdot r_{-1}\right)\right),$$
be a maximal collection of closed balls with pairwise disjoint interiors. Define $\widetilde{\sigma}_1\left(\langle\text{ } \rangle \right)\defeq \langle \widehat{T}^{\widetilde{i}}, 1, A_0\rangle $, and denote by $ r_0\defeq \beta\cdot \alpha\cdot r_{-1}$.
\item Let $p=\langle \langle \widehat{T}^i,k,A_0\rangle,b_0,...,A_{n-1},b_{n-1} \rangle\in \Gamma  $, where $1\le n\in \mathbb{N}$. If $k\ne 1$ or $\widehat{T}^i\ne \widehat{T}^{\widetilde{i}}$, define $\widetilde{\sigma}_1\left(p \right)$ arbitrarily. Otherwise, denote by
$$\closedb\left(x_{n},\alpha \cdot  r_{n-1}\right)\defeq s_I\left(\langle X,\closedb\left(x_{-1}, r_{-1}\right),\closedb\left(x_0, \alpha\cdot r_{-1}\right),b_0,...,\closedb\left(x_{n-1}, \alpha\cdot r_{n-2}\right),b_{n-1} \rangle\right).$$
Let
$$A_n\subseteq \mathcal{C}\left(\closedb\left(x_{-1}, r_{-1}\right),\closedb\left(x_0, \alpha\cdot r_{-1}\right),b_0,...,\closedb\left(x_{n-1}, \alpha\cdot r_{n-2}\right),b_{n-1},\closedb\left(x_{n},\alpha \cdot  r_{n-1}\right)\right),$$
be a maximal collection of closed balls with pairwise disjoint interiors. Define $\widetilde{\sigma}_1\left(p \right)\defeq A_n$, and denote by $ r_n\defeq \beta\cdot\alpha \cdot r_{n-1}$.
\end{itemize}
In words, at each stage, Player 1 chooses the ball that Player I would have chosen in the general $\left(\alpha ,\beta \right)$-Schmidt subgame, by following the strategy $s_I$, and then offers Player 2 a maximal collection of closed balls with pairwise disjoint interiors, out of the admissible responses that Player II could make to that choice. Let $\sigma_{2}\in S_{2}\left(\Gamma \right)$ be any pure strategy of Player 2. Denote
$$\langle \widetilde{\sigma}_1,\sigma_{2}\rangle =\langle \widehat{T}^i, k,\langle A_n,b_n\mid n\in\N\rangle\rangle=\langle \widehat{T}^{\widetilde{i}},1,\langle A_n,b_n\mid n\in\N\rangle\rangle.$$
By the definition of $F=\psi \left[\llbracket s_I\rrbracket\right]$, it follows that 
$$\bigcap_{n\in \mathbb{N}}b_n\subseteq F.$$
Furthermore, by the assumption, since Player 1 chose $\widetilde{i}\in \mathcal{I}$ in stage $0$ of the game, it follows that $\left|A_n\right|\ge m$, for every $\widetilde{n}\le n\in \mathbb{N}$. Thus,
\begin{align*}
p^{\mathcal{H}}_{F}\bigl(\langle \widehat{T}^{\widetilde{i}},1,&\langle A_n,b_n\mid n\in\N\rangle\rangle \bigr)=\\
&=\liminf_{N\rightarrow\infty}\tfrac{1}{N}\sum_{n=0}^{N-1}\log_{\left(\alpha \cdot \beta \right)^{-1}}\left(\left|A_n\right|\right)\\
&\ge \liminf_{N\rightarrow\infty}\tfrac{1}{N}\sum_{n=\widetilde{n}}^{N-1}\log_{\left(\alpha \cdot \beta \right)^{-1}}\left(\left|A_n\right|\right)\\
&\ge\liminf_{N\rightarrow\infty}\tfrac{1}{N}\sum_{k=\widetilde{n}}^{N-1}\log_{\left(\alpha \cdot \beta \right)^{-1}}\left(m\right)\\
&=\log_{\left(\alpha \cdot \beta \right)^{-1}}\left(m\right).
\end{align*}
Since $\sigma_{2}\in S_{2}\left(\Gamma \right)$ is arbitrary, Player 1 can guarantee $\log_{\left(\alpha \cdot \beta \right)^{-1}}\left(m\right)$, and thus by Theorem~\ref{unrestricted_upper_bound}
$$\dim_{\mathcal{H}}\left(\psi\left[\llbracket s_I\rrbracket\right]\right)=\dim_{\mathcal{H}}\left(F\right)\ge \log_{\left(\alpha \cdot \beta \right)^{-1}}\left(m\right),$$
completing the proof.
\end{proof}

\subsection{Proof of Theorem~\ref{unrestricted_upper_bound}}
\label{unrestricted_upper_bound_proof_subsection}
We now turn to prove Theorem~\ref{unrestricted_upper_bound}.\footnote{We thank Omri N. Solan for a helpful discussion which led to the proof of this theorem.} To prove it, we will need the following definition.
\begin{definition}
\label{packingNumberDefinition}
Let $(X,d)$ be a metric space, let $\emptyset\neq A\subseteq X$, and let $r\in(0,\infty)$.
\begin{itemize}
\item The \emph{external covering number of $A$ with closed balls of radius $r$}, denoted by $N^\text{ext}_r(A)$, is the smallest number of closed balls of radius $r$ that is needed to cover the set $A$ (by their union).
\item The \emph{packing number of $A$ with closed balls of radius $r$}, denoted by $N^{\pack}_r(A)$, is the maximal number of closed balls of radius $r$ and with pairwise disjoint interiors that can be packed inside the set $A$ (\emph{i.e.}, such that their union is a subset of $A$).
\end{itemize}
\end{definition}
The following lemma gives a bound on the number of closed balls of radius $r$ with pairwise disjoint interiors that can be packed in a closed ball of radius $3r$, in a doubling metric space.
\begin{lemma}
	\label{packingDoublingLemma}
	Let $(X,d)$ be a doubling metric space, with the doubling constant $D\in\N\setminus\{0\}$. Then for every $x\in X$ and $r\in(0,\infty)$ we have
	\begin{equation*}
		N^{\pack}_r\left(\closedb(x,3r)\right)\leq D^3.
	\end{equation*}
\end{lemma}
\begin{proof}
	Let $x\in X$ and $r\in(0,\infty)$. Since $(X,d)$ is doubling, we have
	\begin{equation*}
		N^\text{ext}_{3r/2}\left(\closedb(x,3r)\right),N^\text{ext}_{3r/4}\left(\closedb\left(x,\tfrac{3}{2}r\right)\right),N^\text{ext}_{3r/8}\left(\closedb\left(x,\tfrac{3}{4}r\right)\right)\leq D.
	\end{equation*}
	Thus, because this is true for every $x\in X$, we have
	\begin{equation}
		\label{CoveringNumberInequalityDoubling}
		N^\text{ext}_{3r/8}\left(\closedb(x,3r)\right)\leq D^3,
	\end{equation}
	since every closed ball of radius $\tfrac{3}{2}r$ that covers $\closedb(x,3r)$ can be covered by at most $D$ closed balls of radius $\tfrac{3}{4}r$, \emph{etc}.
	Suppose by contradiction that
	\begin{equation*}
		N^{\pack}_r\left(\closedb(x,3r)\right)>N^\text{ext}_{3r/8}\left(\closedb(x,3r)\right).
	\end{equation*}
	Let $y_1,...,y_k\in X$ such that
	\begin{equation*}
		k=N^\text{ext}_{3r/8}\left(\closedb(x,3r)\right),
	\end{equation*}
	and
	\begin{equation*}
		\closedb(x,3r)\subseteq\bigcup_{i=1}^k\closedb\left(y_i,\tfrac{3}{8}r\right).
	\end{equation*}
	Additionally, let $(x_i)_{i=1}^\ell$ be a sequence such that $\{\closedb(x_i,r)\}_{i=1}^\ell$ is a maximal packing set of closed balls of radius $r$ inside $\closedb(x,3r)$, with $\ell\in\N\cup\{\infty\}$. By the contradiction assumption, $\ell>k$.
	Thus, by the pigeonhole principle, there exist two centers $x_i\neq x_j$ such that both belong to the same covering ball, \emph{i.e.},
	\begin{equation*}
		x_i,x_j\in\closedb\left(y_1,\tfrac{3}{8}r\right),
	\end{equation*}
	where $y_1$ was chosen without loss of generality. Thus, by the triangle inequality
	\begin{equation*}
		d(x_i,x_j)\leq2\cdot\tfrac{3}{8}r<r,
	\end{equation*}
	a contradiction to the condition of non-intersection of the interiors of the packing closed balls. Hence,
	\begin{equation*}
		N^{\pack}_r\left(\closedb(x,3r)\right)\leq N^\text{ext}_{3r/8}\left(\closedb(x,3r)\right),
	\end{equation*}
	and by Inequality (\ref{CoveringNumberInequalityDoubling}),
	\begin{equation}
		N^{\pack}_r\left(\closedb(x,3r)\right)\leq D^3,
	\end{equation}
	for every $x\in X$ and $r\in(0,\infty)$.
\end{proof}

We will also use the following  lemma.
\begin{lemma}[Mass distribution principle]
\label{Massdistributionprinciple} 
Let $\left(X,d\right)$ be a complete doubling metric space, let $\mu $ be a Borel probability measure on $X$, and let $A\subseteq X$ be a Borel subset such that $\mu \left(A\right)>0$. Suppose there exist $\widetilde C,R,\delta >0$ such that 
$$\mu \left(\closedb\left(x,r\right)\cap A\right)\le \widetilde C\cdot r^\delta,$$
for every $x\in A$ and for every $r\in \left(0,R\right)$. Then $\dim_{\mathcal{H}}\left(A\right)\ge \delta $.
\end{lemma}
For a proof of Lemma~\ref{Massdistributionprinciple}, the reader is referred, for example, to \cite[Proposition 4.2]{FalconerBook}.

\begin{proof}[Proof of Theorem~\ref{unrestricted_upper_bound}]
	Since $\dim_{\mathcal{H}}\left(S\right)\geq 0$, there is nothing to prove in the case $\delta=0$. Suppose therefore that $\delta>0$. Let $\sigma^\delta_1$ be a pure strategy of Player 1 which guarantees $\delta$ in the non-overlapping Hausdorff dimension game associated to $\left(\left\{(T^i,\varphi_i)\right\}_{i\in \mathcal{I}},S\right)$. Denote
	\begin{equation*}
		\langle T^i,k,A_0\rangle \defeq\sigma_1^\delta(\emptyset),
	\end{equation*}
	where $k\in\N\setminus\{0\}$ is the shrinking exponent, and $A_0\subset T^i$ is a set of positions of length $k$. We shall from now on assume that $T^i$ and $k$ are given and omit them from the notation of game's tree and plays. Denote
    \begin{align*}
        M\defeq\max_{\rho >0,x\in X}N^{\pack}_{\beta^k\cdot\rho }\left(\closedb(x,\rho )\right)<\infty,
    \end{align*}
    which is well defined since $(X,d)$ is doubling, and denote
    \[
    \widetilde{M}=\log_{\beta^{-k}}\!\bigl(M\bigr).
    \]
    We use the behavior strategy \(\sigma^{U}_2\) of Player 2 consisting of uniformly randomly choosing a position from $ A_n$ at every position $\hist{n}$, \emph{i.e.},
    \begin{equation*}
    \sigma^{U}_2\bigl(\hist{n}\bigr)\sim \Unif( A_n).
    \end{equation*}
    The pair $(\sigma_1^\delta,\sigma_2^U)$ gives rise to the probability measure $\Prb_{\sigma_1^\delta,\sigma_2^U}$ on $\llbracket\Gamma\rrbracket$, whose support is closed. According to \cite[Corollary 3.31]{BellaicheThesis}, there exists a subtree $\widetilde\Gamma\subset\Gamma$ such that $\llbracket\widetilde\Gamma\rrbracket=\supp\Prb_{\sigma_1^\delta,\sigma_2^U}$, on which the two players play with probability 1 when they are following their strategies $\sigma_1^\delta$ and $\sigma_2^U$. Since $\sigma_2^U$ gives a positive probability for every possible extension at every odd position, we simply have that $\widetilde\Gamma=\Gamma_{\sigma_1^\delta}$, and thus $\llbracket\widetilde \Gamma\rrbracket=\llbracket\sigma_1^\delta\rrbracket$.\\
    Let $\varepsilon\in(0,\min\{\delta,\widetilde{M}\})$. For every play \(\seqAp\in\llbracket \sigma_1^{\delta}\rrbracket\) there exists a minimal
    \[
    n_{\varepsilon}^{\seqAp}\in\N,
    \]
    such that 
    \begin{equation*}
    \frac{1}{N+1}\sum_{n=0}^{N}\log_{\beta^{-k}}\left(\abs{A_n}\right)>\delta-\frac{\varepsilon}{2},
    \end{equation*}
    for every $N\geq n_{\varepsilon}^{\seqAp}$,
    since $\sigma_1^\delta$ guarantees $\delta$. Denote
    \begin{equation*}
    G(\ell)\defeq\Bigl\{\seqAp\in\llbracket \sigma_1^\delta\rrbracket \ \Big\lvert\ 
    n_{\varepsilon}^{\seqAp}\le\ell\Bigr\},
    \end{equation*}
    for every $\ell\in\N$, and note that $G(0)\subseteq G(1)\subseteq G(2)\subseteq...$
    Let \(n_{\varepsilon}\in\N\) be such that
    \begin{equation*}
    \Prb_{\sigma_1^{\delta},\,\sigma^{U}_2}(G(n_\varepsilon))\ge \tfrac{1}{2},
    \end{equation*}
    and denote $G\defeq G(n_\varepsilon)$.
    Note that for every \(n_2>n_1\ge n_{\varepsilon}\) we have that
    \begin{equation*}
    \Prb^{\cyl{\widetilde{\Gamma}_{\seqApn{n_1}}}}_{\sigma_1^{\delta},\,\sigma^U_2}\!\Bigl(\cyl{\widetilde{\Gamma}_{\seqApn{n_2}}}\Bigr)
    =\prod_{n=n_1+1}^{n_2}\frac{1}{\abs{ A_n}}\,,
    \end{equation*}
    where $\Prb^{F}_{\sigma_1^{\delta},\,\sigma^U_2}$ is the conditional probability given the realization of the event $F$.
    Denote
    \begin{equation*}
    n_{\varepsilon}^{*}=\left\lfloor\frac{2\widetilde{M}}{\varepsilon}\cdot( n_{\varepsilon}+1)\right\rfloor\, ,
    \end{equation*}
    and notice that \(n_{\varepsilon}^{*}> n_{\varepsilon}\) since \(\varepsilon<\widetilde{M}\). Thus, every $N\ge n_\varepsilon^*$ satisfies
    \begin{align*}
        N+1\ge n_\varepsilon^*+1>\frac{2\widetilde{M}}{\varepsilon}\cdot( n_{\varepsilon}+1),
    \end{align*}
    so that
    \begin{align}
    \label{upperbound_n_varepsilon}
        \frac{n_\varepsilon+1}{N+1}\cdot\widetilde M<\frac{\varepsilon}{2}.
    \end{align}
    Denote by $\psi:\llbracket \Gamma\rrbracket\to X$ the continuous but not necessarily injective projection of a play $\seqAp$ on the element $x\in X$ such that $\{x\}=\bigcap_{n\in\N}\varphi_i(p_n)$. Denote \(\mu=\Prb_{\sigma_1^{\delta},\,\sigma^U_2}\circ\psi^{-1}\) and \(E=\psi\left[G\right]\cap\supp\mu\), so that
    \begin{align*}
        \mu(E)=\mu(\psi\left[G\right])=\Prb_{\sigma_1^{\delta},\,\sigma^U_2}\left(\psi^{-1}\left(\psi\left[G\right]\right)\right)\geq\Prb_{\sigma_1^{\delta},\,\sigma^U_2}(G)\geq\tfrac{1}{2}.
    \end{align*}
    Let \(x\in E\).
    Note that for every \(N\ge n_{\varepsilon}\), the closed ball
    \(\closedb\left(x,r_N\right)\) can intersect no more than \(N^{\pack}_{r_N}\left(\closedb\left(x,3r_N\right)\right)\) balls with pairwise disjoint interiors and of radius $r_N$, since a center $y$ of a closed ball that intersects $\closedb\left(x,r_N\right)$ must satisfy $d(x,y)\leq 2r_N$, and therefore $\closedb\left(y,r_N\right)\subseteq\closedb\left(x,3r_N\right)$. For each $p\in A_N$, let $x_p\in X$ be such that $\closedb\left(x_p,r_N\right)=\varphi_i(p)$. Hence, by Lemma \ref{packingDoublingLemma}, \(\closedb\left(x,r_N\right)\) intersects no more than $D^3$ sets of the form $\varphi_i(p_N)$, where $\seqApn{N}\in\widetilde{\Gamma}$. This means that there are no more than $D^3$ histories of length \(2N+2\) with continuations in $G$ such that $\varphi_i(p_N)\cap\closedb\!\left(x,r_N\right)\neq\emptyset$.
    \noindent Therefore, denoting $$E^x_N=\left\{\seqApn{N}\,|\,\seqAp\in G\text{ and }\varphi_i(p_N)\cap\closedb(x,r_N)\neq\emptyset\right\},$$ we get
    \begin{align*}
    \mu\!\left(\closedb\left(x,r_N\right)\cap E\right)
    &\le&&
    \sum_{\seqApn{N}\in E^x_N}
    \Prb_{\sigma_1^{\delta},\,\sigma^G_2}\!\left(\cyl{\widetilde{\Gamma}_{\seqApn{N}}}\right)\\[3pt]
    &\le&
    D^3&\cdot \max_{\seqApn{N}\in E^x_N}
    \Prb_{\sigma_1^{\delta},\,\sigma^G_2}\!\left(\cyl{\widetilde{\Gamma}_{\seqApn{N}}}\right)\\[3pt]
    &=&
    D^3&\cdot \max_{\seqApn{N}\in E^x_N}
    \Prb_{\sigma_1^{\delta},\,\sigma^G_2}\!\left(\cyl{\widetilde{\Gamma}_{\seqApn{n_{\varepsilon}}}}\right)\\
    &&&\cdot
    \Prb^{\cyl{\widetilde{\Gamma}_{\seqApn{n_{\varepsilon}}}}}_{\sigma_1^{\delta},\,\sigma^G_2}\!\left(\cyl{\widetilde{\Gamma}_{\langle A_n,p_n\mid n_{\varepsilon}+1\le i\le N\rangle}}\right)\\[3pt]
    &=&
    D^3& \cdot \max_{\seqApn{N}\in E^x_N}
    \Prb_{\sigma_1^{\delta},\,\sigma^G_2}\!\left(\cyl{\widetilde{\Gamma}_{\seqApn{n_{\varepsilon}}}}\right)\\
    &&&\cdot
    \prod_{n=n_{\varepsilon}+1}^{N}\frac{1}{\abs{ A_n}}\\[3pt]
    &\le&
    D^3&\cdot \max_{\seqApn{N}\in E^x_N} \prod_{n=n_{\varepsilon}+1}^{N}\frac{1}{\abs{ A_n}} \,.
    \end{align*}
    
    We would like to find an upper bound on
    \begin{equation*}
    \max_{\seqApn{N}\in E^x_N} \prod_{i=n_{\varepsilon}+1}^{N}\frac{1}{\abs{A_n}},
    \end{equation*}
    of the type \(\widetilde C\cdot r_N^{\delta-\varepsilon}\), to use Lemma~\ref{Massdistributionprinciple}.
    We shall use the fact that for every \(N> n_{\varepsilon}^{*}\) and every \(\seqAp\in G\) we have
    \begin{align*}
    \delta-\frac{\varepsilon}{2}
    &< \frac{1}{N+1}\sum_{n=0}^{N}\log_{\beta^{-k}}\left(\abs{ A_n}\right)\\
    &= \frac{1}{N+1}\!\left(
    \frac{n_{\varepsilon}+1}{n_{\varepsilon}+1}\sum_{n=0}^{n_{\varepsilon}}\log_{\beta^{-k}}\left(\abs{ A_n}\right)
    +
    \frac{N-n_{\varepsilon}}{N-n_{\varepsilon}}\sum_{n=n_{\varepsilon}+1}^N\log_{\beta^{-k}}\left(\abs{ A_n}\right)
    \right)\\[3pt]
    &\le \frac{n_{\varepsilon}+1}{N+1}\,\widetilde{M}
    +\frac{N-n_{\varepsilon}}{N+1}\cdot\frac{1}{N-n_{\varepsilon}}\sum_{n=n_{\varepsilon}+1}^N\log_{\beta^{-k}}\left(\abs{ A_n}\right)\\[3pt]
    &\overset{\eqref{upperbound_n_varepsilon}}{\le} \frac{\varepsilon}{2}
    +\frac{N-n_{\varepsilon}}{N+1}\cdot\frac{1}{N-n_{\varepsilon}}\sum_{n=n_{\varepsilon}+1}^N\log_{\beta^{-k}}\left(\abs{ A_n}\right).
    \end{align*}
    Hence,
    \begin{align*}
    \delta-\varepsilon&\le\frac{N+1}{N-n_{\varepsilon}}\!\left(\delta-\frac{\varepsilon}{2}-\frac{\varepsilon}{2}\right)<\frac{1}{N-n_{\varepsilon}}\sum_{n=n_{\varepsilon}+1}^N\log_{\beta^{-k}}\left(\abs{ A_n}\right)\\
    &=\frac{1}{N-n_{\varepsilon}}\log_{\beta^{-k}}\!\left(\prod_{n=n_{\varepsilon}+1}^N\abs{ A_n}\right)=\frac{-\log_{\beta^{-k}}\!\left(\prod_{n=n_{\varepsilon}+1}^N\abs{ A_n}\right)}{-(N-n_\varepsilon)}\\
    &=\frac{\log_{\beta^{-k}}\!\left(\prod_{n=n_\varepsilon+1}^N\frac{1}{\abs{ A_n}}\right)}{\log_{\beta^{-k}}\left(\beta^{k(N-n_\varepsilon)}\right)}=\log_{\beta^{k(N-n_\varepsilon)}}\!\left(\prod_{n=n_{\varepsilon}+1}^N\frac{1}{\abs{ A_n}}\right).
    \end{align*}
    Therefore,
    \begin{equation*}
    \prod_{n=n_{\varepsilon}+1}^N\frac{1}{\abs{ A_n}}
    \le \left(\beta^{k(N-n_{\varepsilon})}\right)^{\delta-\varepsilon}
    = (\beta^{k(N+1)})^{\delta-\varepsilon}\cdot(\beta^{-k(n_{\varepsilon}+1)})^{\delta-\varepsilon}\, ,
    \end{equation*}
    and consequently,
    \begin{align*}
    \mu\!\left(\closedb(x,r_N)\cap E\right)
    &\le
    D^3\cdot \left(\frac{1}{\beta^{k(n_{\varepsilon}+1)}}\right)^{\delta-\varepsilon}
    \cdot (\beta^{k(N+1)})^{\delta-\varepsilon}\\
    &=D^3\cdot \left(\frac{1}{\beta^{k(n_{\varepsilon}+1)}r(i)}\right)^{\delta-\varepsilon}
    \cdot \left(\beta^{k(N+1)}r(i)\right)^{\delta-\varepsilon}\\
    &=D^3\cdot \left(\frac{1}{r_{n_\varepsilon}}\right)^{\delta-\varepsilon}
    \cdot r_N^{\delta-\varepsilon}.
    \end{align*}
    For every $r\in \left(0,r_{n^*_\varepsilon} \right)$, let $N\in \mathbb{N}$ be such that $r_{N+1}\le r< r_N$. Since $r_{N+1}=\beta^kr_N$, we have that $\beta^kr_N\le r$ and thus $r_N\le\beta^{-k}r$. Therefore,
    \begin{align*}
        \mu\!\left(\closedb(x,r)\cap E\right)
    &\le\mu\!\left(\closedb(x,r_N)\cap E\right)\le D^3\cdot \left(\frac{1}{r_{n_\varepsilon}}\right)^{\delta-\varepsilon}
    \cdot r_N^{\delta-\varepsilon}\\
    &\le D^3\cdot \left(\frac{1}{r_{n_\varepsilon}}\right)^{\delta-\varepsilon}
    \cdot \left(\beta^{-k}r\right)^{\delta-\varepsilon}\\
    &=D^3\cdot \left(\frac{1}{r_{n_\varepsilon+1}}\right)^{\delta-\varepsilon}
    \cdot r^{\delta-\varepsilon}.
    \end{align*}
    Thus, by Lemma~\ref{Massdistributionprinciple},
    \begin{align*}
    \dim_{\mathcal{H}}\left(E\right) \ge \delta-\varepsilon,
    \end{align*}
    and since $E\subseteq S$,
    \begin{align*}
    \dim_{\mathcal{H}}\left(S\right) \ge \delta-\varepsilon.
    \end{align*}
    Finally, since $\varepsilon>0$ is arbitrarily small,
    \begin{align*}
        \dim_{\mathcal{H}}\left(S\right) \ge \delta.
    \end{align*}
\end{proof}

\subsection{Proof of Theorem~\ref{MainSchmidt}}
\label{MainSchmidtProofSubsection}
In the proof of Theorem~\ref{unrestricted_upper_bound}, the doubling property of the space was used to bound the number of finite histories of a fixed length that can lead to the game reaching a given closed ball. As Theorem~\ref{MainSchmidt} is proved without a doubling assumption, we use a different technique, which consists mainly of fixing the center of the ball chosen by Player II once every few rounds, periodically.\footnote{We thank Tushar Das and David Simmons for suggesting this technique.} Fixing the center of the ball once every ``block,'' which is possible in the case of the ``full'' $(\alpha,\beta)$-Schmidt game, creates a separation that gives us a bound on the number of finite histories of a fixed length that can lead to the game reaching a given closed ball. To see why this is the case, we shall make use of the following lemma:
\begin{lemma}
\label{separating_lemma}
Let $(X,d)$ be a metric space, and let $\alpha,\beta\in(0,1)$. 
Let $B$ and $B'$ be two closed balls of radius $\rho>0$ with disjoint interiors.
Let $B_\alpha$ and $B'_\alpha$ be closed balls of radius $\alpha\cdot\rho$ such that
\[
B_\alpha\subseteq B,
\qquad
B'_\alpha\subseteq B'.
\]
Let $y$ and $y'$ be centers of $B_\alpha$ and $B'_\alpha$, respectively. Define
\[
B_1:=\closedb(y,\alpha\beta\rho),
\qquad
B'_1:=\closedb(y',\alpha\beta\rho).
\]
Then
\[
d(B_1,B'_1)\geq\alpha(1-\beta)\rho.
\]
\end{lemma}
\begin{proof}
Since
\[
B_\alpha=\closedb(y,\alpha\rho)\subseteq B
\]
and $\beta<1$, we have
\[
B_1=\closedb(y,\alpha\beta\rho)
\subseteq\openb(y,\alpha\rho)
\subseteq\operatorname{int}(B),
\]
where $\openb(y,\alpha\rho)$ denotes the open ball with center $y$ and radius $\alpha\cdot\rho$. Similarly,
\[
B'_1=\closedb(y',\alpha\beta\rho)
\subseteq\openb(y',\alpha\rho)
\subseteq\operatorname{int}(B').
\]
Take arbitrary points
\[
x\in B_1,
\qquad
x'\in B'_1,
\]
and suppose, toward a contradiction, that
\[
d(x,x')<\alpha(1-\beta)\rho.
\]
Since $x\in B_1=\closedb(y,\alpha\beta\rho)$, we have
\[
d(y,x)\leq\alpha\beta\rho.
\]
Therefore, by the triangle inequality,
\[
d(y,x')
\leq d(y,x)+d(x,x')
< \alpha\beta\rho+\alpha(1-\beta)\rho
=\alpha\rho.
\]
Hence
\[
x'\in\openb(y,\alpha\rho)\subseteq\operatorname{int}(B).
\]
On the other hand, since $x'\in B'_1$, we also have
\[
x'\in B'_1\subseteq\operatorname{int}(B').
\]
This contradicts
\[
\operatorname{int}(B)\cap \operatorname{int}(B')=\emptyset.
\]
Thus
\[
d(x,x')\geq\alpha(1-\beta)\rho
\]
for every $x\in B_1$ and $x'\in B'_1$. Taking the infimum over all such pairs gives
\[
d(B_1,B'_1)\geq\alpha(1-\beta)\rho.
\]
\end{proof}

As in the proof of Theorem~\ref{GeneralAlphaBeta}, we show that also in the setting of Theorem~\ref{MainSchmidt}, every strategy of Player I induces a subset of the metric space with Hausdorff dimension at least $$\log_{\left(\beta \cdot \alpha \right)^{-1}}\left(m\right).$$ Since without the doubling property we cannot make use of Theorem~\ref{unrestricted_upper_bound} directly, to obtain the desired lower bound on the Hausdorff dimension we follow the logic of the proof of Theorem~\ref{unrestricted_upper_bound} while exploiting the separating technique mentioned above, instead of the doubling property.
\begin{proof}[Proof of Theorem~\ref{MainSchmidt}]
Let $T$ be a tree defined as follows:
\begin{itemize}
\item In stage $0$, Player I has a single choice, $X$. This is a dummy stage defined so that Player I is formally the first to make a play in the game.
\item In stage $1$, Player II chooses a closed ball $B_{-1}$ of an arbitrary radius $\rho_{-1}\in \left(0,\infty \right)$.
\item If $2\le n$ is even, and a closed ball $B_{n-3}$ of radius $\rho_{n-3}$ has previously been chosen, Player I chooses in stage $n$ a closed ball $B_{n-2}\subset B_{n-3}$ of radius $\rho_{n-2}\defeq\rho_{n-3}\cdot \alpha$.
\item If $3\le n$ is odd, and a closed ball $B_{n-3}$ of radius $\rho_{n-3}$ has previously been chosen, Player II chooses in stage $n$ a closed ball $B_{n-2}\subset B_{n-3}$ of radius $\rho_{n-2}\defeq\rho_{n-3}\cdot \beta$.
\end{itemize}
Let $\psi\colon \llbracket T\rrbracket\to X$ be the projection of a play $\langle X,\langle B_n\mid n\in \mathbb{N}\cup \left\{-1\right\}\rangle\rangle $ on the element $x\in X$ such that $\{x\}=\bigcap_{n\in\N\cup \left\{-1\right\}}B_n$. As in the proof of Theorem~\ref{GeneralAlphaBeta}, each player has a winning strategy in the general $\left(\alpha ,\beta \right)$-Schmidt subgame with target set $S$ if and only if she has a winning strategy in the game $\left(T,\psi^{-1}\left[S\right]\right)$. Thus, it is sufficient to prove that 
$$\dim_{\mathcal{H}}\left(\psi\left[\llbracket s_I\rrbracket\right]\right)\ge \log_{\left(\alpha \beta \right)^{-1}} \left(m\right)>\dim_{\mathcal{H}}\left(S\right),$$
for every $s_I\in S_I(T)$.\\
Let $s_I\in S_I(T)$ be a strategy of Player I and let $\ell\in\N\setminus\{0\} $ be the first positive integer such that 
$$\left(\alpha \cdot \beta \right)^\ell<\alpha \cdot  \left(1-\beta \right).$$
Let $\ell<L\in \mathbb{N}$, and let $\widetilde{T}\subseteq \left\{0,1,...,m-1\right\}^{<\omega }$ be the subtree consisting of sequences $\langle a_0,...,a_{n-1}\rangle\in \left\{0,1,...,m-1\right\}^{<\omega }$ such that $a_k=0$, for every
$$k=L-\ell,L-\ell+1,...L-1\pmod L.$$
Define collections of closed balls as follows:
\begin{itemize}
\item Let $B^{\langle \text{ }\rangle }$ be a closed ball of radius $\rho_{-1}$ such that for every $\rho\in\left(0,\rho_{-1}\right)$, every closed ball $B\subseteq B^{\langle \text{ }\rangle }$ of radius $\rho$ contains at least $m$ closed balls with pairwise disjoint interiors of radius $\beta \cdot \rho$. Such a closed ball $B^{\langle \text{ }\rangle }$ exists by the assumption of the theorem. Let $A_{-1}\defeq \left\{B^{\langle \text{ }\rangle }\right\}$.
\item Denote by $B_{\alpha}^{\langle\text{ } \rangle }\defeq s_I\left(\langle X,B^{\langle \text{ }\rangle }\rangle \right)$ the closed ball of radius $\alpha\cdot \rho_{-1}$ that Player I chooses according to $S_I$ in response to the choice of $B^{\langle\text{ } \rangle }$. By the assumption, there exist $B^{\langle 0\rangle},...,B^{\langle m-1\rangle }\subseteq B_{\alpha}^{\langle\text{ } \rangle }$, closed balls of radius $\rho_0\defeq\alpha \cdot \beta \cdot \rho_{-1}$ with pairwise disjoint interiors. Denote
$$A_0\defeq\left\{B^{\langle 0\rangle},...,B^{\langle m-1\rangle }\right\}.$$
\item Let $1\le n\in \mathbb{N}$ and suppose that the collections $A_0,...,A_{n-1}$ have been defined. For every $B^{\langle b_0,...,b_{n-1}\rangle }\in A_{n-1}$, denote by
$$B_\alpha^{\langle b_0,...,b_{n-1}\rangle }\defeq s_I\left(\langle X,B^{\langle \text{ }\rangle },B_\alpha^{\langle \text{ }\rangle },B^{\langle b_0\rangle },B_\alpha^{\langle b_0\rangle },...,B^{\langle b_0,...,b_{n-1}\rangle}\rangle \right),$$
the closed ball of radius $\alpha\cdot \rho_{n-1}$ that Player I chooses according to $S_I$ given the history $\langle X,B^{\langle \text{ }\rangle },B_\alpha^{\langle \text{ }\rangle },B^{\langle b_0\rangle },B_\alpha^{\langle b_0\rangle },...,B^{\langle b_0,...,b_{n-1}\rangle}\rangle$.
\begin{itemize}
\item If 
$$n=0,1,...,L-\ell -1\pmod L,$$
define
$$A_n\defeq\bigcup_{B^{\langle b_0,...,b_{n-1}\rangle }\in A_{n-1}}\left\{B^{\langle b_0,...,b_{n-1},0\rangle},...,B^{\langle b_0,...,b_{n-1},m-1\rangle }\right\},$$
where $B^{\langle b_0,...,b_{n-1},0\rangle},...,B^{\langle b_0,...,b_{n-1},m-1\rangle }\subseteq B_\alpha^{\langle b_0,...,b_{n-1}\rangle }$ are $m$ closed balls of radius $\rho_n\defeq\alpha \cdot \beta \cdot \rho_{n-1}$ with pairwise disjoint interiors, which exist for every $B^{\langle b_0,...,b_{n-1}\rangle }\in A_{n-1}$ by the assumption. 
\item Else, if
$$ n = L-\ell\pmod L,$$
denote by $x_\alpha^{\langle b_0,...,b_{n-1}\rangle }$ a center of the ball $B_\alpha^{\langle b_0,...,b_{n-1}\rangle }$, \emph{i.e.},
$$B_\alpha^{\langle b_0,...,b_{n-1}\rangle }=\closedb\left(x_\alpha^{\langle b_0,...,b_{n-1}\rangle },\alpha\cdot \rho_{n-1}\right).$$
Let $B^{\langle b_0,...,b_{n-1},0\rangle}\defeq \closedb\left(x_\alpha^{\langle b_0,...,b_{n-1}\rangle },\rho_n\right)$, where $\rho_n\defeq\alpha \cdot \beta \cdot \rho_{n-1}$, and set
$$A_n\defeq\bigcup_{B^{\langle b_0,...,b_{n-1}\rangle }\in A_{n-1}}\left\{B^{\langle b_0,...,b_{n-1},0\rangle}\right\}.$$
Note that by Lemma \ref{separating_lemma}, $d(B,B')\geq\alpha(1-\beta)\rho_{n-1}$ for every $B\in A_n$ and $B'\in A_n\setminus\{B\}$.
\item Else, if $\ell>1$ and
$$n=L-\ell+1,...,L-1\pmod L,$$
choose a unique closed ball $B^{\langle b_0,...,b_{n-1},0\rangle}\subseteq B_\alpha^{\langle b_0,...,b_{n-1}\rangle}$ of radius $\rho_n\defeq\alpha \cdot \beta \cdot \rho_{n-1}$ and set
$$A_n\defeq\bigcup_{B^{\langle b_0,...,b_{n-1}\rangle }\in A_{n-1}}\left\{B^{\langle b_0,...,b_{n-1},0\rangle}\right\}.$$
Note that since we chose a unique ball for every history in the last few stages, we still have that $d(B,B')\geq\alpha(1-\beta)\rho_{\underline{n}}$ for every $B\in A_n$ and $B'\in A_n\setminus\{B\}$, where $\underline{n}\defeq \left\lceil\tfrac{n}{L}\right\rceil\cdot L-\ell-1$.
\end{itemize}
\end{itemize}
In words, the construction proceeds in blocks of $L$ stages. During the first $L-\ell$ stages of each block, for every closed ball of radius $\rho_{n-1}$ in $A_{n-1}$, we choose, according to $s_I$ a closed ball of radius $\alpha \cdot \rho_{n-1}$ contained in it, and then fix inside this ball $m$ closed balls of radius $\rho_n\defeq\alpha \cdot \beta \cdot \rho_{n-1}$ with pairwise disjoint interiors. At the next stage, namely when $n=L-\ell\pmod L$, for every closed ball of radius $\rho_{n-1}$ in $A_{n-1}$, we again choose, using $s_I$, a closed ball of radius $\alpha \cdot \rho_{n-1}$ contained in it, but then fix inside it the closed ball of radius $\rho_n\defeq\alpha \cdot \beta \cdot \rho_{n-1}$ with the same center. If $\ell>1$, then stages of the block, the construction proceeds without branching, \emph{i.e.}, after applying $s_I$, we choose one arbitrary closed ball of radius $\rho_n\defeq\alpha \cdot \beta \cdot \rho_{n-1}$ contained in it.

As we mentioned above, Lemma \ref{separating_lemma} implies that for every $n\in \mathbb{N}$, every $B\in A_{n\cdot L-1}$, and every $B'\in A_{n\cdot L-1}\setminus\{B\}$ we have that $d(B,B')\geq\alpha(1-\beta)\rho_{n\cdot L-\ell-1}$. Thus, it holds that for every $B\in A_{n\cdot L -1}$, $x\in B,$ and $r\in \left(0,\alpha(1-\beta)\rho_{{n\cdot L-\ell-1}}\right)$,
\begin{equation}
    \label{choiceofell}
    \left|\left\{B^{\prime}\in A_{n\cdot L-1}\ \middle|\ \closedb\left(x,r\right)\cap B^{\prime}\ne \emptyset \right\}\right|=1.
\end{equation}
Since $(\alpha\beta)^\ell<\alpha(1-\beta)$, we have that it is true in particular in the case $r\leq\rho_{n\cdot L-1}$. 

Let $\pi : \llbracket \widetilde{T}\rrbracket \to X$ be the continuous projection that maps every $\langle a_0,a_1,...\rangle\in \widetilde{T} $ to the unique element in
$$\bigcap_{n\in \mathbb{N}}B^{\langle a_0,...,a_{n-1}\rangle }\subseteq X.$$
Define
$$Y_L\defeq\pi \left[\llbracket \widetilde{T}\rrbracket\right]=\bigcap_{n\in \mathbb{N}}\bigcup_{B\in A_n}B.$$
Let $\mu = \Prb_U\circ \pi^{-1}$ be a probability measure on $X$, where $\Prb_U$ is the probability measure on $\llbracket \widetilde{T}\rrbracket$ corresponding to the uniform distribution on
$$\left\{q\in \widetilde{T}\ \middle|\ p\prec q\;\text{ and }\;\len\left(q\right)=\len\left(p\right)+1\right\},$$
for every $p\in \widetilde{T}$, as in the proof of Theorem~\ref{unrestricted_upper_bound}. Note that $\mu \left(Y_L\right)=1>0$. Let
$$\delta_L\defeq\frac{L-\ell }{L}\log_{\left(\beta \cdot \alpha \right)^{-1}}\left(m\right),$$
and let $\varepsilon\in \left(0,\delta_L\right)$. Consider the sequence $\left(c_k\right)_{k\in \mathbb{N}}$ defined by
$$c_k=\begin{cases}
    1, &k=0,1,...,L-\ell -1\pmod L,\\
    0, &k=L-\ell,L-\ell+1,...L-1\pmod L.
\end{cases}$$
Let $n_\varepsilon$ be the first positive integer such that
$$ \left|\frac{1}{N}\sum_{k=0}^{N-1}\left(c_k \cdot \log_{\left(\beta \cdot \alpha \right)^{-1}}\left(m\right)\right)-\delta_L\right|<\varepsilon,$$
for every $N> n_{\varepsilon}$. To apply Lemma~\ref{Massdistributionprinciple}, we prove that
$$\mu \left(\closedb\left(y,\rho_N\right)\cap Y_L\right)\le \widetilde{C}\cdot \rho_N^{\delta_L -\varepsilon},$$
for every $y\in Y_L$ and every $N> n_\varepsilon$, where $\widetilde{C}$ is a global constant. Indeed, let $y\in Y_L$ and let $N> n_\varepsilon$. Let $n\in \mathbb{N}$ be the last non-negative integer such that $n\cdot L\le N$. By \ref{choiceofell}, 
$$\left|\left\{B^{\prime}\in A_{n\cdot L-1}\ \middle|\ \closedb\left(y,\rho_N\right)\cap B^{\prime}\ne \emptyset \right\}\right|=1,$$
and since $N-\left(n\cdot L-1\right)\le L$, it holds that
$$\left|\left\{B^{\prime}\in A_N\ \middle|\ \closedb\left(y,\rho_N\right)\cap B^{\prime}\ne \emptyset \right\}\right|\le m^{L-\ell},$$
and in particular $\closedb\left(y,\rho_N\right)\cap Y_L$ is covered by at most $m^{L-\ell}$ closed balls from $A_N$. Furthermore, for every $p\in \widetilde{T}$ such that $\len\left(p\right)=N+1$, we have
\begin{align*}
\Prb_U\left(\llbracket \widetilde{T}_p \rrbracket\right)&=m^{-\sum_{k=0}^{N}c_k}\\
&=\left(\beta \cdot \alpha  \right)^{\sum_{k=0}^{N}c_k\cdot \log_{\left(\beta \cdot \alpha \right)^{-1}}\left(m\right)}\\
&=\left(\left(\beta \cdot \alpha  \right)^{N+1}\right)^{\frac{1}{N+1}\sum_{k=0}^{N}c_k\cdot \log_{\left(\beta \cdot \alpha \right)^{-1}}\left(m\right)}\\
&\le \left(\left(\beta \cdot \alpha  \right)^{N+1}\right)^{\delta_L -\varepsilon}\\
&= \left(\frac{1}{\rho_{-1}}\right)^{\delta_L -\varepsilon}\cdot \left(\rho_{-1}\left(\beta \cdot \alpha \right)^{N+1}\right)^{\delta_L -\varepsilon}\\
&=\left(\frac{1}{\rho_{-1}}\right)^{\delta_L -\varepsilon}\cdot \rho_N^{\delta_L -\varepsilon}\, .
\end{align*}
Hence,
$$\mu \left(\closedb\left(y,\rho_N\right)\cap Y_L\right)\le \widetilde{C}\cdot \rho_N^{\delta_L -\varepsilon}\, ,$$
where $\widetilde{C}=m^{L-\ell}\cdot \left(\frac{1}{\rho_{-1}}\right)^{\delta_L -\varepsilon}$. Thus, as in the proof of Theorem~\ref{unrestricted_upper_bound}, we may apply Lemma~\ref{Massdistributionprinciple} and deduce that
$$\dim_{\mathcal{H}}\left(Y_L\right)\ge \delta_L -\varepsilon,$$
and since $\varepsilon>0$ is arbitrary,
$$\dim_{\mathcal{H}}\left(Y_L\right)\ge \delta_L.$$
Furthermore, $\psi\left[\llbracket s_I\rrbracket\right]\supseteq Y_L$ by the definition of $Y_L$, and thus
$$\dim_{\mathcal{H}}\left(\psi\left[\llbracket s_I\rrbracket\right]\right)\ge \dim_{\mathcal{H}}\left(Y_L\right)\ge \delta_L=\frac{L-\ell }{L}\log_{\left(\beta \cdot \alpha \right)^{-1}}\left(m\right).$$
Since $L>\ell $ is arbitrary,
$$\dim_{\mathcal{H}}\left(\psi\left[\llbracket s_I\rrbracket\right]\right)\ge\log_{\left(\beta \cdot \alpha \right)^{-1}}\left(m\right),$$
completing the proof.
\end{proof}

\end{sloppypar}

\end{document}